\documentclass{article}
\usepackage{amsfonts}
\usepackage{amssymb}

\newtheorem{theorem}{Theorem}[subsection]

\newtheorem{corollary}[theorem]{Corollary}
\newtheorem{definition}[theorem]{Definition}

\newtheorem{lemma}[theorem]{Lemma}

\newtheorem{proposition}[theorem]{Proposition}
\newtheorem{remark}[theorem]{Remark}

\newenvironment{proof}[1][Proof]{\textbf{#1.} }{\ \rule{0.5em}{0.5em}}
\begin{document}
\title{Hierarchical equilibria of
branching populations }
\date{}
\maketitle
\noindent
D.~A.~DAWSON\footnote{  Research supported by NSERC (Canada) and a Max Planck Award for International Cooperation.\\
$^2$ Research supported by CONACYT grant 37130-E (Mexico).\\
$^3$ Research supported by  DFG (SPP 1033) (Germany).}
\\
 {\it Carleton University}\\
{\it Ottawa, Canada K1S 5B6  (e.mail: ddawson@math.carleton.ca)}\\[.5cm]
L.~G.~GOROSTIZA$^2$ \\
{\it Centro de Investigaci\'on y de Estudios Avanzados}\\
{\it 07000 M\'exico D.F., Mexico (e.mail: lgorosti@math.cinvestav.mx)}\\[.5cm]
A.~WAKOLBINGER$^3$\\
 {\it Goethe-Universit\"at}\\
{\it Frankfurt am Main, Germany (e.mail: wakolbin@math.uni-frankfurt.de)}

\begin{abstract}
The objective of this paper is the study of the equilibrium
behavior of a population on the  hierarchical group $\Omega_N$
consisting of families  of individuals undergoing critical
branching random walk and in addition these families also develop
according to a critical branching process. Strong transience of
the random walk guarantees existence of an equilibrium for this
two-level branching system. In the limit $N\to\infty$ (called the
{\em hierarchical mean field limit}), the equilibrium aggregated
populations in a nested sequence of balls $B^{(N)}_\ell$ of
hierarchical radius $\ell$ converge to a backward Markov chain on
$\mathbb{R_+}$. This limiting Markov chain can be explicitly
represented in terms of a cascade of subordinators which in turn
makes possible a description of the genealogy of the population.

\end{abstract}

\vglue1cm \noindent {\bf Mathematics Subject Classifications
(2000):} Primary 60J80;  Secondary 60J60, 60G60.

\vglue1cm \noindent {\bf Key words:} Multilevel branching,
hierarchical mean-field limit, strong transience, genealogy.

\tableofcontents
\newpage

\section{Introduction}\label{introduction}
\setcounter{equation}{0} Spatial branching processes involve two
basic mechanisms, spatial migration and branching. These two
mechanisms  act in opposite directions: the branching causes
fluctuations of the local population densities which are
counteracted by the smoothing effect of the migration, and a {\em
transient} migration is needed to sustain an equilibrium of a
geographically extended population where each individual has an
offspring of mean one. Multilevel branching systems (see, e.g.
\cite {DH}, \cite{GHW}, \cite{Wu}) involve branching, that is
death and replication, at a collective level. For example, in
two-level branching systems  both individuals and {\em } families
(that is, collections of individuals that trace back to a common
ancestor in the individual branching) reproduce independently. In
such systems, the fluctuations of the population densities are
substantially enhanced compared to systems with branching on the
individual level only, and in two-level branching systems a {\em
strongly transient} migration is needed to sustain an equilibrium.
It is well known that Euclidean random walks are transient if and
only if the dimension is bigger than 2, and strongly transient if
and only if the dimension is bigger than 4. In this sense, 2 is
the critical dimension for one-level branching systems, and 4 is
the critical dimension for two-level branching systems.

In the present paper, we will focus on spatial models with a {\em hierarchical} (or {\em ultrametric}) geographical structure ($N$ islands (blocks of radius one) per archipelago (blocks of radius $2$), $N$  blocks of radius $\ell$ per block of radius $\ell+1$, $\ell > 1$, cf. \cite{SF}).  The migration process then is modelled by so called {\em hierarchical random walks}: at a certain rate depending on $\ell$, an individual jumps to a randomly chosen island in distance $\ell$.
 This ultrametric structure leads to a separation of time scales as $N\to \infty$, and makes the models particularly suitable for a thorough analysis of equilibrium states and cluster formation.  It turns out that in the {\em hierarchichal mean field limit} (with order of $N$ individuals per island and $N \to \infty$) there is a separation of time scales in which the population densities in the blocks of different radii evolve. For a block of radius $\ell$, the natural time scale turns out to be $N^\ell$ in the case of one-level branching (see \cite{DG2}) and $N^{\ell/2}$ in the case of two-level branching. On this time scale, the population density in a block of radius $\ell$
 performs, as $N \to \infty$, a  process whose fluctuations are governed by the branching and whose drift is given by a flow of emigration and immigration from the surrounding block.
 For a sequence of nested blocks, this leads to a hierarchy of branching equilibria whose structure we describe in the next subsection. For the case of two-level branching, the convergence of the population densities in nested blocks towards this hierarchy as $N\to \infty$   is stated in Theorem \ref{commute} and proved in section \ref{sec5}.

Generically, in our hierarchical model  the migration process that sustains an equilibrium is at the border to recurrence in the case of one-level branching, and  at the border to weak transience in the case of two-level branching, as $N\to \infty$. In this sense,  the hierarchical one-level branching equilibria studied in \cite{DG2} correspond to a situation ``near dimension 2'', and the hierarchical two-level branching equilibria studied in our paper correspond to a situation ``near dimension 4''.
Dimension 4 is of
considerable interest because it serves as a critical
dimension not only for  the two-level branching systems studied
in this paper but also for a number of phenomena in statistical
physics including the large scale fluctuations of ferromagnetic
models at the critical temperature.

The structures of the family clusters in equilibrium can be best understood in terms of the genealogy of the branching system, see \cite{DPhist, DP, SW}). We will describe the genealogy of the
equilibrium population in the mean-field limit using a cascade of
subordinators.

\section{Overview}
\subsection{Hierarchies of one--and two--level Feller branching diffusions}\label{HEBP}
\setcounter{equation}{0}
Consider a large population whose size is fluctuating because of
critical reproduction, and which is subject to emigration of
individuals and immigration from a surrounding (still larger)
reservoir of individuals. The immigration rate is given by the
population density in the environment, which
fluctuates on a slower time
scale. Now consider an infinite hierarchy of such populations
$\mathcal P_{\ell}, \, \ell=1,2,..$, where $\mathcal P_{\ell+1}$ acts as an
environment for  $\mathcal P_{\ell}$, and think of an equilibrium
situation. We will study two situations where there is
a sequence of time scales such that, in the limit of large local population sizes,
on the time scale associated with
$\mathcal P_{\ell}$
the population density $\zeta_{\ell+1}$ of
           $\mathcal P_{\ell+1}$ remains
constant, and given $\zeta_{\ell+1}=a$,
the dynamics of the population density $\zeta_{\ell}$ of $\mathcal
P_{\ell}$
is of the form
\begin{equation}\label{fluc}
                      d\zeta_{\ell}({t}) =
dM_{\ell}{(t)}-c_{\ell}(\zeta_{\ell}({t}) -
                      a)\,dt.
\end{equation}
Here $c_{\ell}$ is a positive constant which describes the
migration intensity into and out of $\mathcal P_{\ell}$, and
$M_{\ell}$ is a martingale describing the fluctuations of
$\zeta_{\ell}$.

In subsection \ref{review} we will describe a situation in which the
martingale $M_{\ell}$ has quadratic variation
\begin{equation}\label{FQV}
      d\langle M_{\ell} \rangle{(t)} = \zeta_{\ell}({t}) dt,
      \end{equation}
      hence in this case (\ref{fluc})  specializes to
      \begin{equation}\label{subcF}
                      d\zeta_{\ell}({t}) =
\sqrt{\zeta_\ell(t)} dW_\ell(t)-c_{\ell}(\zeta_{\ell}({t}) -
                      a)\,dt,
\end{equation}
where $W_1, W_2,...$ are independent Wiener processes. For each $\ell$ (\ref{subcF})
      is the stochastic differential equation of a
      subcritical Feller branching diffusion with immigration (\cite{AN, EK}).

                  Later on, we will consider a dynamics of the population
                  density $\zeta_{\ell}$ which is again of the form
                  (\ref{fluc}) but where the fluctuations are governed by a
                  ``family structure'' of the population. More precisely, the
                  martingale $M_{\ell}$ has quadratic variation
                  \begin{equation}\label{SFQV}
                      d\langle M_{\ell} \rangle{(t)} =
\left(\int_{(0,\infty)} x^2
                      \xi_{\ell}(t,x)
                  dx\right) dt.
                  \end{equation}
                 where $\xi_{\ell}(t,x)$
                  measures the rescaled number of families of size $x$
within $\mathcal
P_{\ell}$.
The link between the population density
$\zeta_{\ell}{(t)}$ and the density
$\xi_{\ell}(t,x)$ of family sizes is given by
\begin{equation}\label{aggreg}
                     \zeta_{\ell}{(t)} = \int_{(0,\infty)} x\,
\xi_{\ell}(t,x) dx.
                     \end{equation}
The form (\ref{SFQV}) of the fluctuations of $\zeta_{\ell}$
indicates that we are dealing with  Feller branching diffusions of
{\em families} rather than individuals. This {\em family
branching} shows up in the dynamics which is described by an
absolutely continuous measure-valued process with density
$\xi_{\ell}(t,x)$ satisfying the stochastic partial
differential equation (SPDE)
\begin{equation} \label{SPDEca2}
                     \frac{\partial}{\partial {t}} \xi_{\ell}(t,x) =
                     \sqrt{\xi_{\ell}(t,x)}\,\dot{W}_{t}(x)
                     +  G^\ast \xi_{\ell}(t,x)-c \left(-\frac {\partial}{\partial x}x \xi_\ell(t,x)+a\delta'_{0}(x)\right)
\end{equation}
with $c=c_\ell$ where $G^\ast$ is the adjoint of the generator $G$ of a critical Feller branching diffusion  given by
\begin{equation}\label{G} Gf=\frac{1}{2}x\frac{\partial^2}{\partial
x^2}f,\end{equation}
$\dot{W}$ is space-time white noise and $\delta_0 '$ is the
derivative (in the sense of Schwartz distributions) of the
$\delta$-function at $0$. An equivalent formulation of (\ref{SPDEca2}) is
\begin{equation} \label{SPDEca2test}
\frac{\partial}{\partial {t}} \langle \xi_{\ell}(t),f\rangle = \langle \sqrt {\xi_\ell(t)}\dot W_t, f \rangle + \frac 12
\langle \xi_\ell(t), xf'' \rangle - c\left(\langle \xi_\ell(t), xf'\rangle - a\lim_{\varepsilon \downarrow 0} \langle \frac 1\varepsilon \delta_\varepsilon, f\rangle\right),
\end{equation}
where $\delta_\varepsilon$ is the Dirac measure in $\varepsilon$, and  $f:(0,\infty)\to \mathbb R$ has bounded first and second derivatives with $f(x) \le \mbox{const } x$.
Note that the first term in (\ref{SPDEca2}) (and (\ref{SPDEca2test})) comes from the familiy branching, the second comes from the individual branching, and the ``mean reversion''
term comes from the individual emigration and the immigration of infinitesimally small families.

In addition, (\ref{SPDEca2}) shows that the sizes of
the single families develop independently according to {\em
individual} subcritical Feller branching diffusions. We will
therefore call $\xi_{\ell}$, $\ell=1,2,\ldots$, a hierarchy
of {\em two-level branching diffusions}.

Two-level branching diffusions have been introduced by Dawson and
Hoch\-berg \cite{DH} as superprocesses over Feller branching
diffusions, where all mass accumulating in $0$ is removed.
Therefore, these processes take their values in the measures on
$(0,\infty)$. In fact, it turns out  that for $t>0$ they have
absolutely continuous states. In our context, in addition to the
set-up of \cite{DH}, there is a ``continuous immigration of small
families''. We will see how this fits into the framework of
immigration processes from the boundary studied in Li and Shiga
\cite{LS}. For general background on superprocesses and related
stochastic partial differential equations, see \cite{D} and
\cite{DP}.

The hierarchies of branching equilibria considered in our paper
are motivated through a spatial picture which we describe for the
case of (\ref{FQV}) (``one-level branching'') in subsection
\ref{review} and for the case of (\ref{SFQV}) (``two-level
branching'') in subsection \ref{preview}. The case of a hierarchy
of one-level branching systems was studied by Dawson and Greven
\cite{DG1,DG2} in the context of super-random walks (or
interacting Feller branching diffusions).

For any given  $\theta > 0$ (which in the geographical model will play the role of a ``global population density'') we will construct the hierarchy
\begin{equation}\label{entrancel}
                     (\ldots,
\zeta^{\theta}_{\ell+1}, \zeta^{\theta}_{\ell},\ldots,
\zeta^{\theta}_{2}, \zeta^{\theta}_{1}),
\end{equation}
in terms of an entrance law for a backward Markov chain where the
conditional law of $\zeta_{\ell}$
      given
$\zeta_{\ell+1}=a$  is an equilibrium state of
(\ref{fluc}). More precisely, in subsection \ref{sec4.1} we will show the following result.

\begin{proposition}\label{entrancela}
Let $(c_\ell)_{\ell \ge 1}$ be a sequence of positive numbers, and let us
distinguish two cases:

a) Assume $\sum_\ell c_\ell^{-1} < \infty$. For $a >0$, $\ell \in \mathbb N$, let
$\Pi_\ell^{(1)}(a,.)$ be the equilibrium distribution of (\ref{subcF}).

b) Assume $\sum_\ell c_\ell^{-2} < \infty$. For $a >0$, $\ell \in \mathbb N$ and $c = c_\ell$,  let $\xi_{\ell,a}$ be an equilibrium state of (\ref{SPDEca2}), and $\Pi_\ell^{(2)}(a,.)$ be the distribution of $\int_{(0,\infty)} x \xi_{\ell,a}(x)dx$.

In both cases,  for every $\theta > 0$ the backward Markov chain with transition probability
$$ P(\zeta_\ell \in A|\zeta_{\ell+1}=a)= \Pi_\ell(a,A),$$
where $\Pi$ is either $\Pi^{(1)}$ or $\Pi^{(2)}$,
has a unique \textit{entrance law}
$\{\zeta_{\ell}^{\theta}\}_{\ell =\dots,2,1}$ satisfying
\begin{equation}\label{constexp} \mathbb E\zeta_{j}^{\theta}\equiv \theta \end{equation}
and
\begin{equation}\label{entrancelimit} \lim_{j\to\infty}\zeta_{j}^{\theta}=\theta \quad a.s.\end{equation}
\end{proposition}

\subsubsection{A cascade of subordinators}\label{Cassub}
\setcounter{equation}{0}

To work out parallels between the one- and two-level branching
situations described in subsections \ref{review} and
\ref{preview}, and to discuss aspects relevant for the genealogy
of the hierarchical branching equilibria, we write
$\Pi_{\ell}(a,.)$ for the equilibrium distribution of (\ref{fluc})
in the two cases (\ref{FQV}) and (\ref{SFQV}) (which correspond to
cases a) and b) in Proposition \ref{entrancela}).

In both cases the parameter $a$ enters as a factor into the
immigration rate of a continuous state branching process, hence
$\Pi_{\ell}(a_1+a_2,.)=\Pi_{\ell}(a_1,.)*\Pi_{\ell}(a_2,.)$.
Therefore the $\Pi_{\ell}(a,.)$  are infinitely divisible
distributions on $(0,\infty)$ and there exist subordinators (that
is, processes with stationary independent non-negative increments)
$S_{\ell}(a)$ , $a \ge 0$, such that
\begin{equation}\label{subrepn}
        \mathcal L
(S_{\ell}(a)) = \Pi_{\ell}(a,.).
\end{equation}
          We will see in subsection \ref{SFB} that in case a) the
$S_{\ell}$ are Gamma processes.

In both cases, the L\'evy-Khinchin decomposition (see \cite{K2},
Chapt. 15) $S_{\ell}(a) = \sum_{i} H_{i}$ describes
(asymptotically as $N \to \infty$) the partitioning of the
equilibrium population in $B^{(N)}_{\ell}$ into subpopulations
stemming from one and the same immigrant into $B^{(N)}_{\ell}$,
given that the population density of the surrounding block
$B^{(N)}_{\ell+1}$ is $a$.

Recall from (\ref{fluc}) and Proposition \ref{entrancela}  that
\begin{equation}
             \int x \Pi_{\ell}(a,dx) = a.
\end{equation}
Therefore,
\begin{equation}
             \mathbb E S_{\ell}(a) =   a.
\end{equation}

Let us denote (in either of the two cases) the L\'evy measure of
$S_{\ell}$ by $\mu_{\ell}$, and the second moment of $\mu_{\ell}$
by $m_{\ell}$. An iteration of the L\'evy-Khinchin representation
(which can be interpreted in terms of the genealogy of the
branching hierarchy, see subsection \ref{sec4.4}) will show that
the L\'evy measure of the iterated subordinator
$S_{\ell}(S_{\ell+1}(\ldots S_{j-1}))$ has second moment
$m_{\ell}+\ldots+m_{j-1}$ (see subsection \ref{subs4.2}). Under
the condition
\begin{equation}\label{mfinite}
\sum_{\ell=1}^\infty m_{\ell} < \infty.
\end{equation}
we will prove in subsection \ref{sec4.3} that for each $\theta >0$
the
 limit in distribution
\begin{equation}\label{zetalimit}
                     \zeta^\theta_{\ell}=d\mbox{-}\lim_{j\to
\infty}S_{\ell}(S_{\ell+1}\ldots S_{j-1}(\theta))),
\end{equation}
exists, has expectation $\theta$ and defines an entrance law with
respect to $(\Pi_{\ell}(a,.))$. In particular one has
\begin{equation}\label{entrancelaw}
\zeta^\theta_{\ell}=^d S_{\ell}(S_{\ell+1}\ldots
S_{j-1}(\zeta^\theta_{j}))),  \quad j > \ell.
\end{equation}
For each $j > \ell$ this gives a decomposition of
$\zeta^\theta_{\ell}$, which asymptotically as $N\to\infty$,
stands for the partitioning of the equilibrium population
$\zeta^{(N,\theta)}_\ell$ in $B^{(N)}_{\ell}$ into subpopulations
stemming from one and the same immigrant into $B^{(N)}_{j-1}$.

The summability condition (\ref{mfinite}) is equivalent to the
transience condition (\ref{trans}) in the one-level case, and to
the strong transience condition (\ref{trans2}) in the two-level
case, since we will show that $m_{\ell}= 1/2c_{\ell}$ in the
one-level case, and
 $m_{\ell}= 1/4c_{\ell}^2$ in
the two-level case (see Remarks  \ref{momentsofgamma} and
\ref{momentsnu}).

\subsubsection{Genealogy}
In Section \ref{sec4.4} we develop a genealogy of the jumps
occurring in the cascade of subordinators. The idea is that given
a jump of $S_{\ell+1}(\cdot)$ at time $t_i$ there will be a set of
jumps of $S_\ell(\cdot)$ that occur in the time interval
$(S_{\ell+1}(t_i-),S_{\ell+1}(t_i))$ and these level $\ell$ jumps
will be said to be {\it descendants} of the level $(\ell+1)$ jump.
 In subsections (\ref{sec4.4}), (\ref{sec4.5}) and
(\ref{sampledind})   we use this idea  to develop the full
genealogical structure of the population associated with the
entrance law. This leads to a decomposition of the population into
a countable collection of subpopulations of individuals having a
common ancestor. For the  case of critical individual branching
this was done in \cite{DG2}, for the two-level branching case this
is new. We will work out the parallels between the two cases in a
general framework, which also sheds some new light on the results
of \cite{DG2}.

 Intuitively this genealogy  describes the limiting
genealogical structure  of the spatial branching equilibria with
hierarchical geographic structure described in the Introduction as
the parameter $N\to\infty$ and the analogue of the ``clan
decomposition'' of the equilibrium of super-Brownian motion (e.g.
\cite{DPhist}).

\subsection{Hierarchical geography, random walks and branching equilibria}
\setcounter{equation}{0}

\subsubsection{A class of random walks}\label{classofrw}
In order to give a precise formulation for the spatial system we
now  describe the set $\Omega_N$ of sites on which the spatial
population lives. For fixed $N\in \mathbb N$, let $\Omega_N$ be
the countably many leaves of a tree all of whose inner nodes have
degree $N+1$. In other words, each node in depth $\ell +1$, $\ell
= 0,1,\ldots$, has $N$ successors in depth $\ell$ and one
predecessor in depth $\ell +2$. For two sites $y,\,z \in
\Omega_N$, their hierarchical distance $d_{N}(y,z)$ is defined as
the depth of their closest common ancestor in the tree. Note that
$d_N$ is an ultrametric, that is, $d_N(y,z)\leq
\max\{d_N(y,x),d_N(x,z)\}$. We define the individual migration on
$\Omega_N$ as follows. Let $q^{(N)}_{1}, q^{(N)}_{2},\ldots$ be
positive numbers with $\sum_{\ell}q^{(N)}_{\ell}<\infty$. At rate
$q^{(N)}_{\ell}$ the individual makes a jump of distance $\ell$
(i.e., it waits for an exponentially distributed time with
parameter $\sum_{\ell}q_{\ell}^{(N)}$ and then jumps a distance
$j$ with probability $q_{j}^{(N)}/\sum_{\ell}q_{\ell}^{(N)}$),
choosing its arrival site uniformly among all sites at a distance
$\ell$ from its previous site.

The set $\Omega_N$ can be identified with the set  of sequences in
$\{0,1,\ldots,N-1\}$ almost all of whose entries are equal to
zero. With component-wise addition mod $N$, $\Omega_{N}$ is a
countable Abelian group (the so called {\em hierarchical group of
order} $N$). Note that $d_N(y,z)$ is translation invariant; it
will be written as $|y-x|$. The migration specified above is a
(continuous time) random walk on $\Omega_{N}$ called \textit{
hierarchical random walk}. Hierarchical random walks, in
particular their transience and recurrence properties, are studied
in \cite{DGW2}, \cite{DGW3}.

\subsubsection{A system of branching random walks}
We now introduce a system of branching random walks on $\Omega_N$.
This is given by a system of particles undergoing symmetric random walks
with migration rates $q^{(N)}_\ell$ together with branching.
Branching and migration are assumed independent. We specify
the branching mechanism as simply as possible: after a unit
exponential time, an individual, independently of all the others,
either is removed or splits into two, each case with probability $1/2$.

\begin{remark}\label{onelevelequ} \cite{LMW}, \cite{Gr}
Assume that the migration rates $q^{(N)}_{\ell}$ are such that the
individual random walk is transient. Then, for each $\theta > 0$,
there exists a unique branching equilibrium with a mean number of
$\theta$ individuals per site. This equilibrium is of Poisson
type, i.e. the equilibrium population $\Phi$ is the superposition
$\Phi = \sum_{i}\Phi_{i}$ of a Poisson system
$\sum_{i}\delta_{\Phi_{i}}$ of families (each family consists of a
collection of individuals which in a historical picture can be
traced back to a common ancestor).
\end{remark}

\subsubsection {The hierarchical mean-field limit of a
branching equilibrium ``near dimension two''}\label{review}
Now consider, for a large $N$, the total number
$Z^{(N)}_{\ell}(t)$ of individuals in a closed ball
$B_{\ell}^{(N)}$ of radius $\ell$ at time $t$. Note that
$B_{\ell}^{(N)}$ has $N^\ell$ elements, and look at the time
evolution of the {\em population density} (or {\em block mean})
$A^{(N)}_{\ell}(t) = Z^{(N)}_{\ell}(t)N^{-\ell}$ in
$B_{\ell}^{(N)}$ at time scale $t N^{\ell}$. This corresponds to
the classical Feller branching diffusion approximation \cite{EK}.
However, in order for the immigration into and emigration from
$B_{\ell}^{(N)}$ to produce a nondegenerate drift term in the
limit $N\to \infty$ one must adjust the migration rates. The
appropriate  adjustment is
\begin{equation}\label{migrates}
                     q^{(N)}_{\ell} = c_{\ell-1}N^{-{(\ell-1)}},
\end{equation}
where the $c_{\ell}$ do not depend on $N$ and satisfy
\begin{equation}\label{cbound}
 \limsup_{\ell} \frac {c_{\ell+1}}{c_{\ell}} < \infty.
 \end{equation}
  We will call a random
walk on $\Omega_{N}$ with jump rates $q^{(N)}_{\ell}$ of the form
(\ref{migrates}) a $(1,(c_{\ell}),N)$-{\em random walk}.  The
following proposition is proved in \cite{DGW2}.

\begin{proposition}\label{transhier}
Consider a $(1,(c_{\ell}),N)$-random walk on $\Omega_N$ with rates
given by (\ref{migrates}). Assume that
\begin{equation}\label{cboundN}
                     \limsup_{\ell} \frac {c_{\ell+1}}{c_{\ell}} < N.
                     \end{equation}
                  Then the
random walk is transient if and only if
\begin{equation}\label{trans}
                     \sum_{\ell}c_{\ell}^{-1}<\infty.
                     \end{equation}
\end{proposition}

Now assume that $(c_j)$ satisfies (\ref{cbound}).
Then for each $N >\limsup_{\ell} \frac {c_{\ell+1}}{c_{\ell}}$,
because of Proposition \ref{transhier} and Remark
\ref{onelevelequ} there is an equilibrium  for the system of $(1,(c_{\ell}),
N)$-branching random walks with mean $\theta$ for each $\theta
>   0$ . We now
consider  the corresponding equilibrium population densities
$A_{\ell}^{(N)}$  in $B_{\ell}^{(N)}$, $\ell = 1,2,\dots$ (where
here and below $(B_{\ell}^{(N)})$ denotes a sequence of nested
balls of radii $\ell$ in $\Omega_{N}$). In order to identify the
limit as $N\to\infty$ of this sequence of population densities we
must consider the dynamics of $A_{\ell}^{(N)}$ in its natural time
scale $N^\ell$. Let us first discuss on a heuristic level why, in
the limit $N\to \infty$, the drift term on the r.h.s. of
(\ref{fluc}) arises on the time scale $N^{\ell}$ for the
population density in $B_{\ell}^{(N)}$. Because of the ultrametric
structure of $\Omega_{N}$, an individual in $B_{\ell}^{(N)}$ has
to make a jump of size $\ge \ell +1$ in order to leave
$B_{\ell}^{(N)}$. Because of (\ref{migrates}) and (\ref{cbound}),
for $N$ large, jumps of size $> \ell+1$ happen rarely compared to
jumps of size $\ell+1$ (since $q_{l+k}=o(q_{l+1})$ as $N\to
\infty$ for $k > 1$). Hence the individual emigration rate from
$B_{\ell}^{(N)}$ on time scale $N^{\ell}$ is $c_{\ell}$
(asymptotically as $N \to \infty$). Concerning immigration into
$B_{\ell}^{(N)}$, again because of (\ref{migrates}), it is only
the environment in $B_{\ell+1}^{(N)}$ that is relevant for large
$N$. An individual in $B_{\ell+1}^{(N)}\setminus B_{\ell}^{(N)}$
necessarily has to jump a distance $\ell+1$ in order to make it
into $B_{\ell}^{(N)}$, and on average every $(N-1)$-st of these
jumps will take the individual into $B_{\ell}^{(N)}$ (note that
$B_{\ell+1}^{(N)}$ is $N$ times as large as $B_{\ell}^{(N)}$).
Since the block mean $A^{(N)}_{\ell+1}$ does not change
on the time scale $N^{\ell}$ as $N \to \infty$, the total immigration rate into
$B_{\ell}^{(N)}$ on this time scale is (asymptotically as $N\to
\infty$) of the order
\begin{equation}\label{immrates}
                     Z^{(N)}_{\ell+1}N^{\ell}q^{(N)}_{\ell+1}/N =
(N^{\ell+1}A^{(N)}_{\ell+1})c_{\ell}/N =
c_{\ell}A^{(N)}_{\ell+1}N^{\ell}.
\end{equation}
This suggests that the limiting dynamics of the population
densities $A^{(N)}_\ell$, in their natural time scales, as
$N\to\infty$ is given by (\ref{fluc}) with $a=A^{(N)}_{\ell+1}$.
The separation of time scales on balls of different radii that
underlies the previous discussion is a feature of the hierarchical
random walks, which is due to the ultrametric structure of
$\Omega_N$ (see \cite{DGW2}). This is also explained in more
detail in Remark \ref{timescale} below.

Instead of branching particle systems, Dawson and Greven
\cite{DG2} consider super-random walks (or so-called interacting
Feller diffusions) on $\Omega_{N}$. Note also that the definition
of the random walk in \cite{DG2} is slightly different but
asymptotically equivalent as $N\to\infty$ to the one used in this
paper. In \cite{DG2} the sites to which a jump is made are chosen
with uniform distribution on a ball rather than on a sphere.
However, the  ``interior'' of the ball is asymptotically
negligible as compared to the sphere as $N$ goes to
infinity.

A particle system analogue of Theorem 4(b) of \cite{DG2}  is the
following, which we state without proof.

\begin{proposition}\label{convlevel1} Consider a sequence $(c_{\ell})$
      satisfiying  conditions (\ref{trans}) and (\ref{cbound}) for transience.
For $N$ large enough such that (\ref{cboundN}) is met, let the
individual migration process be a $(1,(c_{\ell}),N)$-random walk. For $\theta > 0$
let $Z_{\ell}^{(N, \theta)}(t)$ denote the total number of individuals in
$B^{(N)}_\ell$ at time $t$ and  $A_{\ell}^{(N, \theta)}(t):=
Z_{\ell}^{(N, \theta)}(t)N^{-\ell}$ be the population density at time $t$
in $B_{\ell}^{(N)}$ in the Poisson type branching equilibrium
population on $\Omega_{N}$ with mean number $\theta$ of
individuals per site (see Remark \ref{onelevelequ}). Let $\{\zeta_\ell^\theta\}$
be the entrance law provided by Proposition \ref{entrancela}, case a).
Then
\[\{A_{\ell}^{(N,\theta)}(0)\}_{\ell\in \mathbb{N}}
\Longrightarrow \{\zeta_{\ell}^{\theta}\}_{\ell
\in\mathbb{N}}\quad \mbox { as }\;N\to\infty ,\]
where $\Longrightarrow$ denotes weak convergence of finite dimensional distributions.
\end{proposition}

Let us now explain in which sense transient $(1,(c_{\ell}),N)$-random walks
can be interpreted as random walks ``near dimension $2$''.

\begin{definition}  \cite{DGW2} Let $Z$ be an irreducible transient random
walk on a countable Abelian group $\Gamma$. Its degree of
transience is defined by
\begin{equation}\label{orderoftrans}
                     \gamma:=\sup \{\eta \ge 0 : \mathbb E_{0}L^\eta < \infty\},
\end{equation}
where $L$ is the last exit time of $Z$ from $0 \in \Gamma$.
\end{definition}

Expressed in more analytic terms, the degree of transience of $Z$
is
$$\gamma=\sup \{\eta \ge 0 :  \int_{1}^\infty t^{\eta} p_{t}(0,0)\, dt <
\infty\},$$ where $p_{t}$ is the transition probability of $Z$
\cite{DGW2, SaW}.

For simple symmetric random walk on the $d$-dimensional lattice
$\mathbb Z^d$, it is well-known that dimension 2 is the borderline
for transience. For $d > 2$, the degree of transience is $d/2-1$,
since the rate of decay of the transition probability is
$p_t(0,0)\sim \textrm{const.} t^{-d/2}$.

\begin{remark} \cite{DGW2} \label{borderstrongrec}
a) Let $0 < c < N$. Then the
              $(1,(c^{\ell}),N)$-random walk on $\Omega_{N}$ is transient if and only if $c>1$.
         In this case its degree of transience is $\log
c/(\log N-\log c)$. Thus for fixed $c$ the transient
$(1,(c^{\ell}),N)$-random walks on $\Omega_{N}$ have degrees of
transience $O(1/\log N)$ and therefore asymptotically as $N\to
\infty$, can be viewed as analogues of  random walks ``near (Euclidean)
dimension $2$''.

b) Assume that $(c_\ell)$ satisfies conditions $(\ref{trans})$ and $(\ref{cbound})$, and put $c := \limsup c_{\ell+1}/ c_{\ell}\ge 1$. Then, for all $N>c$ the
$(1,(c_{\ell}),N)$-random walk on $\Omega_{N}$ is transient with degree of transience in the interval $ [0,
\log
c/(\log N-\log c)]$.
\end{remark}

Since certain properties of systems of branching random walks such
as persistence and structure of occupation time fluctuations
depend only on the degree of transience of the random walks,
branching populations whose migration process is a hierarchical
random walk can give insight into the behavior of a larger class
of branching populations whose random walks have the same degree
of transience.

\subsection{Two-level branching systems in a hierarchical geography}

\subsubsection{Strongly transient migration}\label{STRW}

Whereas the situation described in  subsection \ref{review} gives
an analogue to a  situation ``near dimension 2'', our main focus
later on will be on the analogue to a situation ``near dimension
4''.  In this context we will consider the (stronger) mass
fluctuations induced by a critical reproduction of whole {\em
families} (of mutually related individuals), together with a
(stronger) smoothing caused by a strongly transient migration.
\begin{definition}
                     An irreducible random walk $Z$ on a countable
Abelian group $\Gamma$ is called strongly transient if
\begin{equation}\label{st}
\mathbb E_{0} L < \infty
\end{equation}
where $L$ denotes the last exit time of $Z$ from $0$. A transient
random walk with $\mathbb E_{0} L = \infty $ is called weakly
transient
\end{definition}

Note that strong transience is equivalent to
$$\int_{1}^{\infty}tp_{t}(0,0)dt<\infty$$
and a necessary condition is that the  degree of transience be
equal to or greater than 1  \cite{DGW2}. Moreover, as mentioned
above simple symmetric $d$-dimensional random walk has degree of
transience $d/2-1$, and it is strongly transient iff $d>4$.

In order to introduce a family of strongly transient random walks
on $\Omega_N$ we replace (\ref{migrates}) by

\begin{equation}\label{migrates2}
                     q^{(N)}_{\ell} = c_{\ell-1}N^{-{(\ell-1)/2}},
\end{equation}
where the $c_\ell$ do not depend on $\ell$ and satisfy
(\ref{cbound}).
      We will call a random walk on $\Omega_{N}$ with jump rates
$q^{(N)}_{\ell}$ of the form (\ref{migrates2}) a
$(2,(c_{\ell}),N)${\em -random walk}.

\begin{proposition} \label{twolevelc} \cite{DGW2}
Consider a $(2,(c_{\ell}),N)$-random walk on $\Omega_N$ with rates
given by (\ref{migrates2}).  Assume that
\begin{equation}\label{cbound2N}
                     \limsup_{\ell} \frac {c_{\ell+1}}{c_{\ell}} < N^{1/2}.
                     \end{equation}
Then the random walk is strongly transient if and only if
\begin{equation}\label{trans2}
                     \sum_{\ell}c_{\ell}^{-2}<\infty.
                     \end{equation}

\end{proposition}

\begin{remark} \cite{DGW2}
a) Let $0 < c < N^{1/2}$. Then the
              $(2,(c^{\ell}),N)$-random walk on $\Omega_{N}$ is strongly transient if and only if $c>1$.
              In this case its degree of transience is
$$\frac{\log N +2\log c}{\log N-2\log c}.$$ Thus for fixed $c$  the
strongly transient $(2,(c^\ell),N)$-random walks have degree of
transience $1+O(1/\log N)$ and therefore  asymptotically as $N\to
\infty$ can be viewed as analogues of random walks ``near
(Euclidean) dimension 4''.

b) Assume that $(c_\ell)$ satisfies conditions $(\ref{trans2})$
and $(\ref{cbound})$, and put $c := \limsup c_{\ell+1}/ c_{\ell}
\ge 1$. Then, for all $N>c^2$ the $(2,(c_{\ell}),N)$-random walk
on $\Omega_{N}$ is strongly transient with degree of transience in
the interval $[1, (\log N +2\log c)/(\log N-2\log c)]$.
\end{remark}

\begin{remark}\label{timescale}
The natural time scale for the strongly transient
$(2,(c^\ell),N)$-random walk in $B_\ell^{(N)}$ is
$N^{{\ell}/{2}}$. More precisely, as shown in  \cite{DGW2},
asymptotically as $N \to \infty$ for this random walk on the time
scale $N^{\ell/2}$ only the migrations within the ball
$B_\ell^{(N)}$ and the surrounding ball $B_{\ell+1}^{(N)}$ are
relevant.
 Similarly, this occurs for the transient $(1,(c^\ell),N)$-random walk in time scale $N^{\ell}$
 (see \cite{DGW2} for details). This effect is basic for the limiting hierarchy of branching
 equilibria obtained in this paper.
\end{remark}

\subsubsection{Two-level branching equilibria}\label{tlbe}

A main objective of this paper is to study two-level branching
systems for a migration which is on the border between strong and weak transience
-- recall that a strongly transient migration is required for the existence of a branching equilibrium.
Thus, for Euclidean random walks,  $d=4$  is the critical dimension for
a two-level branching system in the same
  way as dimension $d=2$ is the critical dimension for a one-
level branching system.

We are going to study two-level branching systems on $\Omega_N$.
Consider a system of $(2, (c_{\ell}), N)$-random walks on $\Omega_{N}$ such
that $(c_{\ell})$ satisfies the conditions (\ref{cbound2N})
  and (\ref{trans2}) for strong transience, and recall from Remark \ref{borderstrongrec}  that
  these random walks are close to the border to weak transience at least for large $N$ if $(c_\ell)$ does not grow superexponentially.
 Introduce, in addition to the individual branching and migration,
a {\em family branching}: independently of everything else, after
a unit exponential time each family $\Phi_i$  either vanishes or
reproduces resulting in two identical copies
$\Phi_i^\prime$,$\Phi_i^{\prime\prime}$,  each case with probability $1/2$. After a reproduction event
$\Phi_i^\prime$ and $\Phi_i^{\prime\prime}$ evolve as independent one
level branching systems. This creates the basic {\em two-level
branching system}  $\Psi^{(N)}(t)$ which is
started with the family system at time $t_0$ given by
$\Psi^{(N)}(t_0) =\sum_{i}\delta_{\Phi_{i}}$ described in Remark
\ref{onelevelequ}.

 The
following result is the analogue of \cite{GHW} for two-level
branching systems on $\Omega_N$.

\begin{proposition} \label{4+dequil}  Assume that\newline
(i) the random walk on $\Omega_N$ is strongly transient and
\newline (ii) $\Psi^{(N)}(t_0) =\sum_{i}\delta_{\Phi_{i}}$ where
$\{\Phi_i\}$ corresponds to the family decomposition of an
equilibrium state for the one-level branching random walk with
mean number $\theta$ of individuals per site.\newline Then as
$t_0\to -\infty$, the two-level branching
      system $\Psi^{(N)}(0)$ converges in distribution
      to a translation invariant  equilibrium   $\Psi^{(N,\theta)}(0)$ with a
mean number $\theta$ of individuals per site.
\end{proposition}
\begin{remark} The notation $\Psi^{(N)}(t)$ and $\Psi^{(N,\theta)}(t)$
will be used throughout to denote the two level branching system
and the equilibrium process with mean $\theta$, respectively.

\end{remark}
\begin{remark} Greven and Hochberg \cite{GH} have obtained more general
conditions under which the convergence to equilibrium as in
Proposition \ref{4+dequil} occurs as well as conditions under
which it fails.
\end{remark}
\medskip

       Now we consider a system  of
$(2,(c_\ell),N)$ random walks with $(c_{\ell})$ satisfying
conditions (\ref{trans2}) and (\ref{cbound}) for strong
transience. Then for each $N>\left( \limsup_\ell
\frac{c_{\ell+1}}{c_\ell}\right)^2$ because of Propositions
\ref{twolevelc} and \ref{4+dequil} there is a two-level branching
equilibrium $\Psi^{(N, \theta)}(t)$ with a mean number of $\theta$ individuals per site
for each $\theta > 0$. A main objective of this paper is to study
the equilibrium structure that arises in the limit as
$N\to\infty$ of the corresponding sequence of family structures in
the blocks $B^{(N)}_\ell$.

\subsection{The hierarchical mean-field limit of a two-level
branching equilibrium ``near dimension four''}\label{preview}

\subsubsection{Local normalized family-size process} \label{LNFSP}
Let, for fixed $N$,
$\Psi^{(N,\theta)}(t), t\in \mathbb{R},$ be the equilibrium process of the
two-level branching system  as provided by Proposition
\ref{4+dequil}. Denote the number of families in
$\Psi^{(N,\theta)}(t)$ having $j$ individuals in
$B_{\ell}^{(N)}$ by $n_{\ell}^{(N,\theta)}(t,j)$ and write
$H_{\ell}^{(N,\theta)}(t):=\sum_{j}n_{\ell}^{(N,\theta)}(t,j)\delta_{j}$ for the {\em local family-size process}.
For each fixed $\ell$, $H_{\ell}^{(N,\theta)}(0)$
is a random measure on $\mathbb{N}$ which describes the population
of equilibrium family sizes in the block $B_{\ell}^{(N)}$. We note
that the  process $\{H_\ell^{(N,\theta)}(t)\}_{t\in \mathbb{R}}$ can be
viewed as a branching Markov chain on $\mathbb{Z}_+$ with
instantaneous killing at $0$, the Markov chain on $\mathbb{Z}_+$
being a standard subcritical binary branching process with
immigration. Note that $\sum_{j}jn_{\ell}^{(N,\theta)}(t,j)=Z^{(N, \theta)}_{\ell}(t)$,
the number of individuals in $B^{(N)}_\ell$ at time $t$ . Now
consider the \textit{equilibrium normalized family size process}
defined by
\begin{equation} \label{norm} \eta_{\ell}^{(N,\theta)}(t)=
\sum_{j}N^{-\ell/2}n_{\ell}^{(N,\theta)}(tN^{\ell/2},j)\delta_{jN^{-\ell/2}},
\end{equation}
In other words, for $0<a<b<\infty$,
\begin{equation}\label{norm2}
   \eta_{\ell}^{(N,\theta)}(t)(a,b) = N^{-\ell/2}\eta_{\ell}
(N^{\ell/2}t,N^{\ell/2}(a,b)).
\end{equation}
Note that the natural time scale in which to observe the
subpopulation in $B^{(N)}_\ell$ in this case is $N^{\ell/2}$ and
not $N^\ell$ as was the case for one-level branching (see Remark \ref{timescale}).

For each $\ell, N$ and $t$, $\eta_{\ell}^{(N, \theta)}(t)$ is a
random measure on $(0,\infty)$. More precisely, we take as state
space the set $M^1(0,\infty)$ of Radon measures $\mu$ on
$(0,\infty)$ that satisfy the condition $\int x\mu(dx)<\infty$.
(Note that we do not keep track of families of size $0$.)

    The corresponding normalized
population mass in $B^{(N)}_\ell$ (the ``radius $\ell$ block
average'') is given by
\begin{equation}\label{zeta}
\zeta^{(N,\theta)}_\ell(t)=\int x\eta^{(N,\theta)}_\ell(t,dx)=
\sum_{j}N^{-\ell}jn^{(N,\theta)}_{\ell}(N^{\ell/2}t,j)=
N^{-\ell}Z^{(N, \theta)}_{\ell}(N^{\ell/2}t),
\end{equation}
and in terms of
$\Psi^{(N,\theta)}(t)$:
\begin{equation}\label{zeta1}
\zeta^{(N,\theta)}_\ell(t)=\frac{1}{N^\ell}\int
\sum_{x\in\Omega_N,
|x|\leq\ell}\mu(x)\Psi^{(N,\theta)}(N^{\ell/2}t,d\mu).
\end{equation}

\subsubsection{Convergence theorem}

We now state our main result  that makes precise the sense in
which the entrance law described in section \ref{HEBP}
approximates the two level spatial equilibrium in $\Omega_N$
obtained in Proposition \ref{4+dequil} when the parameter
$N\to\infty$ and the random walk satisfies (\ref{trans2}).

\begin{theorem}\label{commute} (Hierarchical mean-field limit)
Consider a sequence  $(c_\ell)$ satisfying conditions
(\ref{cbound}) and (\ref{trans2}) for strong transience of the
$(2,(c_\ell), N)$-random walk.
 For fixed $N\in \mathbb N$ obeying (\ref{cbound2N}), and a sequence of nested blocks $B_\ell^{(N)}$
in $\Omega_N$, let
$\{\zeta_{\ell}^{(N,\theta)}(0)\}_{\ell\in\mathbb{N}}$
be the {\em radius $\ell$ block averages} (defined in (\ref{zeta}))
of an equilibrium  two-level branching system with an expected number
$\theta$ of particles per site.  Let
$\{\zeta_{\ell}^{\theta}\}_{\ell \in\mathbb{N}}$ be the entrance law provided by
Proposition \ref{entrancela}, case b). Then
\[\{\zeta_{\ell}^{(N,\theta)}(0)\}_{\ell\in \mathbb{N}}
\Longrightarrow \{\zeta_{\ell}^{\theta}\}_{\ell
\in\mathbb{N}}\quad \mbox { as }\;N\to\infty .\]
\end{theorem}

The proof of this result is based on the spatial ergodic theorem
for the equilibrium random field on $\Omega_N$ obtained in section
\ref{sec5.1}, a separation of time scales property derived in
section \ref{sec5.2} and a diffusion limit theorem for  the family
size processes $\{\eta^{(N,\theta)}_\ell(t)\}$ as $N\to\infty$
obtained in section \ref{sec5.3}. Using these results
 the proof of Theorem \ref{commute} is given in section \ref{sec5.4}.

\section{Super subcritical Feller branching}

In this section we continue the investigation of diffusion limits of two-level branching populations {\em without geographical structure}, which were introduced in \cite{DH}. In our case, these are superprocesses whose basic process is a subcritical Feller branching diffusion killed at 0 (this killing corresponds to the removal of void families).   With a view towards the application to the hierarchically structured geographical model, we will concentrate in subsection \ref{DL} on an initial condition of {\em many small families}, which in the diffusion limit corresponds (on a heuristic level) to an intial condition $\infty \delta_0$. In subsection \ref{DLI} we investigate time stationary super subcritical Feller branching processes which arise as diffusion limits of two-level branching populations with a high-rate immigration of individuals. The simplest situation is to think of each immigrant individual founding a new family; in the diffusion limit this leads to {\em super subcritical Feller branching diffusions with immigration} of
$\infty \delta_0$ at a constant rate (abbreviated by SSFBI). Again with a view towards the geographical model, we will consider the situation where (only) every once in a while a newly immigrant individual belongs to an already existing family. If this happens relatively rarely, then the diffusion limt remains to be SSFBI, see Proposition \ref{immprocessconvergence} and Corollary \ref{diffapproxc}.

\subsection{Diffusion limit of two-level branching particle systems}\label{DL}
\setcounter{equation}{0}

For $c >0$ and $\varepsilon >0$, consider the $M^1(0,\infty)$-valued family-size process $\{\tilde H^\varepsilon(t,dx)\}$
of a two-level branching particle system (without geographical structure) with branching rates equal to $1/\varepsilon$
at both levels and subcritical at the individual level with subcriticality parameter $\varepsilon c$. (An example  is the local family size process $H_\ell^{(N,\theta)}$ (defined in subsection \ref{LNFSP}) run at time scale $N^{\ell/2}$ and {\em with immigration suppressed}; here, $c=c_\ell$ and
$\varepsilon = N^{-\ell/2}$.)

Consider the rescaled family-size process
\begin{equation}\label{etatilde}
\widetilde \eta^{\varepsilon}(t,(x_1,x_2)) := \varepsilon \widetilde H^\varepsilon(t,(x_1/\varepsilon,x_2/\varepsilon)), \, t > 0.
\end{equation}

\begin{proposition} \label{convergence} Let
$\widetilde\eta^\varepsilon(t)$ be as in (\ref{etatilde}). Assume
that $\widetilde{\eta}^{\varepsilon}(0) = \varepsilon\lfloor
\frac{a}{\varepsilon}\rfloor\delta_{\varepsilon \lfloor
\frac{x}{\varepsilon}\rfloor}$ where $a>0$ and $x > 0$ are fixed. Then \\
(a) as $\varepsilon\rightarrow 0$,
\begin{equation}\label{convergencewithoutimm}
\{ \widetilde{\eta}^{\varepsilon}(t)\}_{t \ge 0} \Longrightarrow
\{ \xi(t)\}_{ t\ge 0}
\end{equation}
in the sense of weak convergence of $M_f([0,\infty))$-valued
c\`adl\`ag processes, and  $\xi(t)$ is the $M_f([0,\infty))$-valued
superprocess starting in $\xi(0)=a\delta_{x}$, whose motion is the
subcritical Feller branching process with generator $G_{c}$ given
by
\begin{equation}\label{Gc} G_{c}f=\left( \frac{1}{2}x\frac{\partial^2}{\partial
x^2}-c x\frac{\partial}{\partial x}\right)f,
\end{equation}
 acting on functions
$$f \in C^2_{0}([0,\infty)):=\{f \in C^2([0,\infty)):
 \lim_{x\to\infty} f(x)=0,\;f(0)=0 \}.$$ (Here $M_f([0,\infty))$ denotes
 the space of finite measures on $[0,\infty)$.)\\
 (b) The law of the process $\xi(t)$ is uniquely determined by the
Laplace functional as follows:

\begin{eqnarray} \label{Lap}
\mathbb E_{a\delta_x}\left(
\exp(-\int_{\mathbb{R}^{+}}f(y)\xi(t,dy))\right)
&&=\exp(-\int V_tf(y)\,\xi(0,dy))\nonumber\\&&=\exp(- au(t,x)),
\end{eqnarray}
where $V_tf(x)=u(t,x)$ is the unique solution of the
non-linear p.d.e.
\begin{equation}\label{LT}
\frac{\partial u(t,x)}{\partial t}   = G_{c_\ell}u(t,x)-\frac{1}{2}u(t,x)^{2}
\end{equation}
$$
u(0,x)   =f(x).
$$
\end{proposition}
\noindent\textbf{Proof} This  is essentially Theorem 4.1 of
\cite{DH}.
\bigskip

\medskip

\noindent\textbf{Remark} As an application of (\ref {Lap}) and
(\ref{LT}) we obtain the compact support property of $\xi$ appearing in Proposition \ref{convergence}.
      Assume that $\xi(0)$ has compact support and let $\bar
R_t$ denote the range of $\xi(s)$ up to time $t$. Then following
the method of Iscoe as in Theorem 1.8 of \cite{LS} or Theorem A of
\cite{DLM} one can show that $\bar R_t$ is bounded almost surely.
This involves showing (as in \cite{DLM}) that for $c>0$ the
equation
\begin{eqnarray*}
&&\frac{1}{2}x\frac {\partial^2 u}{\partial x^2}=xc
\frac{\partial u}{\partial
x}+\frac{1}{2}u^2,\\
&& \left. u(0)=0,\quad \frac{\partial u}{\partial
x}\right|_{x=0}=\alpha
\end{eqnarray*}
has for any $\alpha>0$ a blow-up at some finite $x$.

\subsubsection{Evolution equation and entrance law}\label{EEL}

In Proposition \ref{immprocessconvergence}  we will prove an
extension  of Proposition \ref{convergence}  which includes
immigration. In this subsection we obtain some properties of the
solution of the evolution equation that will be used there.

Let
\begin{equation}\label{C0}
                     C_{1}(0,\infty):= \{f\in C(0,\infty),\,
                     \lim_{x\to \infty}f(x)=0,\;|f(x)|\leq \mathrm{const}\cdot x\}
\end{equation}
and  $(T_{t})$ be the semigroup of the Feller branching diffusion
with subcriticality parameter $c$, absorbed at zero.

Let
          $V_{t}f$ be the solution of
\begin{equation}\label{Vdifferential}
           \frac{\partial}{\partial t} V_{t}f = \frac 12 x \frac{\partial
           ^2}{\partial x^2} V_{t}f -cx \frac{\partial}{\partial x}
V_{t}f - \frac
           12 (V_tf)^2, \quad V_{0}f = f.
\end{equation}
\begin{lemma}
Let $f\in C_1(0,\infty)$, $f\geq 0$. Then for $t>0$, $V_tf(x)$ is
differentiable at zero and
\begin{equation}
(V_tf)^\prime(0)=\int f(y)\kappa_t(dy)-\frac 12\int_0^t \int
(V_sf)^2(y)\kappa_{t-s}(dy)\, ds
\end{equation}
where $\kappa_t$ is the $(T_t)$-entrance law given by
(\ref{density}) in the Appendix.
\end{lemma}
\begin{proof} First  note if $f\in C^+_1(0,\infty)$ then  $V_tf(x)\leq T_tf(x)$ so that
$V_t$ maps $C^+_1(0,\infty)$ into itself.  Using the evolution form of (\ref{Vdifferential}), \begin{equation}\lim_{\varepsilon \downarrow 0} \frac
{V_tf(\varepsilon)}{\varepsilon}=\lim_{\varepsilon\downarrow 0}
\frac {T_tf(\varepsilon)}{\varepsilon}-\frac 12\int_0^t
\lim_{\varepsilon\downarrow 0}\frac
{T_{t-s}(V_sf)^2(\varepsilon)}{\varepsilon}ds\end{equation} The
result then follows from (\ref{dTz}).
\end{proof}
\begin{proposition}\label{xi0} Let $\xi^\varepsilon(t), t \ge 0$
                     denote the
super subcritical Feller branching diffusion (without immigration)
process starting in $\varepsilon^{-1} {\delta_{\varepsilon}}$ at
time $0$. Then as $\varepsilon \to 0$, $\xi^\varepsilon$ converges
 in the sense of weak convergence of
$M^1(0,\infty)$-valued continuous processes on the time interval
$[t_0,\infty)$ for all $t_0 > 0$ to a measure-valued  diffusion $ \xi^0$ where, for
all $t > 0$, $\xi^0(t)$ is an infinitely divisible random measure
with Laplace functional given by
\begin{equation}\label{Lapfunc}
\mathbb E \exp\left(-\langle  \xi^0(t),f \rangle\right) =
\exp\left(-(V_{t}f)'(0)\right).
\end{equation}
\end{proposition}
\begin{proof}
 Because  $V_{t}f(0) = 0$ we have
\begin{eqnarray}
                     \mathbb E \exp\left(-\langle  \xi^0(t),f
\rangle\right) &=& \lim\limits_{\varepsilon \downarrow 0}\mathbb E
\exp\left(-\langle \xi^{\varepsilon}(t),f \rangle\right) \\
\nonumber &=& \lim\limits_{\varepsilon \downarrow 0}\mathbb
E_{\varepsilon^{-1}\delta_{\varepsilon}} \exp\left(-\langle
\xi(t),f \rangle \right)\\ \nonumber &=& \lim\limits_{\varepsilon
\downarrow 0}
\exp\left(-\varepsilon^{-1}(V_{t}f)(\varepsilon)\right)
\\ \nonumber &=&\lim\limits_{\varepsilon \downarrow 0}
\exp\left(-\varepsilon^{-1}((V_{t}f)(\varepsilon)-(V_{t}f)(0)\right)
\\&=& \nonumber \exp\left(-(V_{t}f)'(0)\right).
\end{eqnarray}
\hfill\end{proof}

\begin{remark}\label{defK}
                     Since $\xi^0(t)$ is infinitely divisible, its
Laplace transform
                     must be of the form
\begin{equation}\label{Laplinfdiv}
                 \mathbb E \exp\left(-\langle  \xi^0(t),f \rangle\right)=
                 \exp\left(-\int(1-e^{-\langle m, f\rangle})K_{t}(dm)\right),
                 \quad f\in C_{1}((0,\infty)).
                 \end{equation}
for some uniquely determined measure $K_{t}$ on $ M(0,\infty)$, the space of Radon measures on $(0,\infty)$.
The measure  $K_{t}$ is the canonical measure of $\xi^0(t)$.
\end{remark}

We can now put these results into the framework of \cite{LS}.

A crucial property of  the entrance law $(\kappa_{t})$ given by (\ref{kappa}), which follows immediately
by partial integration from the density (\ref{density}), is given by the
following lemma.
\begin{lemma}\label{lemmaaboutkappa}
For all bounded continuously differentiable functions $g$ on
$[0,\infty)$ with $g(0) = 0$,
\begin{equation}\label{propertyofkappa}
                     \lim_{t \to 0} \int_{0}^\infty g(x)\kappa_{t}(x)dx
= g'(0).
\end{equation}
\end{lemma}

We fix a strictly positive function $\rho \in D(G_{c})$ with
\begin{equation}\label{rho}
                     \rho(x)=x \mbox{ for } x \in (0,\frac{1}{2}],
\quad \rho(x)= 1 \mbox{ for }
                     x\ge 1.
\end{equation}
Note that such a $\rho$ meets condition [A] in \cite{LS}. We take
as state space $M_{\rho}:=\{\mu\in M(0,\infty):\int
\rho(x)\mu(dx)<\infty\}$.

Following \cite{LS} we put
\begin{equation}\label{Crho}
                     C_{\rho}(0,\infty):= \{f \in C(0,\infty): |f| \le \mathrm{const } \rho,\, \lim_{x\to\infty}f(x)=0\},
\end{equation}

                  \begin{equation} \label{Drho}
                     D_{\rho}(G_{c}):= \{f \in D(G_{c}): f,\,
                     G_{c}f \in C_{\rho}(0,\infty)\}
\end{equation}
                     and
\begin{equation}\label{kappa0}
                    \kappa_{0^+}(g) := \lim_{t \to 0} \int_{0}^\infty
                    g(x)\kappa_{t}(x)dx,
\quad g \in D_{\rho}(G_{c}).
\end{equation}

Combining (\ref{Lapfunc}), (\ref{Laplinfdiv}), (\ref{propertyofkappa}) and (\ref{kappa0}) we obtain
\begin{equation}\label{identif}
                     \int(1-e^{-\langle m, f\rangle})K_{t}(dm) = (V_tf)'(0) =
                     \kappa_{0^+}(V_{t}f),\quad f\in
                     C_{1}((0,\infty)),
\end{equation}
where $V_tf$ is the solution of (\ref{Vdifferential}).

\subsection{Super subcritical Feller branching diffusion with individual immigration}\label{SSFBI}
\setcounter{equation}{0}
\subsubsection{Diffusion limit with immigration}\label{DLI}

We now extend
Proposition \ref{convergence} to include immigration, taking a fixed $t_0$  as origin of time.
Since in our
application the population from which the immigrants come is
structured into families that undergo family branching we
incorporate multitype immigration and label the set of possible
families of immigrants by $I:=[0,1]$.

 Let
$M^{1}(I\times (0,\infty))$ denote the set of Radon measures $\mu$
on $I\times (0,\infty)$ satisfying $\int_{I\times (0,\infty)}
x\mu(dy,dx)<\infty$. We denote the single atom measure
corresponding to one individual of type $y^\varepsilon_k\in I$ by
$\delta_{y^\varepsilon_k,1}$.

Consider the $M^1(I\times (0,\infty))$-valued family-size process
$\{H^\varepsilon_I(t,dy,dx)\}_{ t_0\leq t}$ with branching rates equal to
$1/\varepsilon$ at both levels,
    critical at the family level and
subcritical at the individual level with subcriticality parameter
$\varepsilon c$ (i.e. a mean offspring number  of $1-\varepsilon c$ per branching event),
and with immigration of individuals of type
$y^\varepsilon_k\in I$, given by
$\delta_{y^\varepsilon_k,1},\;k\in \mathbb{N}$, at rate
$ca^\varepsilon_k/\varepsilon^2$ with $\sum_k
a^\varepsilon_k =a$ and $\lim_{\varepsilon\to 0}\sup_k
a^\varepsilon_k = 0$. (The motivation for this comes from our geographical model, we will see in Section \ref{sec5} that this setting corresponds to the situation where the surrounding population, which serves as the source of immigration, is thought to have a frozen family structure.)
Consider the  rescaled process
$\{\eta_I^\varepsilon(t)\}_{t_0\leq t}$ defined by

\begin{equation} \label{rescaledeta}
\eta^{\varepsilon}_I(t,\{y^\varepsilon_k\}\times (x_1,x_2)):=
\varepsilon H^\varepsilon_I(t,\{y^\varepsilon_k\}\times (\frac{
x_1}{\varepsilon},\frac{x_2}{\varepsilon})).\end{equation}
\begin{proposition}\label{immprocessconvergence}
Assume that as $\varepsilon\to 0$,
$\eta_I^{\varepsilon}(t_0)\Rightarrow \mu_0\in M^1(I\times
(0,\infty))$ and $\alpha^\varepsilon \to \alpha$ in the sense of
weak convergence of finite measures on $I$, where $\alpha$ is a
nonatomic measure whose total mass $a = \alpha(I)$ plays the role of the overall immigration rate.
 Then as $\varepsilon \to 0$,
$$\{\eta_I^{\varepsilon}(t,\cdot)\}_{t_0\leq t}  \Longrightarrow \{\xi_I(t)\}_{t_0\leq t} $$  in
the sense of weak convergence of $M^1(I\times (0,\infty))$-valued
c\`adl\`ag processes on the time interval $[t_0,\infty)$, where $\{\xi_I(t)\}$ is the
measure-valued diffusion with generator
\begin{eqnarray}\label{generator}
\mathfrak{ G}F(\mu)&= f'(\langle \mu,\phi\rangle)\langle \mu, G_c
\phi\rangle +\frac{1}{2}f''(\langle \mu,\phi\rangle)\langle
\mu,\phi^2\rangle\\& +f'(\langle
\mu,\phi\rangle)c\int\frac{\partial\phi}{\partial
x}(y,x)|_{x=0}\alpha(dy),\nonumber
\end{eqnarray}
$G_c$ is the operator given by (\ref{Gc})
and $ D(\mathfrak{G})$ denotes the class of functions of the form
$F_{f,\phi}(\mu)=f(\langle\mu,\phi\rangle)$ where $f\in
C_b^3(\mathbb{R})$ and $\phi$ is continuous function on
$I\times(0,\infty)$ with
$\frac{\partial^2\phi(y,x)}{\partial x^2}$  bounded and
continuous on $I\times [0,\infty)$ and $|\sup_y \phi(y,x)|\leq
\mathrm{const}\cdot x$.
\end{proposition}
\begin{proof} $\eta^{\varepsilon}$ is given by a pregenerator
$\mathfrak{G}^\varepsilon$ acting on the class of bounded
continuous functions $F$   on $M^1(I\times(0,\infty))$  given by
\begin{eqnarray}\label{Gepsilon1}
\mathfrak{G}^\varepsilon F(\mu)&&= \sum_k
\frac{ca^\varepsilon_k}{\varepsilon^2}[F(
\mu+\varepsilon\delta_{y^\varepsilon_k,\varepsilon})
-F(\mu)]\\
&&+\frac12 \sum_k \sum_{j=1}^\infty
[F(\mu+\varepsilon\delta_{y^\varepsilon_k,j\varepsilon})-F(\mu)]
\frac{\mu(y^\varepsilon_k,j\varepsilon)}{\varepsilon^2}\nonumber\\
&&+\frac12 \sum_k\sum_{j=1}^{\infty}[F(\mu-\varepsilon\delta_{y^\varepsilon_k,j\varepsilon})-F(\mu)]
\frac{\mu(y^\varepsilon_k,j\varepsilon)}{\varepsilon^2}\nonumber\\
&&+\frac12(1-\varepsilon
c)\sum_k\sum_{j=1}^{\infty}[F(\mu-\varepsilon\delta_{y^\varepsilon_k,j\varepsilon}
+\varepsilon\delta_{y^\varepsilon_k,(j+1)\varepsilon}) -F(\mu)]
\frac{j\mu(y^\varepsilon_k,j\varepsilon)}{\varepsilon^2}\nonumber\\
&&+\frac12(1+\varepsilon
c)\sum_k\sum_{j=1}^{\infty}[F(\mu-\varepsilon\delta_{y^\varepsilon_k,j\varepsilon}
+\varepsilon\delta_{y^\varepsilon_k,(j-1)\varepsilon})) -F(\mu)]
\frac{j\mu(y^\varepsilon_k,j\varepsilon)}{\varepsilon^2}.\nonumber
\end{eqnarray}
Here the first term comes from the immigration, the second and
third from the critical branching at the family level and the
fourth and fifth from the subcritical branching at the individual
level with subcriticality parameter $c>0$.

For $F \in  D(\mathfrak{G})$, $\mathfrak{G}^\varepsilon F$ takes the form
\begin{eqnarray}\label{Gepsilon2}
&&\mathfrak{G}^\varepsilon F(\mu)\\&=& \nonumber
\sum_k\frac{ca^\varepsilon_k}{\varepsilon^2}[f(\langle\mu,
\phi\rangle+\varepsilon\phi(y_k,\varepsilon))-f(\langle\mu,\phi\rangle)]\\
&+&\frac12 \sum_k\sum_{j=1}^\infty
[f(\langle\mu,\phi\rangle+\varepsilon\phi(y_k,j\varepsilon))-f(\langle\mu,\phi\rangle)]
\frac{\mu(y_k,j\varepsilon)}{\varepsilon^2}\nonumber\\
&+&\frac12 \sum_k\sum_{j=1}^{\infty}[f(\langle
\mu,\phi\rangle-\varepsilon\phi(y_k,j\varepsilon))-f(\langle\mu,\phi\rangle)]
\frac{\mu(y_k,j\varepsilon)}{\varepsilon^2}\nonumber\\
&+&\frac12(1-\varepsilon c)\sum_k\sum_{j=1}^{\infty}[f(\langle
\mu,\phi\rangle-\varepsilon\phi(y_k,j\varepsilon)+\varepsilon\phi(y_k,(j+1)\varepsilon)))-f(\langle\mu,\phi\rangle)]
\nonumber\\&& \hspace{2cm}
\cdot\frac{j\mu(y_k,j\varepsilon)}{\varepsilon^2}\nonumber\\
&+&\frac12(1+\varepsilon c)\sum_k\sum_{j=1}^{\infty}[f(\langle
\mu,\phi\rangle-\varepsilon\phi(y_k,j\varepsilon)+\varepsilon\phi(y_k,(j-1)\varepsilon)))-f(\langle\mu,\phi\rangle)]
\nonumber\\&&\hspace{2cm}
\cdot\frac{j\mu(y_k,j\varepsilon)}{\varepsilon^2}\nonumber
\end{eqnarray}

Tightness of the family $\{\eta_I^\varepsilon\}_{0<\varepsilon
\leq 1}$  is proved by a standard argument as in \cite{Wu} or
\cite{DZ}. Using a Taylor expansion for the functions $f$ and
$\phi$ it can be verified that as $\varepsilon \to 0$, for
$F_{f,\phi}\in D(\mathfrak{G})$,  $\mathfrak{G}^\varepsilon
F_{f,\phi}(\mu) \to \mathfrak{G}F_{f,\phi}$.  Then by Lemma
\ref{moments} below, we obtain bounds on the third moments uniform
in $\varepsilon$. Using this and the tightness  it follows that
any limit point of the laws of the family $\{\eta^\varepsilon_I\}$
satisfies the martingale problem associated to the generator
$\mathfrak{G}$. Finally a standard argument (e.g. \cite{DP}, proof
of Theorem 1.1) shows that any solution  of this martingale
problem has Laplace functional given by
\begin{eqnarray}\label{ImmigrationLaplace}
&&\mathbb E_{\mu_0}\exp\big(-\int f(u,x)\xi_I(t,dy,dx)\big)\\&&\quad
=\exp\Big(-\int V_t f(y,x)\mu_0(dy,dx)
-c\int_0^t \int_I   \frac{\partial}{\partial
x}V_{t-s}f(y,x)\Big|_{x=0}\alpha(dy)\,ds \Big),\nonumber
\end{eqnarray} where $V_tf$ is given by the unique solution of the
nonlinear p.d.e.
\begin{eqnarray}\label{Vdifferential2}
           &&\frac{\partial}{\partial t} V_{t}f = \frac 12 x
           \frac{\partial
           ^2}{\partial x^2} V_{t}f -cx \frac{\partial}{\partial x}
V_{t}f - \frac
           12 (V_tf)^2, \\&& V_{0}f = f\in
           C^+_0(I\times(0,\infty)).\nonumber
\end{eqnarray}
Therefore there is a unique limit point and the proof is complete.
\end{proof}

For fixed $c > 0$, $a>0$ and arbitrary atomless measure  $\alpha$
on $I = [0,1]$ with total mass $a$, let $\xi_I(t,dy,dx)$ be as in
Proposition \ref{immprocessconvergence}. Consider the marginal
process
\begin{equation}\label{defSSFBI}
\xi(t,dx):= \int_I\xi_I(t,dy,dx).
\end{equation}
 We call $\xi(t)$ a
{\em super subcritical Feller branching diffusion with individual
immigration} (SSFBI) with initial state $\mu(dx)=\int_I\mu_0(dy,dx)$ (and parameters $a,\,
c $). The expression for the Laplace functional
(\ref{ImmigrationLaplace}) (with $f$ only a function of $x$) shows
that this coincides with the so-called {\em immigration process
with immigration rate $ac$ corresponding to the entrance law} $
\kappa_{t}$ (given by (\ref{cm}) in the Appendix) and starting
from zero measure at time $0$. The  existence of a superprocess with
immigration corresponding to an entrance law was first
established by \cite{LS} (Thm. 1.1). The resulting Laplace
transform of $\xi(t)$ with zero initial measure is given by
\begin{equation}\label{LTofxi}
                 \mathbb E \exp\left(-\langle  \xi(t),f \rangle\right)=
                 \exp\left(-ac\int_{0}^t\kappa_{0^+}(V_{s}f)ds\right),
                 \quad f\in C_{1}(0,\infty),
\end{equation}
see (\ref{identif}).

\begin{corollary}\label{diffapproxc}
(a)
The random measure $$\beta(t,dy):=\int_{(0,\infty)}
x\xi_I(t,dy,dx)=\sum_k b_k(t)\delta_{y_k}$$  is a purely atomic
finite random measure on $I$ in which the  atoms $
b_k(t)\delta_{y_k}$ correspond to the
aggregated mass at time $t$ coming from immigrants of family type $y_k\in I$.\\
 (b) For the corresponding family of stationary processes
$\{\bar\eta_I^{\varepsilon}(t)\}_{t\in \mathbb{R}}$, the random
measures $\bar\eta_I^{\varepsilon}(0)$ converge  to the
equilibrium for the process with generator $\mathfrak{G}$ given by
(\ref{generator}).  The equilibrium random measure, $\xi^a$ has
Laplace functional
\begin{equation}
                 \mathbb E \exp\left(-\langle  \xi^a,f \rangle\right)=
                 \exp\left(-ac\int_{0}^\infty\kappa_{0^+}(V_{s}f)ds\right),
                 \quad f\in C_{1}(0,\infty).
\end{equation}

\end{corollary}
\begin{proof}
 (a) The random measure
$\beta(t,dy)$ on $I$ has independent increments and  no
fixed atoms and is therefore purely atomic (see [K],
Chapt. 7). \\
(b)  Given $t_0<0$,  $\eta_I^{\varepsilon}(0)$ can be decomposed
into two parts - one coming from the initial value at $t_0$ and
one from the immigration in the interval $(t_0,0)$.  From
Proposition \ref{immprocessconvergence} it follows that the immigration parts converge  (in
the sense of weak convergence of probability measures on the space
of c\`adl\`ag functions $D([t_0,0],M^1(I\times (0,\infty))$) to the
diffusion limit with immigration, that is the process with
generator $\mathfrak{G}$. Next note that for each $\varepsilon
>0$ the contribution to the aggregated measure at time $0$
from the state at time $t_0$  is stochastically decreasing to zero
due to the subcriticality and the contribution from immigration on
$(t_0,0)$ is stochastically increasing as $t_0\downarrow -\infty$.
Moreover, using the moment bounds from Lemma \ref{moments}, it
follows that the family of random measures
$\bar\eta_I^{\varepsilon}(0)$ is tight. Therefore we have
convergence to $\bar \xi_I(0)$, the equilibrium state for the
process with generator $\mathfrak{G}$. The representation for the
Laplace functional of $\xi^a$ follows by letting $t\to\infty$ in
(\ref{LTofxi}).
\end{proof}

\begin{remark}\label{bigfamilies}
1. Corollary \ref{diffapproxc} (a)  implies that for all $\delta
>0$, asymptotically as $\varepsilon\to 0$ only a finite number of
immigrant families contribute all but $\delta$ of the mass. Each
atom corresponds to an excursion from zero for the SSFBI process
and consists of descendants of only one immigrant family. In fact,
we will see that asymptotically at the particle level each immigrant family
corresponds to the descendants of one immigrating particle.

2. Note that the assumption $\eta_I^{\varepsilon}(t_0)\Rightarrow
\mu\in M^1(I\times (0,\infty))$ in Proposition \ref{immprocessconvergence}
puts constraints not only on the
aggregated mass but also on the family structure of the
population. To understand what happens if this condition is not
satisfied  consider $\eta_I^{\varepsilon}(t_0)= \sum
a^{\varepsilon}_k\delta_{x^{\varepsilon}_k}$ with $\sum_k
a^{\varepsilon}_k x^\varepsilon_k =\theta$ but $\inf_k
x^\varepsilon_k \to \infty$ as $\varepsilon\to 0$. In this case at
times $t>t_0$, $\eta_I^{\varepsilon}(t)\to 0$ due to  ultimate
extinction of the critical family level branching.  This
observation is used below to prove by contradiction that the
equilibrium populations in $B^{(N)}_\ell$ are asymptotically
composed of families of size $O(N^{\ell/2})$.

3. Similarly, if the immigration mechanism is such that it feeds a
few large families rather than giving small new families a chance,
then in the time stationary process the family branching makes
everything extinct as $\varepsilon \to 0$.
\end{remark}

\subsubsection{Moments}

The following lemma was used in the proof of Corollary \ref{diffapproxc} and will also be needed below.
For the ease of notation  we put $t_0 = 0$, otherwise we would have to replace $t$ by $t-t_0$.
\begin{lemma} \label{moments} For $t>0$, let $\eta^\varepsilon(t,dx):= \int_I\eta_I^\varepsilon(t,dy,dx)$, where
$\eta^\varepsilon_I$ is as in
(\ref{rescaledeta}).   Let
\[m_{j,k}(t)=\mathbb E\left[\langle \eta^\varepsilon(t),
x^j\rangle^k\right], \;j,k\in\mathbb{N},\]
\[M(t)=\mathbb E(\langle \eta^\varepsilon(t),x \rangle\langle \eta^\varepsilon(t),x^2 \rangle),\]
and $o_\varepsilon(1)$ denote a term that is uniformly bounded in
$\varepsilon$ and converges to $0$ as $\varepsilon\to 0$, and
$\tilde o(t)$ denote a term that is uniformly bounded in $t\ge 0$,
and $\tilde o(t)$ converges exponentially fast to $0$ as $t\to
\infty$, and $|\tilde o(t)|\leq \mathrm{const}\cdot t $ for small
$t>0$. Then
\newline (a) \[m_{1,1}(t)=a(1-e^{-ct})+ e^{-ct}m_{1,1}(0), \]
\begin{eqnarray}&& m_{2,1}(t)=\frac{a}{2c}
[1-2e^{-ct}+e^{-2ct} ]+
m_{2,1}(0)e^{-2ct}+\frac{m_{1,1}(0)}{c}(e^{-ct}-e^{-2ct})\nonumber\\&&\qquad\quad\quad+\frac{\varepsilon
a}{2}(1-e^{-2ct} ), \nonumber\end{eqnarray}
\newline
\[m_{3,1}(t)=\frac{a}{2c^2}+o_\varepsilon(1) +\tilde o(t),\]
\[ m_{4,1}(t)=\frac{3a}{4c^3}+o_\varepsilon(1) +\tilde o(t),\]
(b) \begin{eqnarray}&& m_{1,2}(t)= m_{1,2}(0)e^{-2ct}\nonumber\\&&
+\frac{a}{4c^2}\left\{1-4e^{-ct}+2ct
e^{-ct}+3e^{-2ct}\right\}\nonumber\\&& +a^2
\left\{1-2e^{-ct}+e^{-2ct}\right\}\nonumber\\&&
+\frac{m_{1,1}(0)}{c^2}\left\{e^{-ct}-cte^{-ct}+2ac^2e^{-ct}-e^{-2ct}-2ac^2e^{-2ct}\right\}\nonumber\\&&
+\frac{m_{2,1}(0)}{c}(te^{-2ct})+\frac{\varepsilon^2}{2}\left\{
a(1-e^{-2ct})\right\}\nonumber\\&&+
\frac{\varepsilon}{4c}(3a-2e^{-ct}+2ate^{-2ct}+4e^{-ct}(-a+m_{1,1}(0))+ae^{-2ct}-4m_{1,1}(0)e^{-2ct}),\nonumber
\end{eqnarray}
\[ m_{2,2}(t)=\frac{3a}{16c^4} +o_\varepsilon(1) +\tilde o(t),\]
  \begin{eqnarray} M(t)=
\frac{a}{4c^3}+\frac{a^2}{2c}+o_\varepsilon(1) +\tilde o(t),
\nonumber
\end{eqnarray}
 \[
m_{1,3}(t)=\frac{3a^2}{4c^2}+a^3+\frac{a}{12c^4}+o_\varepsilon(1)+\tilde
o(t).\]
\end{lemma}
\begin{proof} The proof is obtained by applying the martingale
problem with the generator given by (\ref{Gepsilon1}) and
  (\ref{Gepsilon2}) to functions of the form
$F(\mu)=f(\langle \mu,\phi\rangle)$ or
$F(\mu)=f(\langle\mu,\phi_1\rangle,\langle\mu,\phi_2\rangle)$ to
derive  the following moment equations:
\[
\frac{dm_{1,1}(t)}{dt}= ca -c m_{1,1}(t),
\]
\[
\frac{dm_{2,1}(t)}{dt}=m_{1,1}(t)-2cm_{2,1}(t)+ca\varepsilon,
\]
\[
\frac{dm_{3,1}(t)}{dt}=ca\varepsilon^2+3m_{2,1}(t)-3cm_{3,1}(t)+o_1(\varepsilon),
\]
\[\frac{dm_{4,1}(t)}{dt}=ca\varepsilon^3+6m_{3,1}(t)-4cm_{4,1}(t)+o_1(\varepsilon),
\]
\[
\frac{dm_{1,2}(t)}{dt}=m_{2,1}(t)-2cm_{1,2}(t)+(2ca+\varepsilon)m_{1,1}(t)+\varepsilon^2ca,
\]
\begin{eqnarray}
\frac{dM(t)}{dt}= cam_{2,1}(t)
+m_{3,1}(t)+m_{1,2}(t)-3cM(t)+o_1(\varepsilon),\nonumber
\end{eqnarray}
\[  \frac{dm_{2,2}(t)}{dt}=m_{4,1}(t)-4cm_{2,2}(t)+o_1(\varepsilon)
\]
\begin{eqnarray} \frac{dm_{1,3}(t)}{dt}= 3cam_{1,2}(t)+3M(t) -3cm_{1,3}(t) +o_1(\varepsilon).\nonumber
\end{eqnarray}
Note that the coefficients of the $o_1(\varepsilon)$ terms only
contain moments lower in the hierarchy and hence are
asymptotically negligible. The results were obtained by solving
this linear system using MAPLE.
\end{proof}

\begin{remark}\label{inhimm}  In the case in which we replace the constant
immigration rate $a$ by a random function of time $a(\cdot)$ the
expression for $m_{1,2}(t)$ becomes
\begin{eqnarray}
m_{1,2}(t)=&&
m_{1,2}(0)e^{-2ct}\\&&+\frac{m_{1,1}(0)}{c^2}\left\{e^{-ct}-cte^{-ct}+2c^2ae^{-ct}-e^{-2ct}-2ac^2
e^{-2ct}\right\}\nonumber
\\&&+\frac{1}{c^2}\int_{0}^t
k_1(t,s)a(s)ds + \int_{0}^t\int_{0}^{s_2} k_2(t,s_2,s_1)
a(s_1)a(s_2)ds_1ds_2\nonumber\\&& + o_\varepsilon(1)\cdot\int_0^t
k_3(t,s)a(s)ds \nonumber
\end{eqnarray}
where $k_i(t,\cdot),\;i=1,2,3$ are bounded non-negative kernels
satisfying \begin{equation}\label{kconditions} \sup_{t}\int_0^t
k_i(t,s)ds <\infty,\;i=1,3,\;\;\sup_{t}\int_0^t\int_0^{s_2}
k_2(t,s_2,s_1)ds_1ds_2 <\infty
\end{equation}
and $ o_\varepsilon(1)\to 0 $ as $ \varepsilon \to 0$.
\end{remark}

\subsubsection{SPDE representation}
Let $\xi$ be an SSFBI process starting from zero measure at time $0$
 as in subsection \ref{DLI}; recall that the Laplace transform of $\xi(t)$ is given by (\ref{LTofxi}).
By an argument similar to that of \cite{LS} (Thm 1.2) it follows
that  there is a unique orthogonal martingale measure $M(ds\,dx)$
on $[0,\infty)\times (0,\infty)$ having quadratic variation
measure $\langle M \rangle (ds\,dx)= ds\,\xi(s,dx)$ such that
\begin{eqnarray}\label{inteq}
                    && \langle \xi(t),f\rangle - \langle \xi(0),f\rangle \\ && \nonumber =
\int_{0}^t
                     [\langle \xi(s),G_{c}f \rangle+ac\,\kappa_{0^+}(f)]ds +
\int_{0}^t\int_{(0,\infty)}
                     f(x)M(ds\,dx), \quad  f \in D_{\rho}(G_{c}).
\end{eqnarray}
Proceeding as in the proof of Theorem 1.7. of \cite{LS} and
tracing the arguments of \cite{KS} one infers that $\xi(t)$ has
absolutely continuous states, that is, $\xi(t,dx)=\xi(t,x)dx$, and
that one can define a time-space white noise $\dot W_{t}(x)$ on an
extension of the original probability space such that
\begin{equation}\label{martmeas} M(ds\,dx) = \sqrt{\xi(s,x)}\,\dot
W_{s}(x)ds\,dx.\end{equation}
 Moreover, $\xi(t,x)$ is almost
surely jointly continuous in $(t,x)\in [0,\infty)\times(0,\infty)$
and is a solution of the SPDE (\ref {SPDEca}) below. Note however
that in contrast to \cite{LS}, $\xi(t,x)$ does not have a finite
limit as $x\downarrow 0$. Indeed, putting
\begin{equation}\label{defZ}
           Z_{t}(x)= \int_{0}^t \int_0^\infty p_{t-s}(y,x) M(ds,dy)
\end{equation}
where $p$ is the transition density of the $c$-FBD  process (see section \ref{SFB} in the Appendix), we
obtain as in \cite{LS} ((4.10), (4.11)):
\begin{equation}\label{equxi}
           \xi(t,x) = Z_{t}(x) + \int_{0}^tac\kappa_{s}(x)ds, \quad x>0.
\end{equation}
       From (\ref{density}) we obtain
\begin{equation}\label{C(t)}
           \int_{0}^tac\kappa_{s}(x)ds
           =\frac{2ac}{x}\exp\left({\frac{-2cx}{1-e^{-ct}}}\right).
           \end{equation}

Thus, $\mathbb E (\int_{\varepsilon}^1 \xi(t,x)dx) =
C(t,\varepsilon)|\log \varepsilon|$ where $C(t,\varepsilon)$ is
uniformly bounded away from $0$, and

\begin{eqnarray}
&&\mathrm{Var}\left(\int_{\varepsilon}^1\xi(t,x)dx\right)\nonumber\\
&&=\mathbb E\left[\int_0^t\int_0^\infty f(s,y)M(ds,dy) \right]^2\nonumber\\
&&=\mathbb E\left[\int_0^t\int_0^\infty f(s,y)^2 \xi(s,y)dy
\right]\nonumber
\end{eqnarray}
where
\[ f(s,y)=\int_\varepsilon^1 p_{t-s}(y,x)dx\leq
\frac{y}{\varepsilon}\wedge 1,\quad \textrm{ and } \mathbb
E\xi(s,y)\leq
\frac{ac}{y}\exp\left({\frac{-2cy}{1-e^{-cs}}}\right),
\]
where we have used a stochastic Fubini theorem (cf. \cite{IW}
Chapt. 3, Lemma 4.1).

Therefore
\[\mathrm{Var}\left(\int_\varepsilon^1 \xi(t,x)dx\right)\leq
C_1(t)+C_2(t) |\log \varepsilon|\]
      for some positive constants
$C_1(t),C_2(t)$.
      We then  obtain for any $\delta >0$
\begin{equation}
P\left( \Big| \int_\varepsilon^1 \xi(t,x)dx-C(t)|\log \varepsilon
| \Big|> \delta C(t)| \log \varepsilon|  \right)\leq
\frac{C_1(t)+C_2(t)|\log \varepsilon|} {(\delta
C(t)|\log\varepsilon|)^2}\nonumber
\end{equation}
and therefore $\int_{\varepsilon}^1 \xi(t,x)dx$ converges in
probability to $\infty$ as $\varepsilon \to 0$. Since
$\int_{\varepsilon}^1 \xi(t,x)dx$ is monotone in $\varepsilon$
this convergence must be a.s., and thus $\xi(t,x)$ is a.s.
unbounded as $x \to 0$.

\bigskip

Finally, recalling (\ref{propertyofkappa}) and (\ref{martmeas}),
we see that (\ref{inteq}) is the integral form of the SPDE

\begin{equation} \label{SPDEca}
                     \frac{\partial}{\partial {t}} \xi(t,x) =
                     \sqrt{\xi(t,x)}\,\dot{W}_{t}(x)
                     +  G^\ast_{c} \xi(t,x)-c a\delta'_{0}(x).
\end{equation}
Let us comment on the meaning of the three terms on the right
hand side of (\ref{SPDEca}) viewed  as the limiting family size
process as $N\to\infty$ in a ball $B^{(N)}_\ell$. The first one
comes from the family branching, the second one incorporates the
individual branching and individual emigration at rate $c=c_\ell$
from $B^{(N)}_{\ell}$ (recall that $G_{c}$ is the generator of a
$c$-subcritical Feller branching diffusion), and the third term
describes immigration of small families into $B^{(N)}_{\ell}$ from
the surrounding medium  at a large rate. In fact, $\delta'_{0}$
can be viewed as the limit as $\varepsilon \to 0$ of a large
number $1/\varepsilon$ of small families of size $\varepsilon$;
note that $(1/\varepsilon)\delta_{\varepsilon}$ converges to
$-\delta'_{0}$ in the sense of Schwartz distributions on the
smooth functions vanishing at $0$.

Weak uniqueness of the solution of (\ref{SPDEca}) follows from
Proposition \ref{immprocessconvergence}; however, as in the case
of \cite{KS} it is an open question whether strong uniqueness
holds.

      The total (or aggregated) population size
      \begin{equation}\label{project}
\zeta(t) = \int_{(0,\infty)}x\, \xi(t,x)\,dx
\end{equation}
solves the equation
\begin{equation}\label{dynzeta}
d\zeta(t) = \sqrt{\int_{(0,\infty)}x^2\xi(t,x)dx}\,
dW_{t}-c(\zeta(t)-a)\,dt,
\end{equation}
which is a
one-dimensional projection of equation
(\ref{SPDEca}).
      Note that the process $\zeta$ is not Markov.
\subsubsection{Equilibrium canonical moments}\label{ECM}
 As $t \to
\infty$, the SSFBI process  $\xi(t)$ with parameters $c$ and $a$
converges in distribution to the infinitely divisible equilibrium
random measure $\xi^{a}$, cf. Corollary \ref{diffapproxc}. Writing
\begin{equation}\label{zetaeq}
\zeta^{a}:= \int_{(0,\infty)}x\xi^a(dx)
\end{equation}
for the ``aggregation'' of $\xi^a$, and $\nu_c$ for the canonical measure of $\zeta^a$,
we obtain from (\ref{LTofxi}) and (\ref{identif}):
\begin{equation}\label{Lambdac}
\nu_c ([b_1,b_2]) := c\int_0^\infty m([b_1,b_2]) K_t(dm) dt, \quad
0 < b_1 < b_2,
\end{equation}
According to Lemma
\ref{moments}, $\zeta^a$
has first and second moments
\begin{equation}\label{momentsofzeta}
\mathbb E[\zeta^{a}] = a, \quad \quad \mathbb E[(\zeta^{a})^2] = a\left(a+\frac1{4c^2}\right).
\end{equation}

\begin{definition}
Let us write $\hat \nu_c$ for the size-biasing of $\nu_c$ (cf subsection \ref{subsec5}).
\end{definition}

Because of the well-known relations (cf. Remark \ref{sbinfdiv})
$$\mathbb E[\zeta^a] = \int_{(0,\infty)}x\nu_c(dx), \quad  \mathbb E[(\zeta^a)^2] = \mathbb E[\zeta^a] \left(\mathbb E[\zeta^a]+\int_{(0,\infty)}x\hat \nu_c(dx)\right)$$
we obtain immediately from (\ref{momentsofzeta}):
\begin{remark}\label{momentsnu}
a) $\int_{(0,\infty)}x\nu_c(dx) = 1$

b) $\int_{(0,\infty)}x^2\nu_c(dx) = \frac {1}{4c^2}.$
\end{remark}

\begin{remark}
We also note that infinitely many immigrant families contribute to
$\zeta^a$.  This follows from
\begin{equation}\label{infinite}
\nu_c(0,\infty)=\infty.
\end{equation}
To see this note that from (\ref{Lambdac}) and (\ref{Laplinfdiv}),
\[ \nu_c(0,\infty)=c\int_0^\infty (-\log P(\xi^0_t=0))dt.\]
For  $\delta >0, \;x >0$ let $ f_\delta (x)=
     \frac {x}{\delta} \wedge 1 $.
Recalling (\ref{Vdifferential}) and (\ref{Lapfunc})  note that
\begin{eqnarray}
P(\xi^0_t=0)=\lim_{\theta\to\infty}\lim_{\varepsilon\to
0}e^{V_t(\theta
f_\delta)(\varepsilon)/\varepsilon}.\nonumber
\end{eqnarray}
Then by a simple modification of \cite{DH}(6.10), for any $\delta
>0$ and $\theta
>1$
\begin{eqnarray}
\lim_{\varepsilon \to 0} \frac{V_t(\theta
f_\delta)(\varepsilon)}{\varepsilon}\geq e^{-t}
\lim_{\varepsilon\to 0}\frac{T_t
f_\delta(\varepsilon)}{\varepsilon}\geq e^{-t}\int_{\delta}^\infty
\kappa_t(x)dx.\nonumber
\end{eqnarray}
Therefore by (\ref{density})
\[ \nu_c(0,\infty)\geq c\lim_{\delta \downarrow 0}\int_0^\infty e^{-t}\big( \int_\delta ^\infty
\kappa_t(x)dx \big) dt = \infty.\]
\end{remark}

\section{The genealogy of jumps in a cascade of subordinators}
\label{sec4} In this section we will carry out the program
outlined in Subsection \ref{Cassub} to obtain a representation for
the sequence $\{\zeta^\theta_\ell\}$ in terms of a cascade of
subordinators and then use this representation to obtain a
genealogical description of the population.

\subsection{Propagation of equilibria}\label{sec4.1}
\setcounter{equation}{0}

In the preceding section we encountered the equilibrium
distribution for an $\mathbb R_{+}$-valued process $\zeta(t)$,
whose dynamics is given by  (\ref{dynzeta}). A simpler situation
is the one corresponding to one-level branching where we have the
equilibrium, $\zeta^a$  of
\begin{equation}\label{FBDa}
d\zeta(t) = \sqrt{\zeta(t)}dW(t)-c(\zeta(t)-a)dt,
\end{equation}
recall subsection \ref{HEBP}. In this section we will derive
structural results which are common to {\em both} situations. We
therefore denote the equilibrium states of (\ref{FBDa}) and
(\ref{dynzeta}) by the same symbol $\zeta^a$.

In both situations the dynamics has two parameters $a$ and $c$,
and the equilibrium distribution is infinitely divisible with
expectation $a$. Therefore as in (\ref{subrepn}) this equilibrium
distribution has a representation as
\begin{equation}
\mathcal L(\zeta^a) =  \mathcal L(S(a)),
\end{equation}
where $S(\tau), \tau\ge 0,$ is a subordinator with $\mathbb E
S(\tau) \equiv \tau, \, \tau \ge 0$.

Let us denote the  L\'evy measure of $S$ by $\mu$, and note that
$\mu$ has expectation 1. Note  that $\mu = \nu_{c}$ (defined in
subsection \ref{ECM}) if $\zeta$ follows the dynamics
(\ref{dynzeta}), and $\mu = \gamma_c$ (given by
(\ref{gammaexplicit})) if  $\zeta$ follows the dynamics
(\ref{FBDa}).

$S(a)$ has a L\'evy-Khinchin representation as a Poissonian
superposition
\begin{equation}\label{superpos}
             S(a) =^d \sum_{i:\tau_{i}\le a} y_{i},
             \end{equation}
where $\sum_{i}\delta_{(\tau_{i},y_{i})}$ is a Poisson population
on $\mathbb R_{+}\times (0,\infty)$ with intensity measure
$d\tau\, \mu (dy)$. Since by (\ref{gammaexplicit}), resp.
(\ref{infinite}),
\begin{equation}\label{infinitemu}\mu(0,\infty)=\infty\end{equation}
in the two cases,
$P(S(a)=0)=0$ if $a>0$ and $S(\cdot)$ has infinitely many jumps in
any open interval.

Since  $\zeta^a$ is a Poisson superposition of immigrant clusters
(recall  (\ref{Lambdac}) and (\ref{gammac})), the representation
(\ref{superpos}) has a genealogical interpretation: the summands
$y_{i}$ measure the size of those parts of $\zeta^a$ which trace
back to one and the same immigrant.

Our  aim in this section is to study the hierarchy (\ref{fluc})
into which (\ref{FBDa}) and (\ref{dynzeta}) are embedded. In both
situations, the parameters of the hierarchy are a sequence
$(c_{\ell})_{\ell=1,2,...}$ of positive numbers. Recall that, for
different levels $\ell$ in the hierarchy, the dynamics of
$\zeta_{\ell}(t)$ run at separated time scales, and the
equilibrium state $\zeta_{\ell+1}$ at level $\ell+1$ acts as
(random) parameter $a$ for the dynamics at level $\ell$. Instead
of one L\'evy measure $\mu$, we now have a sequence of L\'evy
measures $(\mu_{\ell})$ ( which is either $(\gamma_{c_{\ell}})$ or
$(\nu_{c_{\ell}})$).

In this way we obtain a Markovian dynamics which transports the
equilibria down the levels:
\begin{equation}
\mbox{given }\zeta_{\ell+1}=a\, , \zeta_{\ell}\mbox{ is infinitely
divisible with canonical measure }a\mu_{\ell}.
\end{equation}
Since $\mu$ has expectation 1, $(\zeta_{\ell})$ constitutes a
backward martingale. We now turn to the following problems:

a) Find a condition on $(\mu_{\ell})$ which guarantees the
existence of an entrance law, denoted by $\zeta_\ell^\theta$,  for
$(\zeta_{\ell})$ starting in $\theta > 0$ ``at level $\infty$''
and having constant expectation $\theta$.

b) Describe the ``branching genealogy'' underlying such an
entrance law.

We will answer these questions in the next subsections. Later on,
we will give a relation with the asymptotics of the genealogy of
the equilibrium branching populations on $\Omega_N$ as $N\to
\infty$.

\subsection{An entrance law from infinity}\label{sec4.3}
\setcounter{equation}{0}

Let $S_{k}$, $k=1,2,\ldots$ be independent subordinators with
L\'evy measures $\mu_{k}$. We denote the second moment of
$\mu_{k}$ by $m_{k}$. For $j > \ell$ define the random variables
\begin{equation}\label{itsubA}
S_{\ell}^{j}(a):=S_{\ell}(S_{\ell+1}(\dots(S_{j-1}(a)))),
\end{equation} and
write $\Pi_\ell(a,.)$ for the distribution
 of $S_\ell(a), \, \ell =1,2,...,\, a >0$.

 \begin{proposition}\label{existentra} If
\begin{equation}\label{summable}
           \sum_{k=1}^\infty m_{k} < \infty
           \end{equation}
           then, for each $\theta > 0$, the sequence of processes
$$(S^j_{j-1}(\theta), \, S^j_{j-2}(\theta), ...,\,S^j_{2}(\theta),\, S^j_{1}(\theta))$$
 converges as $j\to \infty$ (component-wise) in probability to a sequence
$$(..., \zeta^\theta_2,\, \zeta^\theta_1)$$
which obeys
\begin{equation}\label{flow}
S_\ell(\zeta^\theta_{\ell+1})=\zeta^\theta_{\ell}  \mbox{ a.s. for
all } \ell
\end{equation} and
\begin{equation}\label{convtotheta}
\lim_{\ell \to \infty}\zeta_\ell^\theta = \theta \mbox{ a.s. }
\end{equation}
 In
particular, the distributions $\pi_\ell = \mathcal
L(\zeta_\ell^\theta)$ are an entrance law for the backward Markov
chain with probability transition function
\begin{equation}\label{ent}
P(\zeta_\ell \in A | \zeta_{\ell+1}= a) = \Pi_\ell(a,A),
\end{equation}
 and they
are its unique entrance law with the property
\begin{equation}\label{distconvtotheta}
\pi_\ell \Rightarrow \delta_\theta \mbox{ as } \ell \to \infty.
\end{equation}
\end{proposition}
\begin{proof}
Since $S_k(a)$ is infinitely divisible with canonical measure
$a\mu$, we have
$$\mathrm{Var}(S_k(a)) = am_k.$$
Hence we obtain \begin{eqnarray*} \mathrm{Var}(S^{k+2}_k(a)) &=&
\mathrm{Var}[\mathbb E[S_k(S_{k+1}(a))|S_{k+1}(a)]] + \mathbb
E[\mathrm{Var}[S_k(S_{k+1}(a))|S_{k+1}(a)]] \\&=&
a(m_{k+1}+m_{k}).
\end{eqnarray*}
In the same way we get for all $j>k> \ell$:
\begin{equation}\label{Varjl}
\mathrm{Var}(S^{j}_k(a))= a(m_{j-1}+...+m_k)
\end{equation}
and
$$ \mathbb E(S^j_\ell(a) - S^k_\ell(a))^2 = a(m_{j-1}+...+m_k).$$
Thus, due to (\ref{summable}), for fixed $\ell$ the sequence
$(S^j_\ell(\theta))_{j>\ell}$ is Cauchy in $L^2$. We define
\begin{equation}\label{defzetaell}
\zeta_{\ell}^\theta \equiv L^2-\lim_{j\to\infty}S^j_\ell(\theta).
\end{equation}
Since $a\mapsto S_\ell(a)$ is continuous in $L^1$, we have
\begin{equation}\label{aslimit}S_\ell(\zeta_{\ell+1}^\theta) = S_\ell(\lim_{j\to
\infty}S^j_{\ell+1}(\theta)) =\lim_{j\to
\infty}S^j_{\ell}(\theta)= \zeta_{\ell}^\theta \mbox{
a.s.},\end{equation} which proves (\ref{flow}) and, a fortiori,
implies (\ref{ent}). From (\ref{defzetaell}) and (\ref{Varjl}) it
is clear that
$$\mathrm {Var} \,\zeta_\ell^\theta = \theta\sum_{k=\ell}^\infty
m_k.$$ Since $\mathbb E\zeta_\ell^\theta \equiv \theta$, this
together with (\ref{summable}) implies that $\zeta_\ell^\theta$
converges to $\theta$ in probability as $\ell \to \infty$.
Moreover, since because of (\ref{flow}) $\zeta_{\ell}^\theta$ is a
backwards martingale, this convergence is even a.s., and we have
(\ref{convtotheta}).

It remains to show the claimed uniqueness statement. For this let
$(\pi_k)$ be an entrance law for $(\Pi_k)$ obeying
(\ref{distconvtotheta}), and let $X_k, \, k= 1,2,...$ be random
variables, independent of the subordinators $S_\ell$, with
$$\mathcal L(X_j)= \pi_j.$$ From the entrance law property of
$(\pi_k)$ and the definition of $(\Pi_k)$ we have for all
$j>\ell$:
\begin{equation}\label{transportX}
\pi_\ell = \mathcal L(S^j_\ell (X_j)).
\end{equation}
On the other hand we have by monotonicity of $\tau \mapsto
S^j_\ell (\tau)$:
$$\mathbb E |S^j_\ell (X_j)- S^j_\ell (\theta)| \le \mathbb E|X_j
- \theta|.$$ From this, the claimed identity $\pi_\ell = \mathcal
L(\zeta_\ell^\theta)$  follows by Markov's inequality together
with (\ref{defzetaell}) and (\ref{transportX}).
\end{proof}
\medskip

{\bf Proof of Proposition \ref{entrancela}} \\
Recall that with the notation of Proposition \ref{entrancela},
$\Pi^{(1)}_\ell(a,.)= \mathcal L(S^{(1)}_\ell(a))$ and
$\Pi^{(2)}_\ell(a,.)= \mathcal L(S^{(2)}_\ell(a))$, where
$S^{(1)}_\ell$ and $S^{(2)}_\ell$ are subordinators with L\'evy
measures $\gamma_\ell$ and $\nu_\ell$, respectively. From Remark
\ref{momentsofgamma} we have that the second moment of $\gamma_k$
equals $1/(2c_k)$, and Remark \ref{momentsnu} shows that the
second moment of $\nu_k$ equals $1/(4c_k^2)$. The proof of
Proposition \ref{entrancela} is thus immediate from Proposition
\ref{existentra}.

\subsection{The genealogy of jumps in an iteration of
subordinators}\label{sec4.4} \setcounter{equation}{0}

The composition of subordinators gives rise to a ``genealogy'' of
their jumps. To illustrate this, consider the two subordinators
$S_{1}$, $S_{2}$, where
$$S_{1}(b) = \sum_{t_{i}\le b} (S_{1}(t_{i})-S_{1}(t_{i}-)),$$
$$S_{2}(a) = \sum_{\tau_{n}\le a} (S_{2}(\tau_{n})-S_{2}(\tau_{n}-)).$$
Then
$$S_{1}(S_{2}(a)) = \sum_{\tau_{n}\le a}\quad \sum_{S_{2}(\tau_{n}-)< t_{i}\le
S_{2}(\tau_{n})} (S_{1}(t_{i})-S_{1}(t_{i}-)).$$ In this way, the
jumps of $S_{1}$ are coagulated into families of jumps stemming
from one and the same jump of $S_{2}$.

Iterating this, we obtain from the flow property (\ref{flow}) that
\begin{equation}\label{repzeta}
\zeta^\theta_\ell = \sum_{\tau_i \in
[0,\zeta_j^\theta]}(S^j_\ell(\tau_i)-S^j_\ell(\tau_i-)) \mbox{
a.s. },
\end{equation}
where $\{\tau_i\}$ is the set of all points in $[0,
\zeta_j^\theta]$ in which $\tau \mapsto S^j_\ell(\tau)$ has a
jump. The representation (\ref{repzeta}) induces a partition of
$[0, \zeta^\theta_\ell]$ which we denote by $\mathfrak
P_{j,\ell}$. Note that for fixed $\ell$ the $\mathfrak P_{j,\ell}$
are coalescing (i.e. becoming coarser) as $j$ increases.

The sequence of coalescing partitions $\mathfrak P_{j,\ell}$
induces a graph $\mathfrak G_\ell$ as follows: The set of nodes of
$\mathfrak G_\ell$ is the union $\bigcup_{j\ge\ell}\{j\}\times
\mathfrak P_{j,\ell}$. For $n \in \mathfrak G_\ell$ we call its
first component the {\em
 level} of $n$.
For two nodes $n_1 = (j_1, I_1)$, $n_2= (j_2, I_2)$ of $\mathfrak
G_\ell$ we say that $n_1$ is an {\em ancestor} of $n_2$ if
$j_1>j_2$ and $I_2 \subseteq I_1$, and we say that $n_1$ is the
{\em parent} of $n_2$ if $n_1$ is the ancestor of $n_2$ with $j_1
= j_2+1$. The (directed) edges of $\mathfrak G_\ell$ then are all
the parent-child pairs in $\mathfrak G_\ell \times \mathfrak
G_\ell$. Say that two nodes in $\mathfrak G_\ell$ are {\em
related} if they have a common ancestor.  Then by construction of
the sequence $(\mathfrak P_{j,\ell})$ each equivalence class of
$\mathfrak G_\ell$ is a  {\em tree}, i.e. a directed connected
graph without cycles. Therefore, $\mathfrak G_\ell$ is a {\em
forest}, i.e. a union of pairwise disconnected trees. Finally, we
label each node of $\mathfrak G_\ell$ with the length of the
subinterval of $[0, \zeta^\theta_\ell]$ to which it corresponds,
thus arriving at the  random {\em labelled forest} $\mathfrak
F_\ell$ which we associate with the random sequence of coalescing
partitions $\mathfrak P_{j,\ell}$, and which encodes the genealogy
of the jumps of the process $(..., \zeta^\theta_2,
\zeta^\theta_1)$ constructed in Proposition
\ref{existentra}.

Intuitively, viewing $[0, \zeta^\theta_\ell]$ as a continuum of
individuals,  this means that two individuals $a_1, \, a_2 \in [0,
\zeta^\theta_\ell]$ belong to the same element of $\mathfrak
P_{j,\ell}$ if and only if they descend from a common ancestor (or
equivalently, from one subordinator jump) at some level less or
equal than $j$. Furthermore two individuals $a_1, \, a_2 \in [0,
\zeta^\theta_\ell]$ belong to the same element of the minimal
partition $\mathfrak P_{\infty,\ell}$ if and only if they descend
from a common ancestor at any level higher than $\ell$. Using the
independent increments property in $\theta$, (\ref{infinitemu})
and (\ref{aslimit}) it can be shown that there are countably many
distinct elements in $\mathfrak P_{\infty,\ell}$ each
corresponding to an infinite tree. Therefore we have a
decomposition of the equilibrium population into a countable set
of subpopulations each consisting of individuals having a common
ancestor.

           \subsection{The genealogy in the hierarchichal mean field limit}
           \label{sec4.5}
           \setcounter{equation}{0}
          With the special choice of $(S_\ell)$ described at the beginning
of this section, we have all reasons to conjecture that the random
labelled forests $\mathfrak F_\ell$ defined in the previous
subsection describe the genealogy of the (one- or two-level)
branching population in equilibrium as $N\to \infty$.

To make this more precise, consider a fixed sequence
$B_\ell^{(N)},\, \ell = 1,2,...$ of nested balls in $\Omega_N$,
and let $\mathcal P_\ell^{(N)}$ be that part of the equilibrium
population which lives in $B_\ell^{(N)}$. (Here and below we
suppress the notation of $\theta > 0$ which we keep fixed.) Fix
$N$ and $\ell$ for the moment. For two individuals $I_1,\, I_2$ in
$\mathcal P_\ell^{(N)}$ and $j \ge \ell$ we say that
$$ I_1 \sim_j I_2 $$
if $I_1$ and $I_2$ have a common ancestor in $B_j^{(N)}$. This
induces a partition on $\mathcal P_\ell^{(N)}$ which we denote by
$\mathfrak P_{j,\ell}^{(N)}$. The sequence $\mathfrak
P_{j,\ell}^{(N)}, \, j=\ell, \ell+1,...$ is coalescing, and we can
associate with it a labelled forest $\mathfrak F_\ell^{(N)}$ in
the same way as we associated $\mathfrak F_\ell$ with $\mathfrak
P_{j,\ell}$ in the previous subsection, the only difference being
that now we label the nodes of $\mathfrak F_\ell^{(N)}$ by the
cardinalities of their corresponding sub-populations of $\mathcal
P_\ell^{(N)}$, divided by $N^{\ell}$.

Our main result (Theorem \ref{commute}) suggests to conjecture
that, in a suitable topology,
$$\mathfrak F_\ell^{(N)}\rightarrow \mathfrak F_\ell \mbox { as } N
\to \infty.$$

\subsection{Size-biasing iterated subordinators}\label{subs4.2}
\setcounter{equation}{0}

With a view towards the genealogy of a  sampled individual (see
subsection \ref{sampledind}) we will now prove a representation of
the size-biasing of $\mathcal L (S_{\ell}^{j}(a))$, where
$S_{\ell}^{j}$ is the composition of subordinators defined in
(\ref{itsubA}). To this purpose we first consider a single
subordinator $S$ evaluated at a random argument.

\begin{proposition}\label{argmix}
           Let $S(\tau), \tau \ge 0,$ be a subordinator with L\'evy measure
           denoted by $\mu$, and let $A$ be an $\mathbb R_{+}$-valued random
           variable independent of $S$ and with finite expectation. Then the
           size-biasing of  $\mathcal L  (S(A))$ arises as the distribution
           of $S(\hat A) + \hat Y$, where  $\mathcal L(\hat A)$ is the
           size-biasing of  $\mathcal L (A)$,  $\mathcal L(\hat Y)$ is the
           size-biasing of $\mu$, and $\hat A$ and $\hat Y$ are independent.
           \end{proposition}
           \begin{proof}
         We write $\Pi_{\tau \mu}$ for the distribution of a Poisson point
         configuration on $\mathbb R_{+}$ with intensity measure $\tau \mu$.
         Then
         $$ \mathcal L(S(\tau)) =  \mathcal L_{\tau}(\langle \Psi,
         \textrm{id}_{\mathbb R_{+}}\rangle),$$
         where
         $$ \mathcal L_{\tau}(\Psi) = \Pi_{\tau \mu}.$$
         Writing $\sigma$ for the distribution of $A$, we thus have
         \begin{equation} \label{Apsi}
              \mathcal L(S(A)) =  \mathcal L_{\sigma}(\langle \Psi,
         \textrm{id}_{\mathbb R_{+}}\rangle),
         \end{equation}
         where
         $$ \mathcal L_{\sigma}(\Psi) = \int \Pi_{\tau \mu}(.)\sigma(d\tau).$$
         Our task is to compute the size-biasing of
         $\mathcal L_{\sigma}(\langle \Psi,
         \textrm{id}_{\mathbb R_{+}}\rangle)$ with
         $\mathfrak s(\psi)= \langle \psi,
         \textrm{id}_{\mathbb R_{+}}\rangle$, cf. Definition \ref{sizebiasing}.
         To this end let us first compute
         the size-biasing of $\mathcal L_{\sigma}(\Psi)$ with
             $\langle \psi,
         \textrm{id}_{\mathbb R_{+}}\rangle$, and then project. It follows
         from Corollary \ref{sbMP} in the
         Appendix that the size-biasing of
         $\mathcal L_{\sigma}(\Psi)$ with $\langle \psi,
         \textrm{id}_{\mathbb R_{+}}\rangle$ is $\mathcal L (\Phi +
         \delta_{\hat Y})$, where
         $$\mathcal L(\Phi) = \int \Pi_{\tau \mu}(.)
         \hat \sigma (d\tau),$$
         $\hat \sigma$ is the size-biasing of $\sigma$, $\mathcal L(\hat Y)$ is
         the size-biasing of $\mu$, and $\Phi$ and $\hat Y$ are independent.
         Consequently, the size-biasing of $\mathcal L_{\sigma}(\langle \Psi,
         \textrm{id}_{\mathbb R_{+}}\rangle)$ with
         $\langle \psi,
         \textrm{id}_{\mathbb R_{+}}\rangle$ is
         $$\mathcal L(\langle \Phi + \delta_{\hat Y},
         \textrm{id}_{\mathbb R_{+}}\rangle) = \mathcal L(S(\hat A)+\hat
         Y)),$$
         where $\mathcal L (\hat A) = \hat \sigma$, and $S$, $\hat A$ and
         $\hat Y$ are independent. Together with (\ref{Apsi}) this proves the
         claim.
         \end{proof}
         \begin{remark}\label{sbinfdiv}
             For deterministic $A$, Proposition \ref{argmix} renders the
             well-known fact that the size-biasing of an infinitely divisible
             distribution $\pi$ on $\mathbb R_{+}$ is the convolution of $\pi$
             with the size-biasing of the canonical measure of $\pi$.
             \end{remark}
\begin{corollary}\label{Abias}
           Let $S^j_{\ell}(a)$ be the iteration of subordinators defined
           in (\ref{itsubA}), where the $S_{k}$ are independent subordinators
           with L\'evy measures $\mu_{k}$. Then the size-biasing of $\mathcal
           L(S^j_{\ell}(a))$ arises as the distribution of $\widehat
           S^j_{\ell}(a)$ defined by
           \begin{equation}\label{repAhat}
\widehat{S}_{\ell}^{j}(a)=S_{\ell}^j(a) + \widetilde
S_{\ell}^{j-1}(\widehat{Y}_{j-1})+\widetilde
S_{\ell}^{j-2}(\widehat{Y}_{j-2})+...+\widetilde
S_{\ell}^{\ell+1}(\widehat{Y}_{\ell+1})+ \widehat{Y}_{\ell},
\end{equation}
where $\mathcal L(\hat Y_{k})$ is the size-biasing of $\mu_{k}$,
$\widetilde S^k_\ell$ is distributed as $S^k_\ell$ and all random
variables occurring on the r.h.s. of (\ref{repAhat}) are
independent.
\end{corollary}
\begin{proof}
           For $\ell=j-1$, Proposition \ref{argmix} shows that the
           size-biasing of $\mathcal L(S^j_{j-1}(a)) =
           \mathcal{L}(S_{j-1}(a))$ arises as the distribution of
        $$S_{j-1}(a) + \hat Y_{j-1},$$ where both summands are independent.
One more application of Proposition \ref{argmix} thus gives that
the size-biasing of $\mathcal{L}(S_{j-2}^j(a)) =
        \mathcal{L}(S_{j-2}(S_{j-1}^j(a))$ arises as the distribution of
        $$S_{j-2}(S_{j-1}(a) + \hat Y_{j-1}) + \hat Y_{j-2}$$
        which due to the independence of $S_{j-1}(a)$ and $\hat Y_{j-1}$
        equals in distribution to
        $$S_{j-2}(S_{j-1}(a)) + S_{j-2} ' (\hat Y_{j-1}) + \hat Y_{j-2},$$
        $S_{j-2}'$ being an independent copy of $S_{j-2}$.

        Iterating the argument we arrive at our assertion.
        \end{proof}

        \begin{remark}
           As before, let us  denote the
second moment of $\mu_{k}$ (or equivalently the first moment of
$\hat Y_{k}$)  by $m_{k}$. From (\ref{repAhat}) it follows that
\begin{equation}\label{expAhat}
\mathbb E\widehat{S}_{\ell}^{j}(a) =a+\sum_{k=\ell}^{j-1}m_{k}
\end{equation}
Hence the summability of the $m_k$ is a sufficient condition for
tightness of the $\widehat{S}_{\ell}^{j}(a)$.
\end{remark}

We now turn to the entrance law constructed in Proposition
\ref{existentra}.
\begin{corollary}
From (\ref{defzetaell}) and Corollary \ref{Abias} we obtain that
the size-biasing of $\mathcal L (\zeta_{\ell}^\theta)$ arises as
the distribution of
\begin{equation}\label{repzetahat}
\hat \zeta_{\ell}^\theta \equiv \widehat{Y}_{\ell} + \widetilde
S_{\ell}^{\ell+1}(\widehat{Y}_{\ell+1}) + \widetilde
S_{\ell}^{\ell+2}(\widehat{Y}_{\ell+2})+....+ \zeta_{\ell}^\theta,
\end{equation}
where the random variables $\widehat{Y}_{k}$ and $\widetilde
S_{\ell}^{k}$ are as in Corollary \ref{repAhat}, and all random
variables occurring on the right hand side of (\ref{repzetahat})
are independent.
\end{corollary}

We can go one step further and study the genealogical
relationships underlying the representation (\ref{repzetahat}).

To this purpose, let us study the branching dynamics on the
``populations of jumps'' induced by the composition of the
subordinators $S_k$, resuming the reasoning of subsection
\ref{sec4.4}. For each $k$ we consider a branching dynamics which
takes a counting measure $\phi_{k}$ on $\mathbb R_{+}$ into a
random counting measure $\Phi_{k-1}$ in the following way: If
$\phi_{k}= \sum_{i\in I_{k}}\delta_{y_{i}}$, then
$$\Phi_{k-1}= \sum_{i\in I_{k}} \Psi_{i},$$
where $\Psi_i$ is a Poisson counting measure on $\mathbb R_{+}$
with intensity measure $y_{i}\mu_{k-1}$, and the $\Psi_{i}$ are
independent.

Now fix two levels $j,\ell \in \mathbb N$ with $j > \ell$.
Starting with $\phi_{j} = \delta_{a}$, and iterating the branching
dynamics from level $j$ down to level $\ell$, we construct a
random path of counting measures on $\mathbb R_{+}$ which we
denote by
$$(\delta_{a},\Phi^j_{j-1}(a),\ldots,\Phi^j_{\ell}(a)) =:
H^j_{\ell}(a).$$ By keeping track which atom in $\Phi^j_{k-1}(a)$
stems from which atom in $\Phi^j_{k}(a)$, we can enrich the
history $H^j_{\ell}(a)$ to a tree $T^j_{\ell}(a)$, each of whose
nodes is marked by a non-negative real number. For example, if
$\Phi^j_{k}(a) =  \sum_{i\in I_{k}}\delta_{y_{i}}$, $j \ge k \ge
\ell,$ then the set of nodes of $T^j_{\ell}(a)$ at level $k$
corresponds to the index set $I_{k}$, and $y_{i}$ is the mark (or
``size'') of the node with index $i$. We write
$$\mathfrak s (T^j_{\ell}(a))= \langle \Phi^j_{\ell}(a), \textrm
{id}_{\mathbb R_{+}}\rangle$$ for the {\em total size of the tree}
$T^j_{\ell}(a)$ {\em at level} $\ell$, and note that
$$\mathfrak s (T^j_{\ell}(a))=^d A^j_{\ell}(a),$$
where $S^j_{\ell}(a)$ is defined in (\ref{itsubA}). Proceeding in
a similar way as in the proof of Corollary \ref{Abias} we obtain a
``spinal decomposition'' of the size-biased tree, which we state
here without proof.
\begin{proposition}\label{sizebT}
           The size-biasing of $\mathcal{L}(T^j_{\ell}(a))$ with $\mathfrak s
           (T^j_{\ell}(a))$ arises as the distribution of the superposition
           of $T^j_{\ell}(a)$ and
           $$(\delta_{\hat Y_{j-1}},\ldots, \delta_{\hat Y_{k}},
           \Phi^k_{k-1}(\hat Y_{k}),\ldots, \Phi^k_{\ell}(\hat Y_{k})), \quad k=
           j-1,\, j-2,\ldots,\ell,$$
           where $\mathcal{L}(\hat Y_{k})$ is the size-biasing
           of $\mu_{k}$, $T^j_{\ell}(a)$, $\hat Y_{j-1},\ldots, \hat Y_{\ell}$ are
           independent, and given  $\hat Y_{j-1},\ldots, \hat Y_{\ell},$ the
           $\Phi^k_{r}(\hat Y_{k})$, $k > r$, are independent.
           \end{proposition}

           Under the assumption (\ref{summable})
           of summability of the second moments of
           $\mu_{k}$, $k\in \mathbb N$, the size-biasing of the
           random labelled forest $\mathfrak F^\ell=\mathfrak F^\theta_\ell$
           (constructed in subsection \ref{sec4.4}) with respect
           to its `` size'' $\mathfrak s (\mathfrak
           F^\theta_\ell)= \zeta^\theta_\ell$ arises as the independent superposition of
           $\mathfrak F^\theta_\ell$ and $\hat T_{\ell}^{\infty, \rm{can}}$,
           where $\hat T_{\ell}^{\infty, \rm{can}}$ is constructed as
           follows:

           First build a ``spine'' $(\ldots,\hat Y_{\ell-2}, \hat Y_{\ell-1},
           \hat Y_{\ell})$, and given the spine, superimpose independently
           the trees $T^k_{\ell}(\hat Y_{k}), k \ge \ell$.

\subsection{The genealogy of relatives of a sampled individual}
\label{sampledind} In  subsection \ref{sec4.4} we fixed a ball (or
$\ell$-block) $B_\ell^{(N)}$ from the beginning. Now we take a
different viewpoint and think of an individual sampled from the
equilibrium population within a union of many $\ell$-blocks in
$\Omega_N$ from the beginning. Denote the
           chosen individual by $I$, and the $\ell$-block by $\hat
           B_{\ell}$. Recall that for large $N$ the total number of
individuals in an $\ell$-block is approximately distributed like
$N^\ell \zeta_\ell$ (see Theorem \ref{commute}), we see that the
number of individuals in the chosen block
           is approximately distributed like
           $N^\ell \hat \zeta_{\ell}^\theta$, where $\mathcal L (\hat
           \zeta_{\ell}^\theta)$ is the size-biasing of $\mathcal L (
           \zeta_{\ell}^\theta)$.

           The block $\hat B_{\ell}$ sits in a nested sequence of blocks of
           levels $\ell+1,\ell+2,\ldots$ which we denote by $\hat
B_{\ell+1},\hat
           B_{\ell+2},\ldots$ Consider the population $\hat
           \mathcal{P}_{\ell}$ of all those individuals in $\hat B_{\ell}$
           which have an ancestral family in common with the individual $I$.
           The population $\hat
           \mathcal{P}_{\ell}$ can be decomposed in a natural way according
           to its immigration history into the  $\hat B_j$, $j \ge \ell$.

           For $j>\ell$, denote by $\hat
           \mathcal{P}_{\ell}^j$ the subpopulation of all those individuals
           in $\hat
           \mathcal{P}_{\ell}$
           that have some common ancestor
           with $I$ who lived in $\hat B_{j}$ but none who
           lived in $\hat B_{j-1}$. In other words, $\hat
           \mathcal{P}_{\ell}^j$ consists of all those
           individuals $J$ which obey $J\sim_jI$ but not
           $J\sim_{j-1}I$.

            In this way we obtain a
           decomposition of $\hat
           \mathcal{P}_{\ell}$ according to the hierarchical distance of the
           (geographically) closest ancestors common with the chosen individual $I$:
           for $j>\ell$, the subpopulation $\hat
           \mathcal{P}_{\ell}^{j}$ consists of those individuals in
           $\hat B_{\ell}$ whose geographically closest common ancestor with $I$ has
           hierarchical distance $j$ from $I$.

           The size (i.e. the total number of individuals) of $\hat
           \mathcal{P}_{\ell}^{j}$ is approximately distributed as
           $$N^\ell S^{j-1}_{\ell}(\hat Y_{j-1}) = N^\ell
           S_{\ell}(S_{\ell-1}\ldots(S_{j-2}(\hat Y_{j-1})))$$
           and thus has approximate expectation $N^\ell \frac 1{4c_{j-1}^2}$.
           Hence the summability condition (\ref{summable}) (which
           corresponds to the condition for transience resp.
           strong transience of the hierarchical random walk)
         amounts precisely to an expected
finite number of relatives
           of the chosen individual in the block $\hat B_{\ell}$.

\section{The hierarchical mean field limit of two-level branching systems in equilibrium}
\label{sec5} In this section we investigate the  two level
branching equilibrium $\Psi^{(N,\theta)}(0)$  described in
Proposition \ref{4+dequil} and its limiting behavior as $N \to
\infty$. We assume that the underlying random walk is a
$(2,(c_\ell),N)$-random walk on $\Omega_N$ and  $(c_\ell)$
satisfies the strong transience conditions   (\ref{cbound2N}) and
(\ref{trans2}). Recall that $(B_{\ell}^{(N)})$ denotes a sequence
of nested blocks in $\Omega_N$.

We will see  in Lemma \ref{familysizeconstraint} that in equilibrium, asymptotically as
$N\to \infty$, $\Psi^{(N,\theta)}(0)$ consists of the order of $N^{\ell/2}$
families in $B_{\ell}^{(N)}$,  a typical such family having a random multiple of
$N^{ \ell/2 }$ individuals.

\subsection {A spatial ergodic theorem}\label{sec5.1}
\setcounter{equation}{0}

 In this subsection we first
collect some basic facts about the two level branching systems
$\Psi^{(N)}(t)$ on $\Omega_N$ which were introduced in section
2.2.6 and their equilibria. The main result will be a spatial
ergodic theorem for the aggregated equilibria; this will also be
an ingredient in the proof of Theorem \ref{commute}.

 Let $M_c(\Omega_N)\setminus \{\it o\}$ denote the space of non-zero counting measures on
$\Omega_N$ such that finite sets have finite measure, and let
$M^1(M(\Omega_N)\setminus \{\it o\})$ be the set of measures $\nu$
on $M_c(\Omega_N)\setminus \{\it o\}$ such that
$\int_{M(\Omega_N)\setminus \{\it o\}} \mu(B)\nu(d\mu) <\infty$ if
$B$ is a finite set. Then let $M_c^1(M_c(\Omega_N)\setminus \{\it
o\})$ denote the subspace of counting measures in
$M^1(M_c(\Omega_N)\setminus \{\it o\})$. The set
$M_c(\Omega_N)\setminus \{\it o\}$ carries a complete separable
metric $\rho$  which generates the restriction of the vague
topology on
 $M_c(\Omega_N)$ to  $M_c(\Omega_N)\setminus \{\it o\}$ and makes a subset
 of  $M_c(\Omega_N)$ bounded if and only if it is contained in $\{\mu | \mu(B) > 0 \mbox{ for some finite } B
\subset \Omega_N\}$, see \cite{MKM} Proposition 3.3.2. The set
$M^1(M_c(\Omega_N)\setminus \{\it o\})$ is equipped with the
topology  generated by $\nu \mapsto \int_{M_c(\Omega_N)\setminus
\{\it o\}} \mu(B)\nu(d\mu)$, where $B$ is a finite subset of
$\Omega_N$, and $\nu \mapsto \int F(\mu) \nu(d\mu)$, where
$F:M_c(\Omega_N)\setminus \{\it o\} \to \mathbb R$ is continuous,
bounded and with $\rho$-bounded support.

We can represent $\Psi^{(N)}(t)$  as a
c\`{a}dl\`{a}g  Markov process with state space
 $M_c^1(M_c(\Omega_N)\setminus \{\it o\})$. On $M_c^1(M_c(\Omega_N)\setminus \{\it o\})$, we consider the  class of functions
 $D(\mathfrak{G}^{\Omega_N})=\{H(\nu)    =h\left(  \int
 F(\mu)\nu(d\mu)\right)\}$ with $F(\mu)=f(\langle\mu,\varphi\rangle)$
 where $\varphi$ has finite support and $h,f$ are continuous functions on
 $\mathbb{R}$ with bounded second derivatives. We also define
 \[\frac{\delta
F(\mu)}{\delta\mu(x)} = \frac{d}{d\varepsilon}F(\mu+\varepsilon
\delta_x)\big|_{\varepsilon =0},\quad
\frac{\delta^{2}F(\mu)}{\delta\mu(x)\delta\mu(y)}=\frac {\partial
^2}{\partial \varepsilon_1 \partial \varepsilon_2}
F(\mu+\varepsilon_1\delta_x +\varepsilon_2
\delta_y)\big|_{\varepsilon_1=\varepsilon_2 =0}.\]

Then $\Psi^{(N)}$ is the unique solution to the martingale problem
given by the generator

\begin{eqnarray}\label{spatialMP}
&&\mathfrak{G}^{\Omega_N}H(\nu)  \\&&  =
\sum_{x\in\Omega_{N}}\int\nu(d\mu)h^{\prime}(\int F(\mu)\nu
(d\mu))\nonumber\\&&\qquad\qquad\cdot \left[  \frac{\delta
F(\mu)}{\delta\mu(x)}\left(
\sum_{k=1}^{\infty}\frac{c_{k+1}}{N^{k/2}}(\bar{\mu}_{k}(x)-\mu(x))\right)
+\frac{\delta^{2}F(\mu)}{\delta\mu(x)\delta\mu(x)}%
\mu(x)\right]\nonumber\\&&    +\int h^{\prime\prime}(\int
F(\mu)\nu(d\mu))F^{2}(\mu)\nu(d\mu)\nonumber
\end{eqnarray}
where
$$\bar\mu_k(x) = \frac1{N^k} \sum_{y\in \Omega_N:d_N(y,x)\le k}\mu(y).$$

Proposition \ref{4+dequil} establishes the existence of
a non-trivial equilibrium $\Psi^{(N,\theta)}$ that is spatially
homogeneous (that is, with law invariant under translations in
$\Omega_N$) and mean $\theta$. The
\textit{equilibrium random field} $\{\psi^{(N,
\theta)}_x(t)\}_{x\in \Omega_N}$ is defined
by
\begin{equation}\label{defpsi}
\psi^{(N, \theta)}_x(t) = \int \mu(x)
\Psi^{(N,\theta)}(t,d\mu).\end{equation}
This gives the total number of individuals at site $x$ irrespective of their family memberships.


Recall that $\{\eta^{(N,\theta)}_\ell(t,dx)\}$ defined in
(\ref{norm}) describes the equilibrium normalized family size
process in the block $B^{(N)}_\ell$  with mean number $\theta$
individuals per site.

 Also recall that the normalized equilibrium population mass
in $B_\ell^{(N)}$ is given by (see (\ref{zeta}),(\ref{zeta1})
\begin{eqnarray}\label{familysize}
\zeta_{\ell}^{(N,\theta)}(t)&&= \int
x\eta_{\ell}^{(N,\theta)}(t,dx)
=\frac{1}{N^\ell}\sum_{x\in B^{(N)}_\ell}\psi^{(N,
\theta)}_x(N^{\ell/2}t)\nonumber\\&& =\frac{1}{N^\ell}\sum_{x\in
B^{(N)}_\ell}\int\mu(x)\Psi^{(N,\theta)}(N^{\ell/2}t,d\mu).\end{eqnarray}

\begin{lemma} \label{greenbound}
Let $G_N$ denote the Green operator of the $(2,(c_j),N)$-random
walk with $(c_j)$ satisfying condition (\ref{trans2}) for strong transience, and write, for
$\varphi_1,\varphi_2:\Omega_N\to \mathbb{R}_+$, $$\langle \varphi_1, \varphi_2\rangle = \sum_{x \in \Omega_N} \varphi_1(x) \varphi_2(x).$$ Then
\begin{equation}\label{e1}
  \langle\varphi_1,G_N\varphi_2\rangle \leq \mathrm{const}
  \sum_{x,y\in\Omega_N}\frac{\varphi_1(x)\varphi_2(y)}{N^{|x-y|/2}}
\end{equation}
and
\begin{equation}\label{e2}
\langle\varphi_1,G^2_N\varphi_2\rangle \leq \mathrm{const}
  \sum_{x,y\in\Omega_N}\varphi_1(x)\varphi_2(y)\sum_{j=|x-y|}^{\infty}\frac{1}{c_j^2}
\end{equation}
where the constants do not depend on $N$. Hence
\begin{eqnarray}\label{e3}
&&\langle\varphi_1,G_N\varphi_2\rangle\leq \mathrm{const}\langle
1,\varphi_1\rangle\langle 1,\varphi_2\rangle,\\&&
\langle\varphi_1,G_N^2\varphi_2\rangle \leq \mathrm{const} \langle
1,\varphi_1\rangle\langle 1,\varphi_2\rangle.
\end{eqnarray}

\end{lemma}

\begin{proof}
The transition probability $p_t(x,y)$ of the $(2,(c_j),N)$-random
walk is given by
\begin{equation}
p_t(x,y)= (\delta_{0,|x-y|}-1) \frac{e^{-h^{(N)}_{|x-y|}t}
}{N^{|x-y|}} +(N-1) \sum_{j=|x-y|+1}^\infty
\frac{e^{-h_j^{(N)}t}}{N^j}
\end{equation} where the $h^{(N)}_j$ are positive
numbers (depending on $N$) such that

\[ h^{(N)}_j\geq \mathrm{const} \frac{c_{j-1}}{N^{(j-1)/2}},\quad j\geq 1\]
(see \cite {DGW2}). Hence
 \begin{eqnarray}
&&G_N\varphi(x)=\int_0^\infty \sum_y
p_t(x,y)\varphi(y)dt\nonumber\\&& \leq
(N-1)\sum_y\varphi(y)\sum_{j=|x-y|+1}^\infty
\frac{1}{N^jh^{(N)}_j}\nonumber\\&& \leq \mathrm{const}
(N-1)\sum_y\varphi(y)\sum_{j=|x-y|+1}^\infty
\frac{N^{(j-1)/2}}{N^jc_{j-1}}\nonumber\\&& \leq \mathrm{const} \sum_y
\varphi(y)\sum_{j=|x-y|}^\infty\frac{1}{N^{j/2}c_j}\nonumber\\&&
\leq \mathrm{const} \sum_y
\varphi(y)\sum_{j=|x-y|}^\infty\frac{1}{N^{j/2}}\nonumber\\&& \leq
\mathrm{const} \sum_y \varphi(y) \frac{1}{N^{|x-y|/2}},\nonumber
 \end{eqnarray}
 and then (\ref{e1}) follows.

Similarly (see \cite{DGW2})
 \begin{eqnarray}
&&G^2_N\varphi(x)=\int_0^\infty t\left(\sum_y
p_t(x,y)\varphi(y)\right)dt\nonumber\\&& \leq
(N-1)\sum_y\varphi(y)\sum_{j=|x-y|+1}^\infty
\frac{1}{N^j(h^{(N)}_j)^2}\nonumber\\&& \leq \mathrm{const}
(N-1)\sum_y\varphi(y)\sum_{j=|x-y|+1}^\infty
\frac{N^{j-1}}{N^jc^2_{j-1}}\nonumber\\&& \leq \mathrm{const} \sum_y
 \varphi(y) \sum_{j=|x-y|}^\infty\frac{1}{c_j^2}\nonumber
 \end{eqnarray}
and (\ref{e2}) follows.
\end{proof}
\begin{remark}
In \cite{DGW2} the exponential waiting time of the hierarchical random walk has parameter $1$, whereas here this parameter is $\sum_\ell q_\ell^{(N)}$, see subsection \ref{classofrw}. In the present case, $ \sum_\ell q_\ell^{(N)}= \sum_{j=0}^\infty
\frac{c_j}{N^{j/2}}=: L_N$. This amounts to a time change $t \to L_N t$ in the calculations above, which produces $L_N h_j^{(N)}$ in place of $h_j^{(N)}$ and does not change the results of Lemma \ref{greenbound} (assuming $c_0 > 0$).
\end{remark}

\ The following is the analogue of \cite{DGW1}(2.3.3).

\begin{proposition}\label{secmoment}
Under the conditions of Lemma \ref{greenbound}, the first and second moments of $\zeta^{(N,\theta)}_\ell(0)$ are
given by \[\mathbb E\zeta^{(N,\theta)}_\ell(0)=\theta\] and
\begin{equation}
\mathbb E\left(\zeta^{(N,\theta)}_\ell(0)\right)^2 =
\theta^2+\theta(\langle\varphi_{N,\ell},\varphi_{N,\ell}\rangle+\langle
 \varphi_{N,\ell},G_N\varphi_{N,\ell}\rangle +\frac{1}{4}
 \langle \varphi_{N,\ell},G_N^2\varphi_{N,\ell}\rangle)
 \end{equation}
where $\varphi_{N,\ell}(x)=\frac{\mathbf 1}{N^\ell}1_{B^{(N)}_\ell}(x)$
and  $G_N$ is  the Green operator of the (strongly transient)
$(2,(c_\ell),N)$-random walk.

\end{proposition}

\begin{corollary} \label{secmomentcor}If the $(2,(c_j),N)$-random walk satisfies conditions (\ref{cbound})
and (\ref{trans2}), then for each $\ell$ \newline (a)
\begin{equation}\label{spatial 2ndmoment}\sup_N E((\zeta^{(N,\theta)}_\ell(0))^2) <
\infty.
\end{equation}
 (b)
\begin{equation}\label{spatialvariance}
\mathrm{Var}(\zeta^{(N,\theta)}_\ell(0)) \to 0\;\; \mbox {as}\;
\ell\to\infty
\end{equation}
uniformly in $N$.\\

\end{corollary}
\begin{proof}
(a) follows immediately from Proposition \ref{secmoment} and Lemma \ref{greenbound}.\\
(b) From Proposition \ref{secmoment}
\begin{equation}\label{e5}
\mathrm{Var}(\zeta^{(N,\theta}_\ell(0))=\theta(\langle\varphi_{N,\ell},\varphi_{N,\ell}\rangle+\langle
 \varphi_{N,\ell},G_N\varphi_{N,\ell}\rangle +\frac{1}{4}
 \langle \varphi_{N,\ell},G_N^2\varphi_{N,\ell}\rangle)
\end{equation}
We will show that each of the three terms on the r.h.s. of
(\ref{e5}) converges to $0$ as $\ell\to\infty$ uniformly in $N$.
First $\langle\varphi_{N,\ell},\varphi_{N,\ell}\rangle \leq
\frac{1}{N^\ell}\to 0$ as $\ell\to\infty$ uniformly in $N$.

We have from (\ref{e1}) and using the ultrametric property of $|
\cdot|$,
\begin{eqnarray}
\langle \varphi_{N,\ell},G_N \varphi_{N,\ell} \rangle && \leq
\mathrm{const} \frac{1}{N^{2\ell}} \sum_{|x|,|y|\leq \ell}
\frac{1}{N^{|x-y|/2}}\nonumber    \\&& \leq \mathrm{const}
\frac{1}{N^{2\ell}}\left(N^\ell + \sum_{
  x\ne y,\,
  |x|\leq |y|\leq\ell}
\frac{1}{N^{|y|/2}} \right)\nonumber\\&&\leq \mathrm{const} \frac{1}{N^{2\ell}}\left(N^\ell
+\sum_{k=1}^\ell\frac{N^{2k}}{N^{k/2}}\right)\nonumber
\\&&= \mathrm{const}\frac{1}{N^{2\ell}}\left(
N^\ell+\frac{N^{(\ell+1)3/2}-N^{3/2}}{N^{3/2}-1}\right)\nonumber\\&&
\leq \mathrm{const} \left(\frac{1}{N^\ell}+\frac{1}{N^{\ell/2}}\right)
\longrightarrow 0 \;\; \mathrm{as}\; \ell\to\infty
\nonumber
\end{eqnarray}
uniformly in $N$, and from (\ref{e2}), again using the ultrametric property,
\begin{eqnarray}
\langle \varphi_{N,\ell},G^2_N \varphi_{N,\ell} \rangle && \leq
\mathrm{const} \frac{1}{N^{2\ell}} \sum_{|x|,|y|\leq \ell}
\sum_{j=|x-y|}^\infty \frac{1}{c_j^2}\nonumber
\\&& \leq \mathrm{const}
\frac{1}{N^{2\ell}}\left(N^\ell\sum_{j=0}^\infty\frac{1}{c_j^2}
+\sum_{
  x\ne y,\,
  |x|\leq |y|\leq\ell} \sum_{j=|y|}^\infty \frac{1}{c_j^2}\right)\nonumber\\&&
\leq \mathrm{const} \frac{1}{N^{2\ell}}\left(N^\ell+\sum_{k=1}^\ell
N^{2k}\sum_{j=k}^\infty\frac{1}{c_j^2}\right)\nonumber
\\&& \leq
\mathrm{const}\frac{1}{N^{2\ell}}
\left(N^\ell+\sum_{j=1}^\infty\frac{1}{c_j^2}\sum_{k=1}^{j\wedge\ell}N^{2k}\right)\nonumber
\\&&
\leq \mathrm{const}
\frac{1}{N^{2\ell}}\left(N^\ell+\sum_{j=1}^\infty \frac{1}{c^2_j}N
^{2(j\wedge\ell)}\right)\nonumber\\&&
=\mathrm{const}\left(\frac{1}{N^\ell}+\sum_{j=1}^{\ell-1}\frac{1}{c_j^2}\frac{1}{N^{2(\ell-j)}}
+\sum_{j=\ell}^\infty\frac{1}{c_j^2}\right).
\end{eqnarray}
The first term goes to zero as $\ell\to\infty$, the second term
goes to zero as $\ell\to\infty$ by dominated convergence, and
clearly the last term goes to zero as $\ell\to\infty$. Therefore
\begin{equation}
\langle \varphi_{N,\ell},G_N^2\varphi_{N,\ell}\rangle \to
0\;\;\mathrm{as}\;\; \ell\to\infty\;\;\;\mathrm{uniformly\;in\;}N.
\end{equation}
\end{proof}
\begin{remark}
(a) \begin{equation} \lim_{N\to\infty} \mathrm{Var}
(\zeta^{(N,\theta)}_\ell(0))= \theta \sum_{j=\ell}^\infty
\frac{1}{c_j^2}.
\end{equation}
(b) \begin{equation}
\langle\varphi_{N,\ell},G^2_N\varphi_{N,\ell}\rangle \leq
\mathrm{const}\left(\frac{1}{N^\ell}+\frac{1}{N^2}+\sum_{j=\ell}^\infty\frac{1}{c_j^2}\right)\end{equation}
\end{remark}

\begin{corollary}\label{sumc}
Consider the ``exterior function''
\begin{equation}\label{e7}
\varphi_{N,\ell,\mathrm{ext}}(x)=\sum_{k=1}^\infty
\frac{c_{\ell+k-1}}{N^{\ell+2k-1}}{\mathbf 1}_{B_{\ell+k}^{(N)}}(x).
\end{equation}
Then\begin{eqnarray}   &&  \mathbb E\left[\left(\sum_{x\in \Omega_N}
\varphi_{N,\ell,\mathrm{ext}}(x) \int
\mu(x)\Psi^{(N,\theta)}(0,d\mu)\right)^2\right]\nonumber\\&& =
\mathbb E\left[\left(\sum_{k=1}^{\infty}
\frac{c_{\ell+k-1}}{N^{k-1}}\zeta^{(N,\theta)}_{\ell+k}(0)\right)^2\right]<\infty
\end{eqnarray} uniformly in $N$.
\end{corollary}
\begin{proof}
As in Proposition \ref{secmoment}, (from \cite{DGW1}, comment 2.3.5)
\begin{eqnarray}\label{comment2.3} &&\mathbb E \left[\left(\sum_{x\in \Omega_N} \varphi_{N,\ell,\mathrm{ext}}(x) \int
\mu(x)\Psi^{(N,\theta)}(0,d\mu)\right)^2\right]\\ \nonumber&& =
\theta^2\langle 1,\varphi_{N,\ell,\mathrm{ext}}\rangle^2
 +\theta(\langle\varphi_{N,\ell,\mathrm{ext}},\varphi_{N,\ell,\mathrm{ext}}\rangle\\ \nonumber &&+\langle
 \varphi_{N,\ell,\mathrm{ext}},G_N\varphi_{N,\ell,\mathrm{ext}}\rangle +\frac{1}{2}
 \langle
 \varphi_{N,\ell,\mathrm{ext}},G_N^2\varphi_{N,\ell,\mathrm{ext}}\rangle).\nonumber
 \end{eqnarray}
We will prove that each term on the r.h.s. of (\ref{comment2.3})
is bounded uniformly in $N$.
\begin{eqnarray}
\langle
1,\varphi_{N,\ell,\mathrm{ext}}\rangle&&=\sum_{k=1}^\infty\frac{c_{\ell+k-1}}{N^{\ell+2k-1}}\nonumber\\&&
+\sum_{k=1}^\infty\frac{c_{\ell+k-1}}{N^{\ell+2k-1}}\sum_{n=1}^{\ell+k}(N-1)N^{n-1}\nonumber\\&&
\leq
\sum_{k=1}^\infty\frac{c_{\ell+k-1}}{N^{\ell+2k-1}}+\sum_{k=1}^\infty\frac{c_{\ell+k-1}}
{N^{\ell+2k-1}}N^{\ell+k}\nonumber\\&& \leq \mathrm{const}
\sum_{k=1}^\infty\frac{c_{\ell+k-1}}{N^k}\nonumber\\&&
<\infty\quad\mathrm{uniformly\;\;in\;\;}N.\nonumber
\end{eqnarray}
Then, from (\ref{e3}) the terms
$\langle\varphi_{N,\ell,\mathrm{ext}},G_N\varphi_{N,\ell,\mathrm{ext}}\rangle$
and
$\langle\varphi_{N,\ell,\mathrm{ext}},G_N^2\varphi_{N,\ell,\mathrm{ext}}\rangle$
are bounded uniformly in $N$ and since
$\varphi_{N,\ell,\mathrm{ext}}$ has a uniform bound in $N$ as
well, using (\ref{cbound}) we see that
$\langle\varphi_{N,\ell,\mathrm{ext}},\varphi_{N,\ell,\mathrm{ext}}\rangle$
is also bounded uniformly in $N$ (with $N >
\sup c_{\ell+1}/c_\ell$).
\end{proof}

\begin{theorem} {\bf (Spatial ergodic theorem)} \label{spergthm}
The pointwise ergodic theorem holds on $\Omega_N$, that is,
\begin{equation}\label{ergodic} \lim_{j\to\infty}
\zeta_{j}^{(N,\theta)}(0) = \theta \mbox { a.s. }
\end{equation}
\end{theorem}
\begin{proof}
First note that the equilibrium random field $\{\psi^{(N,
\theta)}_x\}_{x\in \Omega_N}$  defined by  (\ref{defpsi}) is invariant under the action of
the group $\Omega_N$, and $E(\psi^{(N, \theta}_x) =\theta$.
Moreover by (\ref{spatialvariance}) the spatial averages
$\zeta^{(N,\theta)}_\ell(0) =\frac{1}{N^\ell} \sum_{x\in
B^{(N)}_\ell} \psi^{(N, \theta)}_x(0)$ satisfy
\[\lim_{j\to\infty} \mathrm{Var}(\zeta_{j}^{(N,\theta)}(0))=0,\]
so that the convergence in probability follows. To complete the
proof note that $\Omega_N$ is an amenable group and the collection
of balls $\{B^{(N)}_\ell\}_{\ell\in\mathbb{N}}$ is a tempered
F{\o}lner sequence.  The a.s. pointwise convergence then follows
by \cite{Lin}(Theorem 1.2).
\end{proof}

\subsection{ Reduction to two successive scales}\label{sec5.2}
In the section we show that the analysis of the multiscale
behavior can be reduced to the case of two successive scales.

\subsubsection{The equilibrium family size process}
\setcounter{equation}{0}
     Let $\{\zeta_{\ell}^{(N,\theta)}(t),\; t\in \mathbb{R}\}_{\ell =\dots,2,1}$ denote the
collection of equilibrium block average  processes in the nested
sequence of blocks $B^{(N)}_\ell$ (see (\ref{familysize})).

\begin{lemma} \label{familysizeconstraint} The fixed time marginal distributions $\{\zeta_{\ell}^{(N,\theta)}(0)\}$
at time $t=0$ satisfy:
\\
(a) For each $\ell\in \mathbb{N}$ the family
$\{\zeta^{(N,\theta)}_\ell(0)\}_{N\in\mathbb{N}}$ is tight and
every limit point has expected value $\theta$.\\
 (b) The
family-size constraint
\begin{equation}\label{fsc}
\mathbb E\left[ \int_K^\infty x\eta^{(N,\theta)}_\ell(0,dx)\right]\leq
\frac{\mathrm{const}}{K}
\end{equation}
 is satisfied for a suitable constant not depending on $K > 0$ and  $N \ge 2$.
\end{lemma}
\begin{proof} (a) follows immediately from Corollary
\ref{secmomentcor}.\\
(b) Let $\varphi^N_\ell(x)=\frac{1}{N^\ell}1_{B^{(N)}_\ell}(x)$.
Then
\begin{eqnarray}&&
 \mathbb E\left[\int_K^\infty x\eta^{(N,\theta)}_\ell(0,dx)\right] \leq
\frac{1}{K}  \mathbb E\left[ \int_0^\infty
x^2\eta^{(N,\theta)}_\ell(0,dx)\right]\nonumber\\&& \leq
\frac{1}{K}  \mathbb E\left[ \int_0^\infty
x^2\eta^{(N,\theta)}_\ell(0,dx)\right]\nonumber\\&&
=\frac{1}{K}\frac{1}{N^{3\ell/2}} \mathbb E\left[\int(\mu(B^{(N)}_\ell))^2\Psi^{(N,\theta)}(0,d\mu)\right]\nonumber
\\&&
=\frac{1}{K}N^{\ell/2} \mathbb
E\left[\int(\mu(\varphi^N_\ell))^2\Psi^{(N,\theta)}(0,d\mu)\right]\quad
 \nonumber
\\&& =
\frac{N^{\ell/2}}{K}
\left[\langle\varphi^N_\ell,\varphi^N_\ell\rangle
+\frac{1}{2}\langle\varphi^N_\ell,G_N\varphi^N_\ell\rangle\right]\quad
(\mathrm{by}\;\; \cite{DGW1} (2.3.1))\nonumber\\&&
\leq\frac{N^{\ell/2}}{K}\mathrm{const}\left(\frac{1}{N^\ell}+\frac{1}{N^{\ell/2}}\right)
\quad \mathrm{ (by\;the\;proof\;of\;Corollary
\;}\ref{secmomentcor})\nonumber\\&& \leq\frac{\mathrm{const}}{K}.\nonumber
\end{eqnarray}

\end{proof}

\begin{remark}
Part (b) of the Lemma \ref{familysizeconstraint} implies that
asymptotically the restriction of $\Psi^{(N,\theta)}(0)$ to the
ball $B^{(N)}_\ell$ consists of families of size $O(N^{\ell/2})$
or smaller.
\end{remark}

\subsubsection {Distant immigrants} \setcounter{equation}{0}

In order to establish that there is a unique limit law and to
identify it we now return to the dynamical picture. In the next
lemma we show that the expected contribution to the equilibrium
population in a ball $B^{(N)}_\ell$  coming from immigrants who
immigrate directly from outside the ball $B^{(N)}_{\ell +1}$ and
the descendants of the population in the ball $B^{(N)}_\ell$ in
the distant past  are both negligible.

\begin{lemma} \label{exteriorcontribution}
(a) Let $\{\zeta_{\ell,\mathrm{ext}}^{(N,\theta)}(0)\}$ denote
the contribution to the equilibrium block average in the ball
$B^{(N)}_\ell$ at time $0$ coming from individuals immigrating
directly from $(B^{(N)}_{\ell+1})^c$. Then
\[  E\left(\zeta_{\ell,\mathrm{ext}}^{(N,\theta)}(0)\right) \leq \frac{\mathrm{const}}{N^{\frac{1}{2}}}.\]
\\
(b) The expected mass to enter $B^{(N)}_{\ell+1}$ from
$(B^{(N)}_{\ell+1})^c$  in a time interval of length $N^{\ell/2}$
is $O(\frac{1}{N^{1/2}})$.
 \\ (c) The expected contribution to the population in the ball
$B^{(N)}_\ell$ at time $0$ from  the descendants of individuals
alive at time $t_0$ is of order $O(e^{-c_\ell |t_0|})$ as $t_0\to
- \infty$.
\end{lemma}
\begin{proof}
(a) To verify this, we first note that  the total number of
individuals to immigrate from $B^{(N)}_{\ell+k}$ with $k\geq 2$ to
$B^{(N)}_{\ell}$ in the time interval
$[-N^{\ell/2}t,-N^{\ell/2}(t+dt))$ is of the order
\[(\zeta^{(N,\theta)}_{\ell+k}(-t)N^{\ell+k})\times
(\frac{c_{\ell+k-1}}{N^{(\ell+k-1)/2}}) \times \frac{1}{
N^k}\times dt N^{\ell/2}
\]
where the first factor is the number of particles in
$B^{(N)}_{\ell+k}$, the second factor is the rate of migration of
each of these particles to a point chosen randomly in
$B^{(N)}_{\ell+k}$, the third factor is the probability this point
falls in the tagged ball $B^{(N)}_\ell$ and the last factor is the
length of the time period. Recalling that the mass process in
$B^{(N)}_\ell$ is subcritical with parameter
$\frac{c_\ell}{N^{\ell/2}}$,  the expected total mass at time $0$
coming from  immigration from outside $B^{(N)}_{\ell+1}$ in the
time interval $[-N^{\ell/2}T,0)$, $T>0$, is of the order
\begin{eqnarray}\label{farimmigrants}
&&\sum_{k\geq 2} \int_{-T}^0
\frac{c_{\ell+k-1}}{N^{(k-1)/2}}\times N^\ell\times
 e^{c_\ell t} \mathbb E(\zeta^{(N,\theta)}_{\ell+k}(-t))dt\\ &&
\leq \textrm{const}\times \frac {N^\ell}{N^{1/2}},\nonumber
\end{eqnarray}
uniformly in $T$ and $\ell$, where we have used the assumption (\ref{cbound}). Therefore the
expected mass (normalized by $\frac{1}{N^\ell}$) at time $0$  in
$B^{(N)}_\ell$ that immigrated during the time interval
$[N^{\ell/2}t_0,0)$ directly from outside $B^{(N)}_{\ell+1}$ is
$O(N^{-\frac{1}{2}})$ as $N\to\infty$.

(b) and (c)   follow from a first moment calculation.

\end{proof}

\medskip
\subsection{Diffusion limit of the  family size process in two
spatial scales}\label{sec5.3}
 \setcounter{equation}{0}

In this section we consider the asymptotic (as $N\to\infty$) time
development
 of the population occupying
$B^{(N)}_{\ell+1}$ and in particular the subpopulation obtained by
considering individuals that occupy a tagged ball $B^{(N)}_\ell$
in the natural time scale for the population in $B^{(N)}_\ell$ and
assuming that the initial family size processes in
$B^{(N)}_{\ell+1}$ satisfy the  family size constraint of Lemma
\ref{familysizeconstraint}.

 As a result of the previous  lemma, asymptotically as
 $N\to\infty$ and $t_0\to -\infty$
 the equilibrium population
 in  $B^{(N)}_\ell$ consists of the descendants of immigrants
 coming from  $B^{(N)}_{\ell+1}$ during the time interval
 $(-N^{-\ell/2}t_0,0]$ (recall Remark \ref{timescale}).
Another key point is that, as we verify below, in this time
interval the family size process in $B^{(N)}_{\ell+1}$ is
asymptotically constant.  Moreover from Lemma
\ref{familysizeconstraint} all but an arbitrarily small proportion
of  the population in $B^{(N)}_{\ell+1}$ is structured into
$O(N^{(\ell+1)/2})$ families containing $O(N^{(\ell+1)/2})$
individuals.

We next show that in the $N^{\ell/2}$-time scale the total
population structure in the ball $B^{(N)}_{\ell+1}$ is essentially
constant.

\begin{lemma}\label{constant} Let
     $\{\zeta_{\ell+1}^{(N,\theta)}(s)\}$ be the
     equilibrium normalized process (see (\ref{familysize})) in its natural time scale $N^{(\ell+1)/2}$.
\newline
 For $t_0= t_0(N)<0$ such that $\frac{|t_0|}{N^{\frac{1}{2}}}< c$,
\begin{eqnarray}\label{martinc2} && P\left(\sup_{t_0\leq t\leq 0} |
\zeta_{\ell+1}^{(N,\theta)}
(N^{-1/2}t)-\zeta_{\ell+1}^{(N,\theta)}(N^{-1/2}t_0)|> K\right)
\nonumber\\&&\leq \left(\frac{\mathrm{const}}{K^2}\right)c + o(1),
\end{eqnarray}
where const does not depend on $t_0$ and $N$, and $o(1)$ converges
to $0$ as $N \to \infty$.

\end{lemma}
\begin{proof}
Recall that   $\Psi^{(N)}(t)$ is characterized as the unique
solution of the martingale problem  with generator
(\ref{spatialMP}). Applying the generator to the function
$F(\nu)=\int
\langle\mu,\frac{1}{N^{\ell+1}}1_{B^{(N)}_{\ell+1}}\rangle\nu(d\mu)$,
it follows that $\{M_{\ell+1}(t)\}_{t_0\leq t\leq 0}$, defined by
\begin{eqnarray} M_{\ell+1}(t)=&&\zeta_{\ell+1}^{(N,\theta)}
(t)-\zeta_{\ell+1}^{(N,\theta)} (t_0)\nonumber\\&&-\int_{t_0}^t
\sum_{k=1}^{\infty}\frac{c_{\ell +k
}}{N^{k-1}}\left(\zeta^{(N,\theta )}_{\ell+1+k}(N^{-k/2}s)-
\zeta^{(N,\theta)}_{\ell+1}(s)\right)ds\nonumber\end{eqnarray} is
a martingale.

Then
\begin{eqnarray}
&&P \left(\sup_{t_0\leq t\leq 0} | \zeta_{\ell+1}^{(N,\theta)}
(N^{-1/2}t)-\zeta_{\ell+1}^{(N,\theta)}(N^{-1/2}t_0)|>
K\right)\nonumber\\&&\leq P\left(   \int_{\frac{t_0}{N^{1/2}}}^0
\big|\sum_{k=1}^{\infty}\frac{c_{\ell +k
}}{N^{k-1}}\left(\zeta^{(N,\theta)}_{\ell+1+k}(N^{-k/2}s)-
\zeta^{(N,\theta)}_{\ell+1}(s)\right)\big|ds
>\frac{K}{2})\right) \nonumber\\&& + P\left(\sup_{\frac{t_0}{N^{1/2}}\leq t\leq 0}
|M_{\ell+1}(t)|>\frac{K}{2} \right)
\end{eqnarray}

Recall from Corollary \ref{secmomentcor}  that $
\mathbb E\left(\zeta^{(N,\theta)}_\ell(t)\right)^2 $ is bounded uniformly
in $N$ and $t$. Let \[g_N(s):=|\sum_{k=1}^{\infty}\frac{c_{\ell +k
}}{N^{k-1}}\left(\zeta^{(N,\theta)}_{\ell+1+k}(N^{-k/2}s)-
\zeta^{(N,\theta)}_{\ell+1}(s)\right)|\] and note that by
Corollary \ref{sumc} $ \mathbb E(g_N^2(s))$ is uniformly bounded in $s$
and $N$. Then
\begin{eqnarray}
&&P\left(   \int_{\frac{t_0}{N^{1/2}}}^0 g_N(s)ds
>\frac{K}{2}\right)\nonumber\\&&
\leq \frac{\mathrm{const}}{K^2}(\frac{t_0}{N^{1/2}})^2 \mathbb E\left[
\frac{1}{\frac{t_0}{N^{1/2}}}\int^0_{\frac{t_0}{N^{1/2}} }
g_N(s)ds\right]^2\nonumber\\&&\leq \frac{\mathrm{const}}{K^2}
(\frac{t_0}{N^{1/2}})^2 \left[
\frac{1}{\frac{t_0}{N^{1/2}}}\int^0_{\frac{t_0}{N^{1/2}}}
 \mathbb E(g_N(s)^2)ds\right]\nonumber\\&&\leq \frac{\mathrm{const}}{K^2}
\left(\frac{t_0}{N^{1/2}}\right)^2\nonumber
\end{eqnarray}

Next we note that by Lemma \ref{moments}(b) and Remark \ref{inhimm} (with
$\varepsilon =\frac{1}{N}$) we get for $t_0\leq t\leq 0$
\begin{eqnarray*} &&  \mathbb E\Big( [\zeta_{\ell+1}^{(N,\theta)}
 (t)-\zeta_{\ell+1}^{(N,\theta)}
 (t_0)]^2 \Big|\zeta_{\ell+1}^{(N,\theta)}
 (t_0)=m, \sum_{k=1}^{\infty}
\frac{c_{\ell+k}}{N^{k-1}}\zeta^{(N,\theta)}_{\ell+1+k}(\cdot)=a(\cdot)
\Big)\\&& =(m^2e^{-2c_{\ell+1}(t-t_0)}-m^2)\nonumber \\&&
+\frac{m}{c_{\ell+1}^2}
\{c_{\ell+1}(t-t_0)e^{-c_{\ell+1}(t-t_0)}+e^{-2c_{\ell+1}(t-t_0)}+2(t-t_0)c_{\ell+1}^3ae^{-c_{\ell+1}(t-t_0)}
\\ &&\qquad\qquad -e^{-c_{\ell+1}(t-t_0)} +(t-t_0)c_{\ell+1}^2\varepsilon
e^{-c_{\ell+1}(t-t_0)}\} \\&& +\frac{1}{c_{\ell+1}^2}\int_{t_0}^t
k_1(s,t)a(s)ds\quad + \int_{t_0}^t\int_{t_0}^{s_2} k_2(t,s_2,s_1)
a(s_1)a(s_2)ds_1ds_2\\
&& + o(1)\cdot\int_{t_0}^t k_3(t,s)s(s)ds
\end{eqnarray*} where $k_i(t,\cdot)\;i=1,3$ and
$k_2(t,\cdot,\cdot)$ are bounded non-negative kernels that satisfy
(\ref{kconditions}) and $o(1)$ converges to $0$ as $N \to \infty$.
Using Corollary \ref{secmomentcor} we conclude that
\begin{equation}\label{secondmoment}
 \mathbb E\left(
[\zeta_{\ell+1}^{(N,\theta)}
 (N^{-1/2}t)-\zeta_{\ell}^{(N,\theta)}
 (N^{-1/2}t_0)]^2\right) \leq {\mathrm {const}}\cdot
\frac{(t-t_0)}{N^{1/2}} +o_1(\frac{1}{N}).
\end{equation}

 We now apply the $L^2$-martingale
inequality, \begin{eqnarray}&&
P\left(\sup_{\frac{t_0}{N^{1/2}}\leq t\leq 0}
|M_{\ell+1}(t)|>\frac{K}{2}\Big|M_{\ell+1}(t_0/N^{1/2})=0
\right)\nonumber\\&&\leq \frac{1}{K^2}
 \mathbb E\left(|M_{\ell+1}(0)|^2\Big|M_{\ell+1}(t_0/N^{1/2})=0\right)\nonumber\\&&\leq
\frac{\mathrm{const}}{K^2}\frac{|t_0|}{N^{1/2}}+o_1(\frac{1}{N}).\nonumber
\end{eqnarray}

\end{proof}

\bigskip
Let $M^1[[0,\infty)\times(0,\infty)]\backslash \{0\}$ denote the
set of non-zero Radon measures on $[0,\infty)\times (0,\infty)$
such that $\int_{0+}^\infty\int_0^\infty
y\mu(dx,dy)+\int_{0+}^\infty y\mu(\{0\},dy) <\infty $ and
$\int_{0+}^\infty\int_{0}^\infty x\mu(dx,dy)<\infty $.

\begin{proposition}\label{keyproposition}
Let $\ell$ be fixed and
\[ \mathcal{Y}^{(N)}(t;dx,dy)=\frac{1}{N^{(\ell+1)/2}}\int \mathbf{1}_{\{\frac{\nu(B^{(N)}_{\ell})}{N^{\ell/2}}\in dx,
\frac{\nu(B^{(N)}_{\ell+1})}{N^{(\ell+1)/2}}\in dy
\}}\Psi^{(N,\theta)}(N^{\ell/2}t,d\nu)\] Assume that
\[{\mathcal{Y}^{(N)}(t_0;dx,dy)}  \stackrel{N\to\infty}{\Longrightarrow}  \;  \mu^0(dy)\delta_0(dx).\] Then
$\{N^{1/2}\mathcal{Y}^{(N)}(t;dx,dy)\mathbf{1}_{x>0}\}_{t_0<t\leq
0}$ converges weakly as $N\to\infty$ to the \break
$M^1(0,\infty)$-valued diffusion $\{\eta(t)\}_{t_0<t\leq 0}$ with
generator $\mathfrak{G}_\ell$ defined as follows.  Let
$F(\mu_\ell)=f(\langle\mu_\ell,\varphi\rangle)$ with $\varphi\in
C^2((0,\infty))$ with $|\varphi(x)|\leq const|x\wedge 1|$. Then
\begin{eqnarray}
\mathfrak{ G}_\ell F(\mu_\ell)&= f'(\langle
\mu_\ell,\varphi\rangle)\langle \mu_\ell, G_{c_\ell}^{(2)}
\varphi_\ell\rangle +\frac{1}{2}f''(\langle
\mu_\ell,\varphi_\ell\rangle)\langle
\mu_\ell,\varphi_\ell^2\rangle\\& +f'(\langle
\mu_\ell,\varphi_\ell\rangle)c_\ell\int\frac{\partial\varphi_\ell}{\partial
x}(x,y)|_{x=0}y\mu^0(dy)\nonumber.
\end{eqnarray}
where $G_{c_\ell}^{(2)}$ denotes the application of the operator
$G_{c_\ell}$ to the $x$ variable.  That is, $\eta(t)$ is a
two-level branching process with constant multitype immigration
source with immigration from zero of type $y$ at rate given by
$y\mu^0(dy)$ and total immigration rate of $\int_0^\infty
y\mu^0(dy)$ .
\end{proposition}

\begin{proof} The proof is a refinement of the proof of
Proposition \ref{immprocessconvergence}. We begin  by noting that
Lemma \ref{exteriorcontribution} implies that asymptotically as
$N\to\infty$ the contribution of  immigrants into
$B^{(N)}_{\ell+1}$ from $B^{(N)}_{\ell +k},\;k\geq 2$ is
negligible. Moreover  the contribution of immigrants into
$B^{(N)}_{\ell+1}$ from $B^{(N)}_{\ell+2}$ in the time scale
$N^{\ell/2}$ is also asymptotically negligible. Therefore in the
time scale $N^{\ell/2}$ we can restrict attention to the
population in $B^{(N)}_{\ell+1}$.  More precisely,  asymptotically
as $N\to\infty$  the population
$\zeta_{\ell}^{(N,\theta)}(t)=\int\int
xN^{1/2}\mathcal{Y}^{(N)}(t,dx,dy)\mathbf{1}_{x>0}$ consists of
descendants of immigrants entering from $B^{(N)}_{\ell+1}$ in the
time interval $[N^{\ell/2}t_0(N),0]$ provided that we take
$t_0(N)\to -\infty$.

As in the proof of Proposition \ref{immprocessconvergence} a
standard argument yields the tightness of the processes
$\{\zeta^{(N,\theta)}_\ell(t)\}_{t_0\leq t\leq 0}$ and
$\{\zeta^{(N,\theta)}_{\ell+1}(N^{-1/2}t)\}_{t_0\leq t\leq 0}$.
 The
tightness in $C([t_0,0],M^1((0,\infty)))$  of
$N^{1/2}\mathcal{Y}^{(N)}(\cdot)\mathbf{1}_{x>0}$ is also obtained
as in the proof of Proposition \ref{immprocessconvergence}. One
difference is the presence of additional terms in the expressions
for the moments. However these expressions  tend to zero as
$N\to\infty$.

 A first moment
calculation  shows that the expected contribution of individuals
who leave $B^{(N)}_\ell$ and then re-immigrate is $O(\frac{1}{N})$
and therefore asymptotically negligible. This means that we can
treat the population in $B^{(N)}_\ell$ in the time interval
$[t_0,0]$ as a two level branching system critical at the family
level  and subcritical at the individual level with subcriticality
parameter $c_\ell$ and with immigration of individuals from
$B^{(N)}_{\ell+1}$. Moreover by Lemma \ref{constant} the source of
immigrant individuals $\zeta^{(N,\theta)}_{\ell+1}(t)$ is
asymptotically constant in the time interval $(t_0(N),0)$ provided
that
\[\frac{|t_0(N)|}{N^{\frac{1}{2}}}\to 0\] as $N\to\infty$.

The main difference from the proof of Proposition
\ref{immprocessconvergence} is that we cannot assume that each
immigrant belongs to a different family.  Since the individuals
immigrating into $B^{(N)}_{\ell}$ come from families in
$B^{(N)}_{\ell+1}$ subject to family branching, it is necessary to
keep track of the family structure in the ball $B^{(N)}_{\ell+1}$.
The reason for this is that  in principle an individual
immigrating into $B^{(N)}_\ell$ could be a member of a
$B^{(N)}_{\ell+1}$-family already represented in $B^{(N)}_\ell$
and then could not be viewed as the founder of an independent
family in $B^{(N)}_\ell$. Part of the argument below is to verify
that that this effect is negligible.

By the family size constraint  at time $t_0$ (see Lemma
\ref{familysizeconstraint}(b)), the population in
$B^{(N)}_{\ell+1}$ consists of a collection of families whose
sizes are $O(N^{(\ell+1)/2})$. We index these families at time
$t_0$ by $i\in \mathbb{N}$ with masses $y_i N^{(\ell+1)/2}$.

Recall from section \ref{tlbe} that $\Psi^{(N,\theta)}(t) =\sum_j
\psi_j(t,\cdot)$ where    $\psi_j(t,\cdot)$ is a counting measure
on $\Omega_N$ corresponding to the spatial distribution of the
family indexed by $j$ at time $t$. We now give a precise
formulation to the time development of the families simultaneously
in $B^{(N)}_{\ell+1}$ and $B^{(N)}_\ell$. To each family $\psi_j$
and $t\in\mathbb{R}$ we associate a couple $(x_j(t),y_j(t))$ where
\[x_j(t)=\frac {\psi_j(t,B_{\ell}^{(N)})}{N^{\ell/2}}, \quad y_j(t)
=\frac {\psi_j(t,B_{\ell+1}^{(N)})}{N^{(\ell+1)/2}}.\]

Then
$$\mathcal{Y}^{(N)}(t;dx,dy) = \varepsilon_2\sum_i\delta_{(x_i(N^{\ell/2}t),y_i(N^{\ell/2}t))}(dx,dy),$$
where $\varepsilon_2 = N^{-(\ell+1)/2}$. Also, let
$\varepsilon_0=\frac{1}{N}$, $\varepsilon_1=\frac{1}{N^{\ell/2}}$.

By assumption,  the family size distribution in $B^{(N)}_{\ell+1}$
at time $t_0$, asymptotically as $N\to \infty$, is given by
$\mu^0$ with $\int y\mu^0(dy) <\infty$.

It suffices to show that for $t\geq t_0$, as $N\to\infty$
\begin{eqnarray}\label{two-1}&& \mathcal{Y}^{(N)}(t; dx,dy)\bold{1}_{\{x<
N^{-\frac{1}{4}}\}}\Rightarrow \mu^0_{\ell+1}(dy)\delta_0(dx)\\
\label{two-2}&& \mathcal{Y}^{(N)}(t; dx,dy)\bold{1}_{\{x>
0\}}\Rightarrow 0\\ \label{two-3}&&
\{\widetilde{\mathcal{Y}}^{(N)}(t)\}_{ t \ge t_0}\Rightarrow
\{\eta(t)\}_{t \ge t_0}
\end{eqnarray} where  $\widetilde{\mathcal{Y}}^{(N)}(t)$
is the renormalized  family measure given by
\begin{eqnarray*}\widetilde{\mathcal{Y}}^{(N)}(t;dx,dy)&=&N^{1/2}{\mathcal{Y}}^{(N)}
(t; dx,dy)\mathbf{1}_{\{x>0\}} +{\mathcal{Y}}^{(N)}
(t;dx,dy)\mathbf{1}_{\{x< N^{-\frac{1}{4}}\}}\end{eqnarray*} and
where $\eta(t)$ is a $M^1[[0,\infty)\times(0,\infty)]$-valued
diffusion with generator $\mathfrak G_\ell$.

To verify  (\ref{two-1}), first using
 Lemma \ref{exteriorcontribution} and  Proposition \ref{immprocessconvergence}
we can verify that
\[ \sup_{t\in (-t_0,0]}
\int\int (1\wedge y)|\mathcal{Y}^{(N)}(t; dx,dy)-
\mathcal{Y}^{(N)}(t_0; dx,dy)| \to 0 \mathrm{\;\;in\;\;
probability.}
\]
Moreover
\[
\int\int_{N^{-1/4}}^\infty N^{1/2} \mathcal{Y}^{(N)}(t; dx,dy) <
N^{1/4}\int\int_0^\infty x N^{1/2}\mathcal{Y}^{(N)}(t; dx,dy).
\]
Therefore by the tightness of $\int\int_0^\infty
 x(N^{1/2}\mathcal{Y}^{(N)}(\cdot; dx,dy)\bold{1}_{\{ x>0\}})$ and $\{\zeta^{(N,\theta)}_{\ell+1}\}$, as
$N\to\infty$,
\[\int\int_{N^{-1/4}}^\infty \mathcal{Y}^{(N)}(1\wedge y)(t;
dx,dy) \to 0 \mathrm{\;\;in\;\; probability} \] and (\ref{two-1})
and (\ref{two-2}) follow.

 We now turn to
the proof of  (\ref{two-3}). In order to implement the rescaling
we introduce the class of functions of the form
\[ F(\mu)=f(\langle \mu,\varphi))= f(\int\int\varphi(x,y)\mu(dx,dy))=
f(\varepsilon_2\sum_i \varphi(x_i,y_i)),\] where
\[ \varphi(x,y)= \varepsilon_0^{-1/2}\bold{1}_{\{x>0\}} \varphi_\ell(x,y)+
\varphi_{\ell+1}(y)\bold{1}_{\{x < N^{-\frac{1}{8}}\}}.\] We
assume that $|\varphi_\ell|,|\frac{\partial
\varphi_\ell(x,y)}{\partial x} |,|\varphi_{\ell+1}|$ are bounded,
$C^2$, and $|\varphi_\ell(x,y)|\leq \mathrm{const}\cdot x $ for
$x\geq 0$ and $|\varphi_{\ell+1}(y)|\leq const\cdot y$ for $y>0$.

Now define \[\mu_{\ell}(dx,dy)= N^{1/2}\mu(dx,dy)\mathbf{1}_{\{x>0\}}\] and
note that
\[
\int \varphi(x,y)\mu(dx,dy)=\int
\varphi_\ell(x,y)\mu_{\ell}(dx,dy)+\int
\varphi_{\ell+1}(x,y)\mathbf{1}_{\{x <
N^{-\frac{1}{8}}\}}\mu(dx,dy).
\]

The generator of $\mathcal{Y}^{(N)}$ acting on $F$ is given by
\newpage
 \small

\begin{eqnarray*}\nonumber
&&\mathfrak{G}^{(N)} F(\mu)\\&=& \nonumber
c_\ell\varepsilon_0\sum_{j_1=0}^{N^{\frac{2\ell-1}{4}}}\sum_{j_2>
j_1}
[f(\langle\mu,\varphi\rangle+\varepsilon_1(\varphi_\ell((j_1+1)\varepsilon_1,j_2\varepsilon_2)
-\varphi_\ell(j_1\varepsilon_1,j_2\varepsilon_2))
-f(\langle\mu,\varphi\rangle)]\\&&\hspace{2cm}\cdot\frac{(j_2-j_1)\varepsilon_2\mu(j_1\varepsilon_1,j_2\varepsilon_2)}
{\varepsilon_2^2}\nonumber
\\ \nonumber &+&
c_\ell\varepsilon_0\sum_{j_1>{N^{\frac{2\ell-1}{4}}}}\sum_{j_2>
j_1}
[f(\langle\mu,\varphi\rangle+\varepsilon_1(\varphi_\ell((j_1+1)\varepsilon_1,j_2\varepsilon_2)
-\varphi_\ell(j_1\varepsilon_1,j_2\varepsilon_2))
-f(\langle\mu,\varphi\rangle)]\\&&\hspace{2cm}\cdot\frac{(j_2-j_1)\varepsilon_2\mu(j_1\varepsilon_1,j_2\varepsilon_2)}
{\varepsilon_2^2}\nonumber\\ &+&\frac12
\sum_{j_1=0}^\infty\sum_{j_2\ge j_1}^\infty
[f(\langle\mu,\varphi\rangle+\varepsilon_2\varphi_\ell(j_1\varepsilon_1,j_2\varepsilon_2)))
-f(\langle\mu,\varphi\rangle)]
\frac{\mu(j_1\varepsilon_1,j_2\varepsilon_2)}{\varepsilon_1 \varepsilon_2}\nonumber\\
&+&\frac12 \sum_{j_1=0}^{\infty}\sum_{j_2\ge j_1}^\infty[f(\langle
\mu,\varphi\rangle-\varepsilon_2\varphi_\ell(j_1\varepsilon_1,j_2\varepsilon_2)))-f(\langle\mu,\varphi\rangle)]
\frac{\mu(j_1\varepsilon_1,j_2\varepsilon_2)}{\varepsilon_1\varepsilon_2}\nonumber\\
&+&\frac12(1-\varepsilon_1
c_\ell)\sum_{j_1=1}^{\infty}\sum_{j_2\ge j_1}^\infty
\\&&\cdot[f(\langle
\mu,\varphi\rangle-\varepsilon_1\varphi_\ell(j_1\varepsilon_1,j_2\varepsilon_2)
+\varepsilon_1\varphi_\ell((j_1+1)\varepsilon_1,(j_2+1)\varepsilon_2))
-f(\langle\mu,\varphi\rangle)]
\frac{j_1\mu(j_1\varepsilon_1,j_2\varepsilon_2)}{\varepsilon_1\varepsilon_2}\nonumber\\
&+&\frac12(1+\varepsilon_1
c_\ell)\sum_{j_1=1}^{\infty}\sum_{j_2\ge
j_1}^\infty\\&&\cdot[f(\langle
\mu,\varphi\rangle-\varepsilon_1\varphi_\ell(j_1\varepsilon_1,j_2\varepsilon_2)
+\varepsilon_1\varphi_\ell((j_1-1)\varepsilon_1,(j_2-1)\varepsilon_2)))
-f(\langle\mu,\varphi\rangle)]
\frac{j_1\mu(j_1\varepsilon_1,j_2\varepsilon_2)}{\varepsilon_1\varepsilon_2}\nonumber\\
&+&\frac12 \sum_{j_1=0}^{N^{\frac{2\ell-1}{4}}}\sum_{j_2\ge
j_1}^\infty\\&&\cdot[f(\langle
\mu,\varphi\rangle-\varepsilon_2\varphi_{\ell+1}(j_2\varepsilon_1)
+\varepsilon_2\varphi_{\ell+1}((j_2+1)\varepsilon_2))))
-f(\langle\mu,\varphi\rangle)]
\frac{j_2\mu(j_1\varepsilon_1,j_2\varepsilon_2)}{\varepsilon_1\varepsilon_2}\nonumber\\
&+&\frac12\sum_{j_1=0}^{N^{\frac{2\ell-1}{4}}}\sum_{j_2\ge
j_1}^\infty\\&&\cdot[f(\langle
\mu,\varphi\rangle-\varepsilon_2\varphi_{\ell+1}(j_2\varepsilon_2)
+\varepsilon_2\varphi_{\ell+1}((j_2-1)\varepsilon_2)))
-f(\langle\mu,\varphi\rangle)]
\frac{j_2\mu(j_1\varepsilon_1,j_2\varepsilon_2)}{\varepsilon_1\varepsilon_2}\nonumber\\
&+&\frac12\sum_{j_1>{N^{\frac{2\ell-1}{4}}}}^{\infty}\sum_{j_2\ge
j_1}^\infty\\&&\cdot[f(\langle
\mu,\varphi\rangle-\varepsilon_1\varphi_\ell(j_1\varepsilon_1,j_2\varepsilon_2)
+\varepsilon_1\varphi_\ell(j_1\varepsilon_1,(j_2+1)\varepsilon_2))
-f(\langle\mu,\varphi\rangle)]
\frac{j_2\mu(j_1\varepsilon_1,j_2\varepsilon_2)}{\varepsilon_1\varepsilon_2}\nonumber\\
&+&\frac12\sum_{j_1>{N^{\frac{2\ell-1}{4}}}}^{\infty}\sum_{j_2\ge
j_1}^\infty\\&&\cdot[f(\langle
\mu,\varphi_\ell\rangle-\varepsilon_1\varphi_\ell(j_1\varepsilon_1,j_2\varepsilon_2)
+\varepsilon_1\varphi_\ell(j_1\varepsilon_1,(j_2-1)\varepsilon_2)))
-f(\langle\mu,\varphi\rangle)]
\frac{j_2\mu(j_1\varepsilon_1,j_2\varepsilon_2)}{\varepsilon_1\varepsilon_2}.
\end{eqnarray*}
\normalsize

We give a brief explanation of these terms.

1. This term corresponds to the migration into $B^{(N)}_\ell$ of
individuals from families currently  minimally occupying
$B^{(N)}_\ell$. Asymptotically as $N\to \infty$, this is
\[ c_\ell f'(\langle \mu,\varphi\rangle)\int_{[0,\varepsilon_0^{1/4})}\int\big( \frac{\partial}{\partial
x}\varphi_\ell(x,y)\big)\, y\,\mu_{\ell+1}(dx,dy),\] whose limit
(by (\ref{two-1})) is
$$c_\ell f'(\langle \mu,\varphi_\ell \rangle)\int \big(\frac{\partial}{\partial
x}\varphi_\ell(x,y)|_{x=0}\big)\, y\,\mu^0(dy).$$

2. The second term is similar to the first except that here only
families having more than minimal mass in $B^{(N)}_\ell$ appear.
In the limit $N \to \infty$ we have asymptotically
\begin{eqnarray} &&c_\ell f'(\langle
\mu,\varphi\rangle)\int_{x>\varepsilon_0^{1/4}}\int\frac{\partial}{\partial
x}\varphi_\ell(x,y)\frac{y}{N^{1/2}}\,\mu_{\ell, N}(dx,dy)\\&&\leq
c_\ell|f'(\langle \mu,\varphi\rangle)|\sup|
\frac{\partial}{\partial
x}\varphi_\ell(x,y)|\nonumber\\&&\qquad\Big(\int\int_{y>
N^{\frac{1}{8}}}y\mu(dx,dy) + \frac{N^{1/8}}{N^{1/2}}
\int_{{N^{-1/4}}}^\infty\int_0^{N^{1/8}}
\mu_{\ell,N}(dx,dy)\Big)\nonumber \\&&\leq c_\ell|f'(\langle
\mu,\varphi\rangle)|\sup| \frac{\partial}{\partial
x}\varphi_\ell(x,y)|\nonumber\\&&\qquad\Big(\int\int_{y>N^{\frac{1}{8}}}y\mu(dx,dy)
+ N^{-1/8}\int_0^\infty\int_0^\infty x \mu_{\ell,N}(
dx,dy)\Big)\nonumber\\&&= o(1)+O(\frac{1}{N^{1/8}})
\end{eqnarray} since $\int
\int_{x>0}\int x\,\mu_{\ell,N}(dx,dy) = O(1).\nonumber$

 3. and 4.  These  two terms arise
from the family level branching. If $|\varphi_\ell(x,y)|\leq
\mathrm{const}\cdot x$, then asymptotically they yield
\begin{eqnarray} &&\frac{1}{2N^{1/2}}f''(\langle \mu,\varphi\rangle)\int
\int\varphi^2(x,y)\mu(dx,dy)\nonumber\\&&=\frac{1}{2N^{1/2}}f''(\langle
\mu,\varphi\rangle)\int\int \varphi_{\ell+1}^2(y)\mu(dx,dy)
\nonumber\\&&+f''(\langle \mu,\varphi\rangle)\int_0^{N^{-1/4}}\int
\varphi_{\ell+1}(y)\varphi_\ell(x,y)\mu(dx,dy) \nonumber\\&&
+f''(\langle \mu,\varphi\rangle)\int
\varphi_{\ell}^2(x,y)\mu_\ell(dx,dy)\nonumber\end{eqnarray}

5. and 6. This corresponds to the critical birth and death of
individuals in $B^{(N)}_{\ell}$. The limiting term is
\[ f'(\langle \mu,\varphi_\ell \rangle)\langle
\mu,G^{(2)}_{c_\ell}\varphi_\ell\rangle,\] where
$G_{c_\ell}^{(2)}$ denotes the application of the operator
$G_{c_\ell}$ to the $x$ variable.

7. and 8. These terms correspond to the birth and death of
individuals in $B^{(N)}_{\ell+1}$ in families that do not or
sparsely  occupy $B^{(N)}_\ell$. Asymptotically we obtain
\[ \frac{1}{N^{1/2}}f'(\langle \mu,\varphi\rangle)\langle
\mu, y\frac{\partial^2}{\partial y^2}\varphi_{\ell+1}(y)\rangle\]

9. and 10. These terms correspond to the birth and death of
individuals in $B^{(N)}_{\ell+1}$ in families which have at least a minimal
number of members in  $B^{(N)}_\ell$. Asymptotically we obtain
\[ f'(\langle \mu,\varphi\rangle)\int_{x>\varepsilon_0^{1/4}}
\frac{1}{N^{1/2}}\mu_\ell(dx,dy) y\frac{\partial^2}{\partial
y^2}\varphi_{\ell}(x,y)=O(\frac{1}{N^{1/4}})\]

Collecting the limiting terms  as $N\to\infty$, we obtain
$\mathfrak{ G}^{(N)}F(\mu) \to \mathfrak{G}_\ell F(\mu)$
 where
\begin{eqnarray}
\mathfrak{ G}_\ell F(\mu)&= f'(\langle \mu,\varphi\rangle)\langle
\mu_\ell, G_{c_\ell}^{(2)} \varphi_\ell\rangle
+\frac{1}{2}f''(\langle \mu,\varphi_\ell\rangle)\langle
\mu_\ell,\varphi_\ell^2\rangle\\& +f'(\langle
\mu,\varphi_\ell\rangle)c_\ell\int\frac{\partial\varphi_\ell}{\partial
x}(x,y)|_{x=0}y\mu^0(dy)\nonumber.
\end{eqnarray}

We conclude that for any limit point of the probability laws of
$\mathcal{Y}^{(N)}\mathbf{1}_{x>0}$,
\[ \tilde F(\eta(t))-\int_{t_0}^t\mathfrak{G}_\ell\tilde F(\eta(s))ds\] is a
 martingale where
   $\tilde F\in C(M^1((0,\infty)))$ is defined by
    $\tilde F(\mu)= f(\langle\mu,\varphi_\ell\rangle+\langle
\mu^0,\varphi_{\ell+1}\rangle)$ and $f,\,\varphi_\ell$ and $
\varphi_{\ell+1}$ satisfy the same conditions as above. But this
coincides with the martingale problem of the two level branching
diffusion with constant multitype immigration source with
immigration rate of type $y$ given by $y\mu^0(dy)$ and total
immigration rate of $\int_0^\infty y\mu^0(dy)$  (recall
Proposition \ref{immprocessconvergence}) which is well posed and
determines a $M^1(0,\infty)$-valued diffusion process. This
completes the proof of the proposition.

\end{proof}

\subsection {The limiting multiscale transition
function}\label{sec5.4} The main result of the  section is the
identification of the limiting multiscale structure, which serves
to complete the proof of Theorem \ref{commute}.

\begin{theorem}\label{immigration} Fix $j \ge 1$ and $\zeta_{j+1}^* > 0$.
  Conditioned on $\zeta^{(N,\theta)}_{j+1}(0)=\zeta_{j+1}^{*(N)}$, and provided that
  $\zeta_{j+1}^{*(N)} \to \zeta_{j+1}^*$ as $N\to \infty$, then the
normalized equilibrium masses $\{\zeta_{\ell}^{(N,\theta)}(0)\}_{\ell =j+1,j,\dots,1}$ in a sequence
  of nested blocks $\{B^{(N)}_\ell\}_{\ell=1,\dots,j+1}$
 converge
as $N \to \infty$ in distribution to the backward Markov chain
$\{\zeta_{\ell}^{}\}_{\ell =j+1,j,\dots,1}$ with transition kernel
\[ P(\zeta_{\ell}\in A|\zeta_{\ell+1}=a)=\Pi^{(2)}_\ell(a,A)\] where
$\Pi^{(2)}_\ell$ is as in Proposition \ref{entrancela} b), and
$\zeta_{j+1}^{}=\zeta_{j+1}^*$.
\end{theorem}
\begin{proof}
The principal step in proving this is given by Proposition
\ref{keyproposition}. Then to verify that the limiting equilibrium
distribution as $N\to\infty$ of $\zeta^{(N,\theta)}_\ell$
conditioned on $\zeta^{(N,\theta)}_{\ell+1}=a$ is given by
$\Pi_\ell^{(2)}(a,dx)$ we then follow the argument in Corollary
\ref{diffapproxc}(b). The convergence of
$\{\zeta_{\ell}^{}(0)\}_{\ell =j+1,j,\dots,0}$ to the Markov chain
then follows by recursion  and the continuity of the mappings $
a\to \Pi_\ell^{(2)}(a,dx)$.

\end{proof}

Combining  Theorem \ref{immigration}
with the spatial ergodic theorem (Theorem \ref{spergthm})   and Proposition \ref{entrancela} b), we see that the interchange of the limits
$N\to \infty$ and $j \to \infty$ is justified and completes the proof of Theorem \ref{commute}.

\subsection{The particle level picture}
\setcounter{equation}{0}

In the previous section the main result is obtained using the
convergence of the solutions of the appropriate sequence of
martingale problems. Some additional understanding of the limiting
process can be obtained by examining the particle picture in both
$B^{(N)}_{\ell+1}$ and $B^{(N)}_\ell$ which we now briefly sketch
in an informal way. First note that the change in the population
in $B^{(N)}_{\ell+1}$ due to movement of particles between
$B^{(N)}_{\ell+1}$ and its exterior in time scale $N^{\ell/2}$ has
expected value $o(\frac{1}{N})$ and therefore is negligible.

For the moment we fix $t_0<0$ and consider the contributions to
$\zeta^{(N,\theta)}_\ell(0)$ coming from immigrants arriving in
the interval $[t_0N^{\ell/2},0)$.  Next we recall that by
(\ref{familysizeconstraint}) all but $O(\frac{1}{K})$ of the
population in $B^{(N)}_{\ell+1}$ at time $0$ is contained in
families with sizes in $(0,KN^{(\ell+1)/2})$. Therefore we can
subdivide the population at time $0$ into $O(N^{(\ell+1)/2})$
independent subpopulations of size $O(N^{(\ell+1)/2})$  where if
necessary we group together smaller families  to form
subpopulations of size $O(N^{(\ell+1)/2})$. (This ensures
independent level two branching in $B^{(N)}_\ell$ for clusters
coming from distinct families in this subdivision.) Since the
families undergo critical branching these  are the descendants of
$O(N^{1/2})$ ancestral families at time $t_0N^{\ell/2}$ and in the
time interval $(t_0N^{\ell/2},0)$ these undergo binary branching
with $n_N(s)$ branches at time $sN^{\ell/2}$. Now consider a
single family in $B^{(N)}_{\ell+1}$ containing $O(N^{(\ell+1)/2})$
individuals. By the analogue of Lemma \ref{constant} the
normalized population size of this family is constant in the time
scale $N^{\ell/2}$.

We subdivide the interval $[N^{\ell/2}t_0,0)$ into
$M(\in\mathbb{N})$ subintervals of equal length $\frac{1}{M} t_0
N^{\ell/2}$. Let $t_{k}:=\frac{M-k}{M}t_0$ and consider an
interval $(t_{k+1}N^{\ell/2},t_{k}N^{\ell/2})$.
      Recall that
the rate of migration of individuals into $B^{(N)}_\ell$ from
$B^{(N)}_{\ell+1}$  is $c_{\ell}N^{-\ell/2}$.  Therefore in the
time interval  the number of individuals to immigrate into the
ball $B_{\ell}^{(N)}$ from a family of size $ mN^{(\ell+1)/2}$ in
$B_{\ell+1}^{(N)}$ is Poisson with mean
\begin{eqnarray}&&mN^{(\ell+1)/2}\times \frac{c_{\ell}}{N^{\ell/2}}%
((t_{k}-t_{k+1})N^{\ell/2})\times \frac{1}{N}\nonumber\\&& =
c_\ell m(t_{k}-t_{k+1})N^{(\ell-1)/2}.
\end{eqnarray}
 Now consider the question of which of
these have descendants alive in $B^{(N)}_\ell$ at time zero. Let
$U:= \min\{s:{n_N(s)=2}\}$. In order to determine this recall that
from the
 structure of the genealogy of the critical branching cluster
(see e.g. \cite{Du},\cite{F}) the  random variable $U$ is uniform
distributed on $[0,1]$. Moreover the probability that an initial
individual produces a non-empty set of descendants in a time
interval of length $U |t_k|N^{\ell/2}$ is asymptotically as
$N\to\infty$
\[
\frac{2c_{\ell}e^{-c_{\ell}Ut_{k}}}{N^{\ell/2}(1-e^{-c_{\ell}Ut_{k}})}.
\]
(see e.g. \cite{AN}(Chapt. 3)). Therefore the probability that
that any of the $O(N^{(\ell-1)/2})$ immigrants arriving from a
family of size $N^{(\ell+1)/2}$ in $[t_{k+1}N^{\ell/2}
,t_{k}N^{\ell/2})$ have descendants alive at time $0$ is no larger
than \[ \frac{2c_{\ell}e^{-c_{\ell}u_0(N)
t_{k}}}{N^{1/2}(1-e^{-c_{\ell}u_0(N)t_{k}})}+ P(U<u_0(N))\]
 for any
$0<u_0(N)<1$. Again by (\cite{Du},\cite{F})  the number of
(family) branches $n_N(t_k)$ at time $t_kN^{\ell/2}$ converges in
distribution as $N\to\infty$ for each of the $O(N^{1/2})$
ancestral families and the family trees from the different
ancestral families are independent. Choosing $u_0(N)\to 0$ such
that $N^{1/2}u_0(N)\to \infty$ we conclude that there is a Poisson
[with mean of order $O((t_k-t_{k+1})e^{-c_\ell t_k})$] number of individuals
who immigrate in $(N^{\ell/2}t_0,N^{\ell/2}t_k]$ producing
descendants at time $0$, and these all come from different
branching trees.

 Therefore the
population at time zero coming from immigrants arriving in
$(N^{\ell/2}t_0,N^{\ell/2}t_k)$ is asymptotically composed of
$O(1)$ two level families each originating from one individual and
these all come from different independent subpopulations  in
$B^{(N)}_{\ell+1}$. In particular the expected mass coming from
clusters containing two or more immigrants from the same
subpopulation in $B^{(N)}_{\ell+1}$ is $O(\frac{1}{N^{1/2}})$.
Each two level cluster coming from an immigrating individual
develops by two level branching, namely, family level branching
inherited from the family branching in $B^{(N)}_{\ell+1}$ and
subcritical individual level branching in $B^{(N)}_\ell$. Finally,
in both cases the total migration rate of individuals into
$B_{\ell
}^{(N)}$ from $B_{\ell+1}^{(N)}$ is%
\[
(\zeta_{\ell+1}^{(N,\theta)}(N^{\ell/2}t_0)N^{\ell+1})\times
(c_\ell N^{-\ell/2})\times N^{-1}
=c_\ell N^{\ell/2}\zeta_{\ell+1}^{(N,\theta)}(N^{\ell/2}t_0)%
\]
where on the left hand side the first factor is the number of
particles in $B^{(N)}_{\ell+1}$, the second is the individual
migration rate to a point chosen randomly in $B^{(N)}_{\ell+1}$
and the last factor is the probability that the tagged ball is
chosen. Therefore, asymptotically as $N\to\infty$, the population
at time zero consists of clusters of descendants of individuals
that immigrate into $B^{(N)}_\ell$ during the time interval
$(-\infty,0]$ at rate $c_\ell
N^{\ell/2}\zeta_{\ell+1}^{(N,\theta)}(N^{\ell/2}t_0)$ and
subsequently  undergo two level branching. In the limit
$N\to\infty$ these clusters correspond to the jumps of the
subordinator $S_\ell(\cdot)$ defined in section \ref{sec4.3}.

\begin{remark}
One can also gain some understanding of the convergence to
equilibrium from a spatially homogeneous initial population (more
general than that addressed in Proposition \ref{4+dequil}) with
intensity $\theta >0$.  Two ingredients are involved in the
convergence to an equilibrium with intensity $\theta$.  The first
is the strong transience condition on the random walk.  The other
feature is the structure of the local family sizes.  We see from
the above that the property that the contribution of families in
$B^{(N)}_{\ell+1}$ containing a number of individuals larger than
$ O(N^{(\ell+1)/2})$ to the equilibrium population in
$B^{(N)}_\ell$ is asymptotically negligible and this property is
then inherited by $B^{(N)}_\ell$. However if the initial family
sizes are too large this iteration can degenerate due to the
family level critical branching and the limiting population is
locally degenerate. For a more detailed analysis of this
phenomenon, see \cite{GH}.

\end{remark}

\section{Appendix}
\subsection{Size-biasing and Palm distributions}\label{subsec5}
\setcounter{equation}{0}

\begin{definition}\label{sizebiasing}
If $\pi$ is a measure on some measurable space $M$, and
$\mathfrak{s}$ is a nonnegative measurable function on $M$ with
$0< \bar s:=\int\mathfrak{s}(z)\pi(dz)<\infty$, then we call the
probability measure $\hat{\pi}$ given by \[
\hat{\pi}(dz):=\frac{\mathfrak{s}(z)}{\bar s}\pi(dz) \] the
size-biasing of $\pi$ with respect to $\mathfrak{s}$.  (Here, we
think of $\mathfrak{s}(z)$ as measuring the size of the object
$z$.)
\end{definition}

\begin{remark}\label{explainsb}
a) An  example of size-biasing which is important in our context
arises as follows. Let $\pi = \mathcal L (\eta)$ be the
distribution of a random measure $\eta$ (on some Polish space $E$,
say). Denote the intensity measure of $\eta$ by $\lambda$, and fix
a nonnegative measurable function $f$ on $E$ with $0 < \langle
\lambda, f \rangle < \infty$. Define the size of a measure $m$ on
$E$ by $\mathfrak{s}(m) := \int f(x) m(dx) = \langle m, f \rangle
$, and denote by $\pi_{f}$ the size-biasing of $\pi (dm)$ with
$\langle m, f \rangle$.

b) Assuming that $\lambda$ is locally finite,  choosing $f =
1_{B}$ where $B$ is a ball in $E$, and letting $B$ shrink gives
the family of Palm distributions $\pi_{x}, x \in S$. Formally,
these arise as the disintegration of the measure $\pi(dm)m(dx)$
with respect to its second marginal $\mathbb E \eta$, that is
\begin{equation}\label{Campbell}
           \mathbb E G(\eta)\langle \eta, h \rangle = \int h(x) G(m)
           \pi_{x}(dm) \lambda(dx).
\end{equation}
See \cite{K}, chapter 10, for more background; there, a random
measure whose distribution is the size-biasing of $(\mathcal
L(\eta))(dm)$ with $\langle m, f \rangle$  is denoted by
$\eta_{f}$.

c) The following fact (\cite {K}, formula (10.6)) is immediate
from (\ref{Campbell}):

The size-biasing $\pi_{f}$ of $\pi(dm)$ with $\langle m, f
\rangle$ is
\begin{equation}\label{sizebiaseta}
           \frac 1{\langle \lambda, f\rangle} \int f(x) \pi_{x}(.)
           \lambda(dx) = \mathbb E \pi_{\hat X},
           \end{equation}
           where $\hat X$ is a random element in $E$ whose distribution is
           the size-biasing of $\lambda$ with $f$.

d) If in the just described situation $E$ consists of one element
only, then the finite random measures on $E$ are
$\mathbb R_{+}$-valued random variables. When speaking of the
size-biasing of a measure $\pi$ on $\mathbb R_{+}$ without
specifying the size function, we always mean the size-biasing of
$\pi(dx)$ with $x$.
\end{remark}

Let us write $\Pi_{\lambda}$ for the distribution of a Poisson
random counting measure on $E$ with intensity measure $\lambda$.
It is well known (see e.g. \cite {K}, beginning of chapter 11)
that the Palm distributions of $\Pi_{\lambda}$ arise as the
distributions of $\Phi + \delta_{x}$, $x \in E$, where
$\mathcal{L}(\Phi) = \Pi_{\lambda}$.

Now let $\sigma$ be a probability measure on $\mathbb R_{+}$ with
$m_{\sigma} := \int \tau \sigma (d\tau) \in (0, \infty)$, and
write
$$\Pi_{\sigma,\lambda} := \int \Pi_{\tau \lambda} \sigma(d\tau)$$ for the
{\em mixed Poisson distribution} with mixing measure $\sigma$.

The following lemma, whose proof we include for convenience, is
part of a characterization theorem (\cite{K}, Theorem 11.5) of
mixed Poisson processes.
\begin{lemma}\label{MPPalm}
           The Palm distributions of $\Pi_{\sigma, \lambda}$ arise as the
distributions of $\Phi+\delta_{x}$, $x \in E$, where $\mathcal
L(\Phi)=\Pi_{\hat \sigma, \lambda}$, and $\hat \sigma$ is the
size-biasing of $\sigma$.
\end{lemma}
\begin{proof} Using the above mentioned form of the Palm
distributions of $\Pi_{\tau\lambda}$, we obtain for all
nonnegative measurable $h$ and $G$ defined on $E$ and $M_{c}(E)$, the space of locally finite counting measures on $E$,
respectively:
\begin{eqnarray}\label{CMP}
           \int G(\psi)\langle \psi,h\rangle \Pi_{\sigma, \lambda}(d\psi)
&=& \int\int
           G(\psi) \langle \psi,h\rangle \Pi_{\tau \lambda}(d\psi)
           \sigma(d\tau) \\ \nonumber
           &=& \int \int \tau \lambda(dx) h(x) G(\psi+\delta_{x}) \Pi_{\tau
\lambda}(d\psi)
           \sigma(d\tau) \\ \nonumber &=& \int m_{\sigma} \lambda(dx) h(x) \int
           G(\psi+\delta_{x}) \Pi_{\tau \lambda}(d\psi) \frac 1{m_{\sigma}}
           \tau \sigma(d\tau)
           \end{eqnarray}
           \end{proof}
\begin{corollary}\label{sbMP}
          Assume $0 < \langle \lambda, h \rangle < \infty$.
Then the size-biasing of $\Pi_{\sigma, \lambda}(d\psi)$  with
$\langle \psi, h \rangle$ arises as the distribution of $\Phi +
\delta_{\hat X}$, where $\mathcal{L}(\Phi) = \int \Pi_{\tau
\lambda} \hat \sigma(d \tau)$, $\mathcal{L}(\hat X)$ is the
size-biasing of $\lambda$ with $h$, and $\Phi$ and $\hat X$ are
independent.
\end{corollary}
\begin{proof}
           This can be seen either by combining Remark \ref{explainsb} c) and
           Lemma \ref{MPPalm}, or by dividing (\ref{CMP}) through $\mathbb
           E_{\sigma, \lambda}\langle \Psi, h\rangle = m_{\sigma}\langle
\lambda, h
           \rangle$.
           \end{proof}

\subsection{Subcritical Feller branching}\label{SFB}
\setcounter{equation}{0}

Let us fix  $c > 0$. In the following, $X$ will denote a
{\em $c$-subcritical Feller branching diffusion ($c$-FBD) process}. In
other words, $X$ is an $[0,\infty)$-valued diffusion process
satisfying
\begin{equation}\label{FBdyn}
                     dX_{t} = \sqrt {X_{t}}\,dW_{t} -c X_{t} dt
\end{equation}
where $W_t$ is Brownian motion.      For $\varepsilon > 0$, let $X^{\varepsilon}$ be the
$c$-FBD process starting in $\varepsilon$ at time $0$.

       From well-known results on Galton-Watson processes conditioned to
survival \cite{AN, Ge} and cluster sizes in continuous-state
branching processes \cite{D}, one expects that the conditional law
$\mathcal L(X^{\varepsilon}_{t}| X^{\varepsilon}_{t}>0)$ converges
to an exponential distribution as $\varepsilon \to 0$. The
following calculation verifies this and identifies the parameter.

The Laplace functional of $X^{\varepsilon}_{t}$ is given by
\begin{equation}\label{laplaceX}
\mathbb E(e^{-\lambda X^{\varepsilon}_t}) =e^{-\varepsilon
v(t,\lambda)},\quad \lambda \ge 0,
\end{equation}
where $v=v(t,\lambda)$ is the solution of
\begin{equation}\label{evolvv}
\frac{\partial v(t,\lambda)}{\partial t}
=-cv(t,\lambda)-\frac{1}{2} v^{2}(t,\lambda),\quad  v(0,\lambda)
=\lambda.
\end{equation}
The solution of (\ref{evolvv}) is given by
\begin{equation}\label{formofv}
v(t,\lambda)=\frac{2\lambda ce^{-ct}}{\lambda(1-e^{-ct})+2c}.
\end{equation}
Combining (\ref{laplaceX}) and (\ref{formofv}) one obtains by a
straightforward calculation
\begin{equation}\label{probnonextinct}
          \mathbb
P(X^{\varepsilon}_t\neq 0)=1-\exp\left({\frac{-2\varepsilon
ce^{-ct}}{1-e^{-ct}}}\right),
\end{equation}
\begin{equation}\label{expectofX}
\mathbb E X^\varepsilon_t = \varepsilon e^{-ct}, \quad \varepsilon
> 0
\end{equation}
and
\begin{equation}\label{cond Laplace}
          \lim_{\varepsilon\rightarrow0}\mathbb E(e^{-\lambda
X^{\varepsilon}_t}|X^{\varepsilon}_t\neq0)=\frac{2c}
{\lambda(1-e^{-ct})+2c}.
\end{equation}
Writing $\mathrm {Exp}(u)$ for the exponential distribution with
parameter $u$, we obtain immediately:

\begin{lemma}\label{aboutX} Fix $t > 0$. \\ a)
                     \begin{equation}\label{convofXepsilon}
                   \mathcal L (X^{\varepsilon}_{t} | X^{\varepsilon}_t > 0) \Rightarrow
\mathrm {Exp}\left(\frac{2c}{1-e^{-ct}}\right) \quad
\mbox{  as } \varepsilon  \to 0.
                     \end{equation}
                     b) \begin{equation}\label{survprob}
\varepsilon^{-1} \mathbb P [X^{\varepsilon}_{t}>0] \to
\frac{2ce^{-ct}}{1-e^{-ct}}\quad \mbox{  as } \varepsilon  \to 0.
\end{equation}
\end{lemma}

For $\varepsilon > 0$, let
\begin{equation}\label{Xbarepsilon}
                     \bar X^{\varepsilon} =
\sum_{i=1}^{N_{\varepsilon}}X^{\varepsilon,i},
\end{equation}
where $N_{\varepsilon}$ is a Poisson(${\varepsilon}^{-1}$)-random
variable and $X^{\varepsilon,1},X^{\varepsilon,2},\ldots$ are
independent copies of $X^{\varepsilon}$. The following lemma is an
easy consequence of Lemma \ref{aboutX}.
\begin{lemma}\label{aboutXbar}
                     Fix $t >0$. Then
                 \begin{equation}\label{convtozetaXbar}
                     \bar X^{\varepsilon}_{t} \Rightarrow  \bar X^0_{t} \mbox{ as }
                     \varepsilon \to 0,
                  \end{equation}
where $\bar X^0_{t}$ is infinitely divisible with canonical
measure
\begin{equation}\label{cm}
                     \kappa_{t}:=  \frac{2ce^{-ct}}{1-e^{-ct}}{\rm
                     Exp}\left(\frac{2c}{1-e^{-ct}}\right).
                  \end{equation}
\end{lemma}

Let us now define
the measure
\begin{equation}\label{gammac}
                     \gamma_{c}:= c\int_{0}^\infty \kappa_{s}
                     \,ds.
\end{equation}
An
elementary calculation based on (\ref{cm}) and (\ref{gammac})
shows that
\begin{equation}\label{gammaexplicit}
                    \gamma_{c}(dx) = 2c\frac 1x e^{-2cx} dx, \quad x > 0.
\end{equation}
This identifies $\gamma_{c}$ as the canonical measure of the ${\rm
Gamma}(2c,2c)$-distribution, and goes along with the well-known fact that
      the equilibrium distribution of (\ref{FBDa}) is the ${\rm
Gamma}(2ca,2c)$-distribution. We note
that $\gamma_{c}$ is the L\'evy measure of the Gamma subordinator
$S(\tau)$, $\tau \ge 0$, with scale parameter $2c$ and
$\mathbb E S(1) = 1$. The following is obvious from
(\ref{gammaexplicit}):
\begin{remark}\label{momentsofgamma}
                   a) $\int_{0}^\infty x \gamma_c(dx) = 1.$\\
                     b) $\int_{0}^\infty x^2 \gamma_c(dx) = 1/2c.$
\end{remark}

Finally we consider the semigroup  $(T_{t})$ be the semigroup of
the c-FBD process. Recalling (\ref{convofXepsilon}), (\ref{survprob}) and (\ref{cm}),
the $(T_{t})$-entrance law $(\kappa_{t})$ from $0$
is  given by
\begin{equation}\label{kappa}
                     \kappa_{t}= \lim_{\varepsilon \to 0} \frac 1\varepsilon
\delta_{\varepsilon}T_{t}  = \lim_{\varepsilon \to 0} \frac
1\varepsilon \mathcal L(X^{\varepsilon}_{t};
X^{\varepsilon}_{t}\neq 0 ),  \quad t > 0,
\end{equation}
where $(X^{\varepsilon}_{t})$ is the $c$-FBD-process starting in
$\varepsilon$ at time $0$. Because of (\ref{cm}), $\kappa_{t}(dx)$
has density
\begin{equation}\label{density}
\kappa_{t}(x)=
\frac{(2c)^2e^{-ct}}{(1-e^{-ct})^2}\exp\left(-\frac{2cx}{1-e^{-ct}}\right),
\quad x \in (0,\infty).
\end{equation}
Then for $f\in C_1((0,\infty))$ and $t>0$,
\begin{equation}\label{dTz}
\frac{d}{dx}T_tf(x)\Big|_{x=0}=\lim_{\varepsilon \to 0}
\frac{1}{\varepsilon} T_{t}f(\varepsilon)=\int f(y)
\kappa_{t}(dy).
\end{equation}

\bigskip

{\bf Acknowledgements} The authors thank the hospitality of The
Fields Institute (Toronto, Canada), Carleton University (Ottawa,
Canada), the Center for Mathematical Research (CIMAT, Guanajuato,
Mexico), Erwin Schr\"odinger Institute (Vienna, Austria) and the
Johann Wolfgang Goethe University (Frankfurt, Germany), where
mutual working visits took place. L.G.G. also thanks the Institute
of Mathematics, National University of Mexico (UNAM), where he
spent a sabbatical during 2002.

\end{document}